\documentclass[a4paper, oneside]{amsart}

\usepackage{amsthm, amsmath, amssymb, amstext, amsfonts, eucal}
\usepackage[top=2cm, bottom=2cm]{geometry}
\usepackage{hyperref}
\usepackage{enumitem}
\usepackage[svgnames]{xcolor}
\usepackage{accents}
\usepackage{booktabs}
\usepackage{multirow}
\usepackage{pdflscape}
\usepackage[capitalize, noabbrev]{cleveref}
\usepackage{caption}
\usepackage{subfiles}

\makeatletter
\renewcommand{\paragraph}{%
  \@startsection{paragraph}{4}%
  \z@{.5\linespacing\@plus.7\linespacing}{-.5em}%
  {\normalfont\bfseries}%
}
\makeatother

\usepackage{tabu}
\usepackage{array, makecell, graphicx}

\definecolor{smoked}{RGB}{216, 212, 204}
\definecolor{mauve}{RGB}{200, 55, 171}
\definecolor{apricot}{RGB}{250, 144, 4}
\definecolor{sky}{RGB}{66, 169, 244}
\definecolor{plum}{RGB}{76, 0, 102}

\definecolor{lightmauve}{RGB}{232, 173, 220}
\definecolor{lightapricot}{RGB}{253, 211, 155}
\definecolor{lightsky}{RGB}{178, 221, 251}
\definecolor{lightplum}{RGB}{184, 153, 192}

\definecolor{darksmoked}{RGB}{198, 194, 176}

\hypersetup{
   colorlinks,
    linkcolor={mauve},
    citecolor={plum},
    urlcolor={plum}
}
\usepackage{graphicx}
\usepackage{tikz,tqft}
\usetikzlibrary{matrix,arrows,bending,cd,decorations.markings,calc,patterns,angles,quotes,backgrounds,math}
\tikzset{>=latex}

\theoremstyle{definition}
\newtheorem{defn}{Definition}[section]
\newtheorem{rem}[defn]{Remark}
\newtheorem*{notation}{Notation}
\newtheorem*{terminology}{Terminology}

\theoremstyle{plain}
\newtheorem{lem}[defn]{Lemma}
\newtheorem{prop}[defn]{Proposition}
\newtheorem{thm}[defn]{Theorem}
 \newtheorem{thma}{Theorem}
 
\newtheorem{cor}[defn]{Corollary}

\newtheorem{fact}[defn]{Fact}

\theoremstyle{remark}
\newtheorem{claim}[defn]{Claim}

\newenvironment{proofclaim}
 {\proof}
 {\endproof}

\numberwithin{equation}{section}

\DeclareMathOperator{\id}{id}
\let\Re\relax
\DeclareMathOperator{\Re}{Re}
\let\Im\relax
\DeclareMathOperator{\Im}{Im}
\DeclareMathOperator{\res}{res}
\DeclareMathOperator{\diag}{diag}
\newcommand{\reg}{\mathrm{reg}}

\newcommand{\Rnum}[1]{\uppercase\expandafter{\romannumeral #1\relax}}
\renewcommand{\leq}{\leqslant}
\renewcommand{\geq}{\geqslant}

\DeclareMathOperator{\stab}{stab}

\DeclareMathOperator{\R}{\mathbb{R}}
\DeclareMathOperator{\C}{\mathbb{C}}
\DeclareMathOperator{\Z}{\mathbb{Z}}

\DeclareMathOperator{\CP}{\mathbb{CP}}

\DeclareMathOperator{\CH}{\mathbb{CH}}

\newcommand{\PU}{\operatorname{PU}}
\newcommand{\SU}{\operatorname{SU}}
\newcommand{\U}{\operatorname{U}}

\DeclareMathOperator{\Hom}{Hom}

\DeclareMathOperator{\Mod}{Mod}
\DeclareMathOperator{\PMod}{PMod}
\DeclareMathOperator{\Aut}{Aut}
\DeclareMathOperator{\Out}{Out}
\DeclareMathOperator{\FixedPts}{\mathcal{F}}

\title{Pure mapping class group fixed points in $\PU(2,1)$}

\author{Samuel Bronstein}
\address[S.~Bronstein]{Max Planck Institute for Mathematics In the Sciences, Leipzig, Germany}
\email{bronstein@mis.mpg.de}

\author{Thomas Le Fils}
\address[T. Le Fils]{IRMAR, Universit\'e de Rennes, Rennes, France}
\email{thomas.le-fils@cnrs.fr}

\author{Arnaud Maret}
\address[A.~Maret]{Université de Neuchâtel, Neuchâtel, Suisse}
\email{arnaud.maret@unine.ch}

\date{\today}

\begin{document}
\begin{abstract}
We classify all the fixed points for the action of the pure mapping class group on $\PU(2,1)$ character varieties of surface groups. Up to trivial peripheral monodromies, we prove that the irreducible fixed points are either pair of pants representations, or restrictions of reduced Burau representations to a once-punctured torus and restrictions of reduced Gassner representations to a $4$-punctured sphere. On the other hand, we prove that the only non-trivial reducible fixed points appear in genus~$0$ where they map all but at most three peripheral loops in the centralizer of the image, and in genus~$1$ where they are valued in the quaternion group~$H_8$.

As consequences, we deduce (i) that an irreducible representation in genus~$0$ with at least four non-trivial peripheral monodromies is rigid if and only if it is a fixed point, and (ii) a partial classification of irreducible representations of pure braid groups into $\PU(2,1)$.
\end{abstract}
\maketitle

\section{Introduction}
\subsection{Main results}
Among the different axes of studies around mapping class group dynamics on character varieties of surface groups representations, we focus on investigating finite orbits, and more precisely fixed points. While we now have a complete panorama of mapping class group fixed points in rank $2$ (that is, for representations into $\mathrm{SL}_2\C$)~\cite{lisovyy-tykhyy,LLL,bronstein-maret}, a systematic classification in higher rank remains largely open. As a first step in that direction, this paper classifies all mapping class group fixed points for representations into~$\PU(2,1)$. Throughout the paper, $S_{g,n}$ will denote an oriented and connected surface of genus~$g\geq 0$ and with~$n\geq 0$ punctures. Its \emph{pure mapping class group}, denoted by~$\PMod(S_{g,n})$, consists of all isotopy classes of orientation-preserving homeomorphisms of $S_{g,n}$ that fix each puncture individually.

\subsubsection{Irreducible fixed points}
We say that a representation $\rho\colon\pi_1S_{g,n}\to\PU(2,1)$ is an~\emph{irreducible fixed point} if the image of $\rho$ is irreducible (Definition~\ref{defn:irreducible}) and if $\PMod(S_{g,n})$ fixes the conjugacy class of $\rho$. Moreover, we say that $\rho$ is \emph{corked} if all peripheral loops (i.e.~loops freely homotopic to a puncture of $S_{g,n}$) have non-trivial images (Definition~\ref{def:corked}).
Our main result (Theorem~\ref{thm:intro-classification-reduced-irreducible}) shows that the only instances of corked irreducible fixed points are restrictions of braid group and pure braid group representations (Sections~\ref{sec:braid-group-representations}--\ref{sec:pure-braid-group-representation} and~\ref{sec:Burau-fixed-points-torus}).
More precisely, we prove that the $1$-parameter family of reduced Burau representations of~$\PMod(S_{1,2})\cong B_4/Z(B_4)$ and the $4$-parameter family of reduced Gassner representations of~$\PMod(S_{0,5})\cong PB_4/Z(B_4)$ induce all non-elementary corked irreducible fixed points by restrictions to surface groups along the point-pushing map from the Birman exact sequence~\eqref{eq:Birman-short-exact-sequence}.

\begin{thma}[Lemma~\ref{lem:irreducible-fixed-point-once-punctured-torus} and Corollary~\ref{cor:irreducible_S04}]\label{thm:intro-classification-reduced-irreducible}
    A corked representation~$\rho\colon\pi_1S_{g,n}\to\PU(2,1)$ is an irreducible fixed point if and only if one of the following holds:
    \begin{itemize}
        \item $(g,n)=(0,3)$ and $\rho$ is any irreducible representation;
        \item $(g,n)=(1,1)$ and $\rho$ is a \emph{Burau fixed point}, i.e.~it is conjugate to the restriction of a reduced Burau representation $B_4/Z(B_4)\to\PU(2,1)$ along the point-pushing map $\pi_1S_{1,1}\hookrightarrow B_4/Z(B_4)$;
        \item $(g,n)=(0,4)$ and $\rho$ is a \emph{Gassner fixed point}, i.e.~it is conjugate to the restriction of a reduced Gassner representation $PB_4/Z(B_4)\to\PU(2,1)$ along the point-pushing map $\pi_1S_{0,4}\hookrightarrow PB_4/Z(B_4)$.
    \end{itemize}
\end{thma}

The case $(g,n)=(0,3)$ is immediate since $\PMod(S_{0,3})$ is trivial. The classification of irreducible fixed points in the case $(g,n)=(1,1)$ as restrictions of reduced Burau representations is a consequence of Formanek's classification of braid group representations~\cite{formanek} and Lemma~\ref{lem:regular-and-abelian-fixed-points-are-restrictions}. We give an alternative proof of it in Lemma~\ref{lem:irreducible-fixed-point-once-punctured-torus} using the hyperelliptic involution of $S_{1,1}$ to reduce to the 4-punctured sphere case, previously treated in Corollary~\ref{cor:irreducible_S04}. Falbel--Garden--Guilloux have also independently obtained this classification for the once-punctured torus by a different argument in a recent work~\cite[Proposition~3.1]{FGG26}. The two main new contributions of Theorem~\ref{thm:intro-classification-reduced-irreducible} are (i) the classification in the case $(g,n)=(0,4)$, and (ii) the proof that corked irreducible fixed points do not exist for other topologies. 

If we allow peripheral loops in the kernel of the representation, then all irreducible fixed points arise from one of the examples in Theorem~\ref{thm:intro-classification-reduced-irreducible} by adding extra trivial punctures. This is the content of the following result.

\begin{thma}[Theorem~\ref{thm:irreducible-fixed-points-genus-zero} and Theorem~\ref{thm:fixed-points-genus-at-least-1}]\label{thm:intro-classification-unreduced-irreducible}
    A representation~$\rho\colon\pi_1S_{g,n}\to\PU(2,1)$ is an irreducible fixed point if and only if one of the following holds:
    \begin{itemize}
        \item $g=0$ and $\rho$ maps the peripheral loops around exactly $n-4$ punctures to the identity and, after filling the trivial punctures, is a Gassner fixed point;
        \item $g=0$ and $\rho$ maps the peripheral loops around exactly $n-3$ punctures to the identity and, after filling the trivial punctures, is an irreducible pair of pants representation;
        \item $g=1$ and $\rho$ maps the peripheral loops around exactly $n-1$ punctures to the identity and, after filling the trivial punctures, is a Burau fixed point.
    \end{itemize}
    In particular, if $g\geq 2$, then there are no irreducible fixed points.
\end{thma}

\subsubsection{Reducible fixed points}
Representations $\rho\colon\pi_1S_{g,n}\to\PU(2,1)$ with reducible image (Definition~\ref{defn:irreducible}) whose conjugacy class is fixed by $\PMod(S_{g,n})$ are called \emph{reducible fixed points}. When $g=0$, examples of reducible fixed points include all reducible representations $\rho\colon\pi_1S_{0,n}\to\PU(2,1)$ that map the peripheral loops around at least $n-3$ punctures into the centralizer of $\rho(\pi_1S_{0,n})$ in $\PU(2,1)$. We call such representations \emph{pants-like} (Definition~\ref{defn:pants-like}). Note that pants-like representations include all abelian representations and all representations of a pair of pants with extra punctures mapped to the identity (the latter form of pants-like representations may be irreducible). We prove that there are no other reducible fixed points.

\begin{thma}[Theorem~\ref{thm:reducible-fixed-points-are-alomst-abelian}]
    A reducible representation~$\rho\colon\pi_1S_{0,n}\to\PU(2,1)$ is a fixed point if and only if $\rho$ is pants-like.
\end{thma}

When $g=1$, there are examples of reducible fixed points with image the quaternion group $H_8$ defined in~\eqref{eq:quaternion-group} which we call \emph{(extended) quaternionic} fixed points (Definitions~\ref{def:quaternionic} and~\ref{def:extended-quaternionic-fixed-point}). These turn out to be the only non-trivial examples of fixed points when $g\geq 1$.

\begin{thma}[Theorem~\ref{thm:fixed-points-genus-at-least-1}]
    Assume~$g\geq 1$. A representation~$\rho\colon\pi_1S_{g,n}\to\PU(2,1)$ is a reducible fixed point if and only if either $\rho$ is the trivial representation or $g=1$, $n\geq 1$, and $\rho$ is an extended quaternionic fixed point.
\end{thma}

\subsection{Consequences}
We point out two immediate consequences of our classification.

\subsubsection{Rigid representations}
The \emph{peripheral monodromies} of a representation $\rho\colon\pi_1S_{g,n}\to G$ are the $n$ conjugacy classes of the images of counter-clockwise oriented peripheral curves around the $n$ punctures of $S_{g,n}$. The space of all conjugacy classes of representations with the same peripheral monodromies is a \emph{relative character variety}. The representation $\rho$ is said to be \emph{rigid} if its conjugacy class is an isolated point of the relative character variety that contains it. When $g=0$ and $G=\PU(2,1)$, relative character varieties may a priori contain several isolated points (although we have not encountered an example yet). On the other hand, when $G=\mathrm{U}(2,1)$, then Katz's result~\cite[Theorem~1.1.2]{katz} and Claim~\ref{claim:conjugate-in-U(2,1)-iff-conjugate-in-GL_3C} imply that an irreducible representation $\rho\colon\pi_1S_{0,n}\to \U(2,1)$ is rigid if and only if the relative character variety that contains its conjugacy class is a singleton. Therefore, since an irreducible rigid representation $\rho\colon\pi_1S_{0,n}\to \PU(2,1)$ can always be lifted to an irreducible rigid representation $\widehat{\rho}\colon\pi_1S_{0,n}\to \U(2,1)$, the conjugacy class of $\rho$ is always fixed by $\PMod(S_{0,n})$. In other words, irreducible rigid representations into $\PU(2,1)$ are always fixed points. On the other hand, all the corked irreducible fixed points $\rho\colon\pi_1S_{0,4}\to\PU(2,1)$ are restrictions of Gassner representations by Theorem~\ref{thm:intro-classification-reduced-irreducible} which can be obtained by a single application of Katz's middle convolution (introduced in~\cite{katz}) to a non-trivial rank-$1$ local system (see~\cite[Remark~5.2]{dettweiler-reiter}) and are therefore rigid. We thus have the following characterization.
\begin{thma}\label{thm:intro-fixed-point-iff-rigid-spheres}
    Let $n\geq 4$ and $\rho\colon\pi_1S_{0,n}\to\PU(2,1)$ be an irreducible representation with at most $n-4$ trivial peripheral monodromies. Then $\rho$ is rigid if and only if $\rho$ is a fixed point. In that case, $\rho$ has exactly $n-4$ trivial peripheral monodromies.
\end{thma}
Note that the hypothesis of having at most $n-4$ trivial peripheral monodromies in Theorem~\ref{thm:intro-fixed-point-iff-rigid-spheres} is required as irreducible representations $\pi_1S_{0,3}\to\PU(2,1)$ with $n-3$ extra punctures mapped to the identity are always fixed points, but generically not rigid. Also note that there are examples of rigid irreducible representations $\pi_1S_{0,3}\to\PU(2,1)$ (with two regular semisimple and one pseudo-reflection peripheral monodromies). 

The case $g=1$ is worth a remark as well. By Theorems~\ref{thm:intro-classification-reduced-irreducible} and Theorem~\ref{thm:intro-classification-unreduced-irreducible}, any irreducible fixed point $\rho\colon\pi_1S_{1,n}\to\PU(2,1)$ is the restriction of a Burau representation with exactly $n-1$ extra punctures mapped to the identity. The only non-trivial peripheral monodromy of $\rho$ is non-scalar (Lemma~\ref{lem:properties-Burau-fixed-points-torus}) and therefore $\rho$ is not rigid by~\cite[Corollary~1.3.1]{katz}.

\subsubsection{Pure braid group representations}
Formanek produced a classification of all irreducible representations of the braid group $B_n$ into $\mathrm{PGL}_r\C$ when $r\leq n-1$. For instance, when $r=n-1$ and $n\geq 7$, only the reduced Burau representations are irreducible~\cite[Theorem~22]{formanek}. When $n=4$ and $r=3$, there are two kinds of irreducible representations $B_4\to \mathrm{PGL}_3\C$ other than the reduced Burau representations~\cite[Theorem~13]{formanek}\footnote{The first one factorizes through the special homomorphism $B_4\to B_3$ and its restrictions along the point-pushing maps $\pi_1S_{0,4}\hookrightarrow PB_4/Z(B_4)\subset B_4/Z(B_4)$ and $\pi_1S_{1,1}\hookrightarrow B_4/Z(B_4)$ are trivial. The second one is a one-parameter family whose restriction to $\pi_1S_{0,4}$ is trivial and whose restriction to $\pi_1S_{1,1}$ gives a non-trivial reducible fixed point in the $\mathrm{PGL}_3\C$ character variety (by Lemma~\ref{lem:restriction-of-PMod-representations-are-fixed-points}), but that is not a fixed point in the $\PU(2,1)$ character variety. This is explained by the fact that this special irreducible representation $B_4\to\mathrm{PGL}_3\C$ does not preserve a Hermitian form of signature $(2,1)$, even though its restriction to $\pi_1S_{1,1}$ does. Those representations of $B_4$ thus do not contribute to new mapping class group fixed points in $\PU(2,1)$ character varieties.} (see also~\cite[Theorem~4.3]{FGG26}).

On the contrary, the theory of representations of pure braid groups remains poorly understood and a classification of their representations out of reach. A consequence of our classification of mapping class group fixed points in $\PU(2,1)$ character varieties, in combination with Lemma~\ref{lem:regular-and-abelian-fixed-points-are-restrictions}, is the following.
\begin{thma}
    The representations $PB_4\to \PU(2,1)$ whose restrictions to $\pi_1S_{0,4}$ along the point-pushing map are irreducible and corked are exactly the reduced Gassner representations.
\end{thma}

We also prove that every fixed point $\pi_1S_{0,n}\to\PU(2,1)$ (irreducible or not) is the restriction of a representation $PB_n/Z(B_n)\to\PU(2,1)$ along the point-pushing map (Theorem~\ref{thm:all-fixed-points-lift-PU(2,1)}). Using that a subgroup of $\PU(2,1)$ is either abelian or has abelian centralizer (Corollary~\ref{cor:non-abelian-intersection-centralizers}), Lemmas~\ref{lem:regular-and-abelian-fixed-points-are-restrictions} and~\ref{lem:fibers-restriction-map-over-abelian-centralizers} provide a complete description of the fibres of representations $PB_n/Z(B_n)\to\PU(2,1)$ that restrict to the same fixed point $\pi_1S_{0,n}\to\PU(2,1)$. So, classifying all representations $PB_n/Z(B_n)\to\PU(2,1)$ comes down to describing explicitly all these fibres.

\subsection{Motivations}
Our study began when we encountered what we call the \emph{Casella fixed point} (Definition~\ref{defn:Casella-fixed-point}). It is a corked irreducible fixed point $\rho\colon\pi_1S_{0,4}\to\PU(2,1)$ with discrete image contained in $\PU(2,1;\Z[\omega])$, where $\omega$ is a primitive third root of unity. The existence of corked irreducible fixed points for $4$-punctured spheres contrasts with the case of rank 2 where all fixed points are either abelian or have peripheral curve in their kernel~\cite[Lemma~3.3]{bronstein-maret} (see also Remark~\ref{rem:pants-like-rank-2}). We then observed that the Casella fixed point naturally fits into the larger family of corked irreducible fixed points given by reduced Gassner representations of the pure braid group $PB_4$. This was our original motivation to verify that all irreducible fixed points arise in this way and lead to Theorems~\ref{thm:intro-classification-reduced-irreducible} and~\ref{thm:intro-classification-unreduced-irreducible}.

\subsubsection{Quantum representations}
Another source of examples of mapping class group fixed points in character varieties of surface group representations into projective pseudo-unitary groups are obtained by restrictions of \emph{quantum representations} of (pure) mapping class groups along the point-pushing map. Koberda--Santharoubane showed that this gives examples of fixed points with infinite image for almost every surface~\cite[Theorem~1.1]{koberda-santharoubane}. See also Santharoubane's lecture notes for a broader introduction~\cite{santharoubane-notes}. More recently, Deroin--Marché proved that the so-called Fibonacci quantum representations $\Mod(S_{0,5})\to\PU(2,1)$ and $\Mod(S_{1,2})\to\PU(2,1)$, which induce mapping class group fixed points for $S_{0,4}$ and $S_{1,1}$ by restriction, are holonomies of complex hyperbolic structures on some compactifications of the corresponding moduli spaces of curves~\cite{deroin-marche}.

\subsubsection{Isomonodromy}
Finite mapping class group orbits, in particular fixed points, are in classical correspondence with algebraic solutions to systems of isomonodromic differential equations, such as Painlevé~\Rnum{6} and Garnier systems. A precise statement of the correspondence was established by Cousin for punctured spheres~\cite[Theorem~A]{cousin} and by Cousin--Heu for punctured surfaces of all genera~\cite[Theorem~A]{cousin-heu}. In rank $2$, mapping class group fixed points for $4$-punctured spheres are exactly the singular points of the relative character varieties (or the Fricke cubic surfaces), and correspond to Riccati solutions of Painlevé~\Rnum{6}~\cite[Theorem~5]{iwasaki}, see also~\cite[Section~9]{IIM}. Our classification is thus a first effort towards a complete picture of the algebraic solutions to rank-$3$ Garnier systems with monodromy in $\PU(2,1)$.

\subsection{Plan of the paper}
The classification of mapping class group fixed points $\pi_1S_{g,n}\to\PU(2,1)$ is divided into two parts: we treat the case $g=0$ first (Sections~\ref{sec:reducible-fixed-points-genus-0}--\ref{sec:irreducible-fixed-points-genus-0}) and then proceed to the case $g\geq 1$ (Section~\ref{sec:fixed-points-genus-at-lesast-one}).

\begin{figure}[h]
\includegraphics[scale=.4]{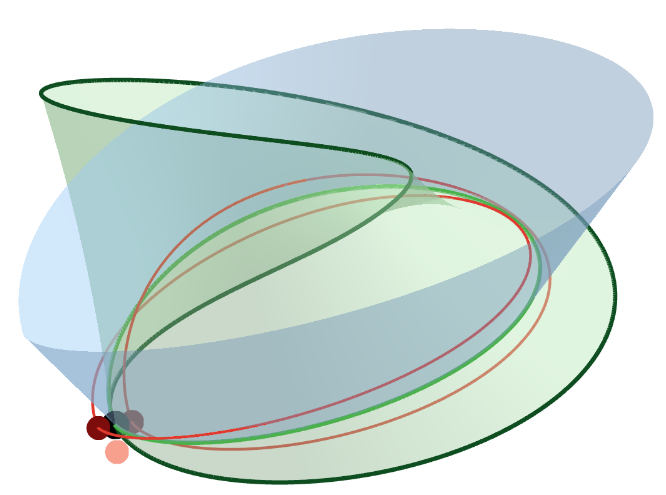}
\caption{The space of conjugacy classes in $\PU(2,1)$ as a semi-algebraic variety.}\label{fig:conjugacy-classes-intro}
\end{figure}

When $g=0$ and $n\leq 3$, every representation is a fixed point (because $\PMod(S_{0,n})$ is trivial in that case). When $n=4$, and unlike for rank-$2$ representations, there are irreducible fixed points with non-trivial peripheral monodromies arising as restrictions of reduced Gassner representations of $PB_4$ (Section~\ref{sec:fixed-points-genus-zero}). In order to classify all fixed points when $n=4$, we use that $\PMod(S_{0,4})$ is generated by two Dehn twists and reformulate the condition of being fixed by $\PMod(S_{0,4})$ in terms of intersections of centralizers in $\PU(2,1)$ (Lemma~\ref{lem:general-lemma-fixed-points}). We thus invest some time listing all possible intersections of two centralizers in $\PU(2,1)$ (see Table~\ref{tab:intersections-centralizers}) after describing all seven non-trivial types of conjugacy classes (Section~\ref{sec:isometries-CH2} and Table~\ref{tab:conjugacy-classes}) and illustrating them (see Figures~\ref{fig:conjugacy-classes-intro} and~\ref{fig:conjugacy-classes}). 
The essential notions and tools from complex hyperbolic geometry are recalled in Section~\ref{sec:complex-hyperbolic-geometry}. This enables us to classify all reducible fixed points $\pi_1S_{0,n}\to\PU(2,1)$ in Section~\ref{sec:reducible-fixed-points-genus-0} by treating the case $n=4$ first, and then inducting on $n$. The classification of irreducible fixed points happens in two steps. First, when $n=4$, we prove that peripheral monodromies of irreducible fixed points must be \emph{pseudo-reflections} (Definition~\ref{def:pseudo-reflections}), see Section~\ref{sec:reducing-to-pseudo-reflections}. Then, we show that the conditions of being irreducible and having only pseudo-reflections as peripheral monodromies (without assuming it is a fixed point) force the representation to be a Gassner fixed point in Section~\ref{sec:irreducible-fixed-points-genus-0}. We sketch a different argument in Appendix~\ref{apx:horn-pseudo-reflections} that shows that being irreducible and having pseudo-reflections as peripheral monodromies forces the representation to be rigid when $n=4$ (Proposition~\ref{prop:apx-irreducible-and-pseudo-reflections-imply-rigid}).

To classify the fixed points $\pi_1S_{g,n}\to\PU(2,1)$ when $g\geq 1$, we start with the case of the once-punctured torus (Section~\ref{sec:once-punctured-torus}). Using its hyperelliptic involution, we reduce the study of fixed points to the case of $4$-punctured spheres and prove that an irreducible fixed point is necessarily of Burau type. We also introduce the only non-trivial example of reducible fixed point: the quaternionic fixed points (Definition~\ref{def:quaternionic}). The key observation that allows us to deal with surfaces of genus $g\geq 1$ with arbitrarily many punctures is the following: a fixed point maps all non-separating simple closed curves to conjugate elements (because $\PMod(S_{g,n})$ acts transitively on them). This allows for efficient induction on the topology of the surface and a reduction of the study to punctured spheres and once-punctured tori (Sections~\ref{sec:twice-punctured-torus}--\ref{sec:higher-genus}).

\subsection{Acknowledgements}
This project benefited from fruitful discussions with Grace Garden and Stephan Tillmann. 
TLF thanks Geordie Williamson for inspiring discussions about Burau representations.
The authors would like to thank the Sydney Mathematical Research Institute and the University of Neuchâtel for their hospitality. SB acknowledges funding from the European Research Council (ERC) under the European Union's Horizon 2020 research and innovation programme (grant agreement No.~101018839). TLF acknowledges support from a SMRI Postdoctoral Fellowship, the CNRS UMR 6625, and the French Agence Nationale de la Recherche through the project ANR-25-CE40-1360 IsoMoDyn. AM was partially supported by the European Research Council (ERC) under the European Union’s Horizon 2020 research and innovation program (Grant agreement No.~101096550).

\subsection{AI statement}
AI agents have been used in this project to code illustrations from sketch, look for references, and check for typos and computation mistakes.
All the original mathematical ideas in this paper are our own.

\section{Complex hyperbolic space}\label{sec:complex-hyperbolic-geometry}
\subsection{Definition}
The group $\SU(2,1)$ is the group of $3\times 3$ complex matrices with determinant~$1$ that preserve a given Hermitian form of signature $(2,1)$. Depending on the context, we will consider different explicit Hermitian forms. The group $\SU(2,1)$ has a non-trivial centre that consists of scalar matrices, i.e.~$Z(\SU(2,1))\cong \Z/3\Z$.
\begin{defn}
The \emph{complex hyperbolic space} of (complex) dimension $2$, denoted by $\CH^2$, is the symmetric space associated to $\SU(2,1)$. In other words, 
\[
\CH^2 = \SU(2,1)/\mathrm{S}(\U(2)\times \U(1)).
\]
\end{defn}

Classic references on complex hyperbolic geometry include Goldman's book~\cite{goldman-book} or Parker's notes~\cite{parker-notes}. We will simply recall some notions that are useful for our work.

\begin{defn}
The maximal totally geodesic $1$-dimensional complex subspaces of $\CH^2$ are called \emph{$\C$-lines}. 
\end{defn}

\subsection{Models}
There are different models for $\CH^2$ like there are different models of the hyperbolic plane. We will use two of them, depending on the context.
\subsubsection{Ball model}
A standard choice of Hermitian form is 
\begin{equation}\label{eq:hermitian-form-diagonal}
J=\begin{pmatrix}
    1 & 0 & 0\\ 0 & 1 & 0\\ 0 & 0 & -1
\end{pmatrix}.
\end{equation}
The complex hyberbolic space $\CH^2$ corresponds to the projectivization of the space of vectors with negative $J$-norm in $\C^3$. A such, it can be naturally identified with the complex ball of $\C^2$:
\[
\CH^2=\{(z_1,z_2)\in\C^2 : |z_1|^2+|z_2|^2<1\}.
\]
The boundary at infinity $\partial\CH^2$ is the unit sphere in $\C^2$.

\subsubsection{Siegel model}
If we instead work with the anti-diagonal Hermitian form
\begin{equation}\label{eq:hermitian-form-antidiagonal}
J'=\begin{pmatrix}
    0 & 0 & 1\\ 0 & 1 & 0\\ 1 & 0 & 0
\end{pmatrix},
\end{equation}
the corresponding model of complex hyperbolic space is known as the \emph{Siegel model}:
\[
\CH^2=\{(z_1,z_2)\in\C^2: 2\Re(z_1)+|z_2|^2<0\}.
\]
The boundary at infinity in this case is given by
\[
\partial\CH^2 = \{(z_1,z_2)\in\C^2:2\Re(z_1)+|z_2|^2=0\} \cup \{\infty\},
\]
where $\infty$ is the projectivization of $(1,0,0)\in \C^{3}$. Letting $\zeta = z_2/\sqrt{2}$ and $v=\Im(z_1)$, we obtain an identification of $\partial\CH^2$ with $(\C\times \R)\cup\{\infty\}$ in which $(z_1,z_2)$ is identified with 
\[
(\zeta,v)=\begin{pmatrix}
    -|\zeta|^2+iv\\
    \sqrt{2}\zeta\\
    1
\end{pmatrix}.
\]

\subsubsection{Change of models}
It is convenient to have the explicit matrix that allows to switch between the Hermitian form $J$ and $J'$. One possible choice is
\begin{equation}\label{eq:change-of-Hermitian-form-PJ'P=J}
P=\frac{1}{\sqrt{2}}\begin{pmatrix}
    1 & 0 & 1\\ 0 & \sqrt{2} & 0\\ 1 & 0 & -1
\end{pmatrix}.
\end{equation}
One can easily verify that
\[
PJ'P=J\quad \text{and}\quad \overline{P}=P^t=P^{-1}=P.
\]

\subsection{Isometries}\label{sec:isometries-CH2}
The group of holomorphic isometries of $\CH^2$ naturally identifies with 
\[
\PU(2,1)=\SU(2,1)/Z(\SU(2,1)).
\]
\begin{notation}
    We will write $\PU(2,1)$ (or $\SU(2,1)$) when we refer to those groups with no particular Hermitian form in mind, and we will write $\PU(J)$ or $\PU(J')$ to highlight a particular Hermitian form. If $H$ is a subgroup of $\U(2,1)$ (typically $\U(2)\times\U(1)$, $\U(1)\times \U(1,1)$, or $\U(1)\times\U(1)\times\U(1)$), then by $P(H)$ we mean the quotient $H/Z(\U(2,1))$ (and not the quotient of $H$ by its centre).
\end{notation}

The space $\CH^2\cup\partial\CH^2$ being a closed ball, every element of $\PU(2,1)$ fixes at least one point in that space. 
\begin{defn}
    An element of $\PU(2,1)$ is called:
    \begin{itemize}
        \item \emph{elliptic} if it fixes at least one point in $\CH^2$;
        \item \emph{parabolic} if it fixes exactly one point in $\partial \CH^2$;
        \item \emph{hyperbolic} if it fixes exactly two points in $\partial\CH^2$.
    \end{itemize}
\end{defn}
It is not hard to see that an element $A$ of $\PU(2,1)$ that fixes at least three points in $\partial\CH^2$ must be elliptic. Indeed, if the fixed points are contained in the boundary of the same $\C$-line $L$, then $A$ is the identity on $L$, hence a complex reflection. Otherwise, $A$ is the identity. Also note that elliptic elements are exactly those contained in a compact subgroup of $\PU(2,1)$.

Our next task consists in describing all different kinds of conjugacy classes in $\PU(2,1)$. A detailed treatment of this task was first achieved by Chen--Greenberg in~\cite[Section~3 and Theorem~3.4.1]{chen-greenberg}. We summarize the list of conjugacy classes and our terminology in Table~\ref{tab:conjugacy-classes}. The original contribution of the section is the illustration of the space of conjugacy classes in $\PU(2,1)$ (i.e.~the non-Hausdorff quotient $\PU(2,1)/\text{conj}$) provided by Figure~\ref{fig:conjugacy-classes}.

\subsubsection{Elliptic}
Elliptic elements of $\PU(2,1)$ are those that fix
at least one point in $\CH^2$. They all belong to a compact subgroup of $\PU(2,1)$. It is convenient to work in the ball model and express elliptic elements in $\PU(J)$. Up to conjugation, they are of the form
\begin{equation}\label{eq:elliptic-elements}
    E(\alpha,\beta)=\begin{pmatrix}
        e^{i\alpha} & 0& 0\\
        0&e^{i\beta} &0\\
        0&0&1
    \end{pmatrix},
\end{equation}    
with $\alpha, \beta$ in $[0,2\pi)$. Their conjugacy classes are parametrized by the unordered pairs $\{\alpha,\beta\}$, identifying the space of conjugacy classes of elliptic elements in $\PU(2,1)$ with a Möbius band. Excluding the identity, there are three different families of conjugacy classes of elliptic elements in $\PU(2,1)$.

\paragraph{(RE) Regular elliptic elements} Those are the elements that fix a unique point $X\in\CH^2$ and preserve exactly two orthogonal $\C$-lines passing through $X$. They fix no points in $\partial\CH^2$. Up to conjugation, they are of the form $E(\alpha,\beta)$ with $\alpha\neq \beta$ in $(0,2\pi)$. 
\begin{itemize}
    \item \textbf{Parametrization.}
Two regular elliptic elements $E(\alpha, \beta)$ and $E(\alpha', \beta')$ are conjugate if and only if $\{\alpha,\beta\}=\{\alpha',\beta'\}$. The conjugacy classes of regular elliptic elements are thus parametrized by $\{(\alpha,\beta)\in (0,2\pi)^2 : \alpha < \beta\}$.
    \item \textbf{Centralizer.}
The centralizer of a regular elliptic elements consists of all elements in $\PU(2,1)$ that stabilize its two invariant $\C$-lines (or, in other words, that have the same eigenlines). For instance,
\[
Z(E(\alpha,\beta))=\{E(\alpha',\beta'): (\alpha',\beta')\in [0,2\pi)^2\}.
\]
Note that $Z(E(\alpha,\beta))$ is isomorphic to the 2-dimensional torus $\rm{P}(U(1)\times U(1)\times U(1))$.
\end{itemize} 

\paragraph{(SR) Scalar rotations}
They fix a unique point $X\in\CH^2$ and globally preserve every $\C$-line passing through $X$. They fix no point in $\partial \CH^2$. Scalar rotations are, up to conjugation, of the form
\[
E(\alpha,\alpha)=\begin{pmatrix}
    e^{i\alpha} &0&0\\
    0&e^{i\alpha}&0\\
    0&0&1
\end{pmatrix},
\]
with $\alpha\in (0,2\pi)$. All the scalar rotations that fix the same point $X\in\CH^2$, together with the identity, form a subgroup
\[
\mathcal{SR}_X\subset P(U(2)\times U(1))\subset\PU(2,1)
\]
isomorphic to $\U(1)$.

\begin{itemize}
    \item \textbf{Parametrization.}
Two scalar rotations $E(\alpha, \alpha)$ and $E(\alpha', \alpha')$ are conjugate if and only if $\alpha=\alpha'$. Conjugacy classes of scalar rotations are thus parametrized by $\alpha\in (0,2\pi)$.
    \item\textbf{Centralizer.}
The centralizer of a scalar rotation with fixed point $X$ consists of all elements of $\PU(2,1)$ that fix $X$. For instance, 
\[
Z(E(\alpha,\alpha))=\mathrm{P} (\U(2)\times\U(1)).
\]
\end{itemize}

\paragraph{(CR) Complex reflections}
These are the elliptic elements that fix pointwise a unique $\C$-line $L\subset\CH^{2}$ called the \emph{mirror}. They also preserves all the $\C$-lines that are perpendicular to their mirror.
The set of fixed points in $\partial \CH^2$ correspond to set of boundary points $\partial L$ of $L$.
Complex reflections are, up to conjugation, of the form
\[
E(\alpha,0)=\begin{pmatrix}
    e^{i\alpha}&0&0\\
    0&1&0\\
    0&0&1
\end{pmatrix},
\] 
with $\alpha\in (0,2\pi)$. All complex reflections fixing the same complex $\C$-line $L$, together with the identity, form a subgroup
\[
\mathcal{CR}_L\subset P(U(1)\times U(1,1))\subset \PU(2,1)
\]
isomorphic to $\U(1)$.

\begin{itemize}
    \item \textbf{Parametrization.}
Two complex reflections $E(\alpha, 0)$ and $E(\alpha', 0)$ are conjugate if and only if $\alpha=\alpha'$. Conjugacy classes of complex reflections are thus parametrized by $\alpha\in (0,2\pi)$.

\item \textbf{Centralizer.}
The centralizer of a complex reflection with mirror $L$ corresponds to the stabilizer of its mirror, i.e.~$\stab(L)$. For instance,
\[
Z(E(\alpha,0))=\mathrm{P}(\U(1)\times \U(1,1)).
\]
\end{itemize}

\subsubsection{Parabolic}
Parabolic elements of $\PU(2,1)$ are those that fix a unique point on $\partial\CH^2$. Like for elliptic elements, we distinguish three different kinds of conjugacy classes of parabolic elements. Two of them are made of unipotent elements (regular parabolic and vertical translations), and the last kind consists of non-unipotent elements (screw parabolic). It is convenient to work in the Siegel model and express parabolic elements in $\PU(J')$.

\begin{defn}
The parabolic elements of $\PU(2,1)$ that are unipotent are called \emph{Heisenberg translations}. All Heisenberg translations fixing the same point $Y\in\partial\CH^2$, together with the identity, form a subgroup 
\[
\mathcal{N}_Y\subset\PU(2,1).
\]
\end{defn}
For instance, the Heisenberg translations fixing $\infty$ in the Siegel model are of the form
\begin{equation}\label{eq:Heisenberg-translations}
A(\zeta,v) = \begin{pmatrix}
    1 & -\sqrt{2}\,\overline{\zeta} & -|\zeta|^2+iv\\
    0 & 1 & \sqrt{2}\zeta\\
    0 & 0 & 1
\end{pmatrix}
\end{equation}
with $(\zeta,v)\in \C\times \R$. Note that $\mathcal{N}_Y$ is isomorphic to the classical Heisenberg group. Indeed, when $Y=\infty$, 
\[
A(\zeta,v)A(\xi,t)=A\big(\zeta+\xi, v+t+2\Im(\overline{\xi}\zeta)\big).
\]

\paragraph{Horosphere foliations}
All parabolic elements fixing the same point $Y\in\partial \CH^2$ preserve each leaf of a foliation $H=\{H_u\}_{u\in\R}$ of $\CH^2$ into horospheres based at $Y$. When $Y=\infty$, the horospheres can be described as follows. For every $u\in\R_+$, define $H_u$ as the set of all $(z_1,z_2)\in \CH^2$ such that $2\Re(z_1)+|z_2|^2=-2u$. In other words,
\[
H_u=\{(-|\xi|^2-u+it, \sqrt{2}\xi)\in\CH^2:\xi\in\C, t\in\R\}.
\]
The horospheres $H_u$ foliate $\CH^2$ in \emph{horosphere coordinates} and identify $\CH^2$ with $\C\times\R\times\R_+$. Note that $A(\zeta,v)$ maps $(\xi,t,u)$ to $(\xi+\zeta, t+v+2\Im(\overline{\xi}\zeta), u)$. The $\C$-lines through $\infty$ are indexed by $z_0\in\C$ and can be written as follows in horosphere coordinates:
\[
L_{z_0}=\{(z_0,t,u):t\in\R, u\in\R_+\}.
\]
The stabilizer of each horosphere $H_u$ of the foliation $H$ based at $Y\in\partial\CH^2$ consists of all Heisenberg translations fixing $Y$, all complex reflections through each $L_{z_0}$, and all screw parabolic elements fixing $Y$ (and the identity). It can be expressed as
\[
\stab(H)=\mathcal{N}_Y\rtimes \mathcal{CR}_{L},
\]
where $L$ is any $\C$-line through $Y$.

\paragraph{(RP) Regular parabolic} These are the Heisenberg translations that fix a unique point $Y\in\partial\CH^2$, but preserve no complex $\C$-lines. Those fixing $\infty$ are of the form $A(\zeta,v)$ with $\zeta\in \C\setminus \{0\}$ and $v\in \R$. Their minimal polynomial is $(X-1)^3$.

\begin{itemize}
    \item \textbf{Representatives.} All Heisenberg translations $A(\zeta,v)$ with $\zeta\neq 0$ are conjugate. A representative of their conjugacy class is 
\[
\begin{pmatrix}
    1 & -1 & -\frac{1}{2}\\
    0 & 1 & 1\\
    0 & 0 & 1
\end{pmatrix}.
\]
None of them can be separated from the identity in $\PU(2,1)$.

\item \textbf{Centralizer.} The centralizer of a regular parabolic element is the set of all Heisenberg translations fixing the same boundary point and with horizontal translations ``in the direction $\zeta$''. For instance,
\[
Z(A(\zeta,v))=\{A(\xi,t): \xi=\lambda \zeta, \lambda\in\R\}\cong \R\times \R.
\]
\end{itemize}

\paragraph{(VT) Vertical translations}
All elements of $\PU(2,1)$ that are conjugate to a Heisenberg translations $A(0,v)$ with $v\neq 0$ are called \emph{vertical translations}. They preserve all the $\C$-lines through their unique fixed point on $\partial\CH^2$. Their minimal polynomial is $(X-1)^2$. All the vertical translations fixing the same point $Y\in\partial\CH^2$ together with the identity form a subgroup
\[
\mathcal{VT}_Y\subset \mathcal{N}.
\]

\begin{itemize}
\item \textbf{Representatives.} There are exactly two conjugacy classes of vertical translations determined by the sign of $v$. The standard representatives are
\[
\begin{pmatrix}
    1 & 0 & \pm i\\
    0 & 1 & 0\\
    0 & 0 & 1
\end{pmatrix}.
\]

\item \textbf{Centralizer.} The centralizer of a vertical translation consists of all elements of $\PU(2,1)$ that fix the same fixed point and stabilize the same horosphere foliation, i.e.~it is equal to $\stab(H)$. For instance,
\[
Z(A(0,v))=\{A(\xi,t):\xi\in\C, t\in\R\}\rtimes \left\{\begin{pmatrix}
    1 & 0 & 0\\
    0 & e^{i\theta} & 0\\
    0 & 0 & 1
\end{pmatrix} : \theta\in\R/2\pi\Z\right\} \cong \mathcal{N}_Y\rtimes \U(1).
\]
Moreover, $\mathcal{VT}_Y$ is also the centralizer of $\stab(H)$ in $\PU(2,1)$.
\end{itemize}

\begin{defn}
    The parabolic elements of $\PU(2,1)$ that are \emph{not} unipotent are called \emph{screw parabolic}.
\end{defn}

\paragraph{(SP) Screw parabolic}
They fix a unique point $Y\in \partial\CH$ and also preserve exactly one $\C$-line through $Y$. The ones fixing $\infty$ are of the form
\[
A(\zeta, v, \theta)=\begin{pmatrix}
    1 & \sqrt{2}\,\overline{\zeta}(1-e^{i\theta}) & -2|\zeta|^2(1-e^{i\theta})+iv\\
    0 & e^{i\theta} & \sqrt{2}\zeta (1-e^{i\theta})\\
    0 & 0 & 1
\end{pmatrix},
\]
with $\theta\in (0,2\pi)$ (if $\theta=0$, they are vertical translations) and $(\zeta,v)\neq (0,0)$ (else, they would be complex reflections). 

\begin{itemize}
    \item \textbf{Representatives.} We can conjugate $A(\zeta,v,\theta)$ to have $\zeta=0$, and further conjugate to have $v=\pm 1$. So, the only conjugacy invariants are $\theta$ and the sign of $v$. In other words, each conjugacy class of screw parabolic is represented by an element of the form
\[
\begin{pmatrix}
    1 & 0 & \pm i\\
    0 & e^{i\theta} & 0\\
    0 & 0 & 1
\end{pmatrix},
\]
for $\theta\in (0,2\pi)$.

\item \textbf{Centralizer.} The centralizer of a screw parabolic consists of all elements fixing the same point $Y\in\partial\CH^2$ and stabilizing the same $\C$-line and the same horosphere foliation $H$. For instance,
\[
Z(A(\zeta,v,\theta))=\{A(\zeta, t, \sigma): t\in\R, \theta\in\R/2\pi\Z\}.
\]
Note that this centralizer contains elliptic elements, namely complex reflections through the invariant $\C$-line.
\end{itemize}

\subsubsection{Hyperbolic}
Hyperbolic elements correspond to those
fixing exactly two points~$Y_1,Y_2$ in~$\partial\CH^2$.
They also preserve the unique $\C$-line asymptotic to $Y_1$ and $Y_2$.
As for parabolic elements, we favour the Siegel model and the Hermitian form $J'$ to work with hyperbolic elements.
The two fixed points give two $J'$-isotropic eigenlines
in~$\C^{2,1}$ and so they are diagonalisable,
admit a representative of the form
\[
H(\lambda)=
\begin{pmatrix}
    \lambda & 0 & 0\\
    0& 1 & 0\\
    0&0&\overline\lambda^{-1}
\end{pmatrix},
\]
with $\lambda\in\C\setminus\{0\}$. The two eigenlines correspond to the two fixed points of $H(\lambda)$ in $\partial\CH^2$. The uniqueness of the two isotropic eigenlines is equivalent to $|\lambda|\neq 1$. If we write $\lambda = re^{-i\theta}$ with $r\neq 1$, then, in $\PU(2,1)$, $H(\lambda)$ is further equal to\footnote{The equivalent form of $H(r,\theta)$ for the Hemritian form $J$ is 
\[
\begin{pmatrix}
    \cosh(r) & 0 & \sinh(r)\\
    0 & e^{i\theta} & 0\\
    \sinh(r) & 0 & \cosh(r)
\end{pmatrix}.
\]
} 
\[
H(r,\theta)=\begin{pmatrix}
    r & 0 & 0\\
    0 & e^{i\theta} & 0\\
    0 & 0 & r^{-1}
\end{pmatrix}.
\]

\begin{itemize}
    \item \textbf{Parametrization}
Two hyperbolic elements $H(\lambda)$ and $H(\lambda')$ are conjugate if and only if $\lambda=\lambda'$ or $\lambda\overline{\lambda'}=1$. Conjugacy classes of hyperbolic elements are thus parametrized by $\lambda\in \C$ with $|\lambda|> 1$ (or $r>1$ and $\theta\in [0,2\pi)$ for $\lambda=re^{-i\theta}$), i.e.~the complement of the closed unit disk in $\C$. The quantity $\log(r)$ is the \emph{translation length}.

\item \textbf{Centralizer}
The centralizer of a hyperbolic element corresponds to all the elements of $\PU(2,1)$ that fix $Y_1$ and $Y_2$. In other words, it consists of all hyperbolic elements with fixed points $Y_1$ and $Y_2$, as well as all complex reflections through $L$ (and the identity). In other words,
\[
Z(H(\lambda))=\big\{H(\lambda') : \lambda'\in \C\setminus\{0\}\big\}.
\]
\end{itemize}

\subsubsection{Other families of isometries}
We have just seen that there are seven different kinds of conjugacy classes in $\PU(2,1)$ if we exclude the identity: three elliptic (RE, SR, CR), three parabolic (RP, VT, SP), and one hyperbolic (H). There are two sub-families of isometries that are useful to introduce: Heisenberg similarities and pseudo-reflections.

\paragraph{Heisenberg similarities}
\begin{defn}\label{def:Heisenberg-similarities}
    The elements of $\PU(2,1)$ that fix a given point in $\partial\CH^2$ are called \emph{Heisenberg similarities}.
\end{defn}
The group of Heisenberg similarities fixing the same point in $Y\in\partial\CH^2$ is denoted by $\stab(Y)\supsetneq \stab(H)$ and is isomorphic to
\[
(\R_+\times \U(1))\ltimes \mathcal{N}_Y,
\]
where the factor $\R_+\times \U(1)\cong \C^\times$ corresponds to the stabilizer of two distinct points in $\partial\CH^2$. It consists of all Heisenberg translations fixing $Y$, all complex reflections across a $\C$-line through $Y$, and all hyperbolic elements fixing $Y$. 

\paragraph{Pseudo-reflections}
\begin{defn}\label{def:pseudo-reflections}
An element $X$ of $\mathrm{GL}_3\C$ is a \emph{pseudo-reflection} if $X-1$ has rank $1$. Similarly, an element $X\in\PU(2,1)$ is a \emph{pseudo-reflection} if it admits a lift $\widehat{X}$ to $\mathrm{U}(2,1)$ that is a pseudo-reflection.
\end{defn}

It is not hard to see that the pseudo-reflections of $\PU(2,1)$ consist of all complex reflections, scalar reflections, and vertical translations.

\subsubsection{Summary and illustration}
\begin{figure}[hb]
\includegraphics[scale=.5]{fig/conjugacy-classes-2.png}
\includegraphics[scale=.5]{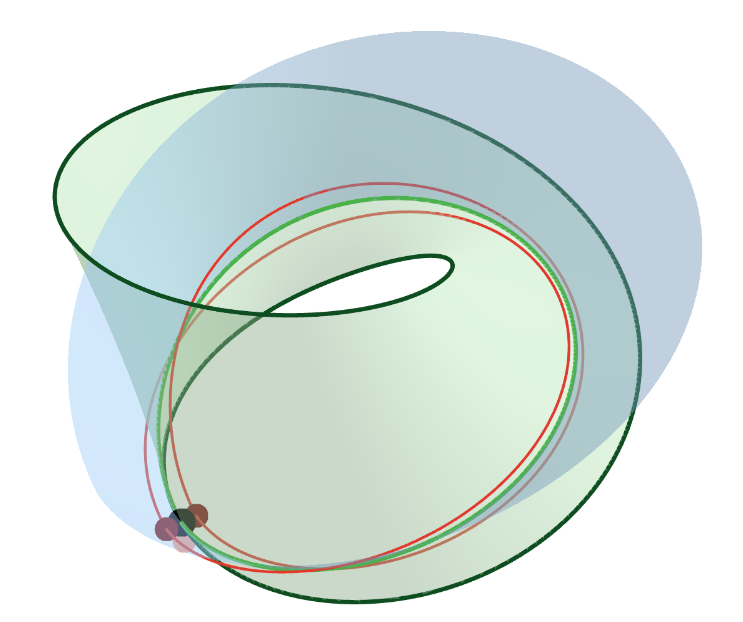}
\caption{The topological quotient $\PU(2,1)/\text{conj}$.}\label{fig:conjugacy-classes}
\end{figure}

The classification of conjugacy classes in $\PU(2,1)$ is summarized in Table~\ref{tab:conjugacy-classes}. It is possible to visualize the topological quotient $\PU(2,1)/\text{conj}$ (in other words, the character variety of representations of $\Z$ into $\PU(2,1)$) as follows.
\begin{enumerate}
    \item Start with the elliptic conjugacy classes. They form a Möbius band (in green on Figure~\ref{fig:conjugacy-classes}). Its boundary (dark green) corresponds to scalar rotations. Inside it passes a circle (lighter green) that corresponds to complex reflections. It meets the boundary at the identity (black dot).
    \item Next, the hyperbolic conjugacy classes form a cylinder (blue) attached to the Möbius band along the inner circle of complex reflections. It does not intersect the Möbius band at any other points (which is why $\PU(2,1)/\text{conj}$ can only be immersed in 3-dimensional space).
    \item Finally, the parabolic conjugacy classes consist first of two circles (red) that correspond to screw parabolic elements and cannot be separated from the circle of complex reflections. The two points that cannot be separated from the identity (red dots) are the two classes of vertical translations. The regular parabolic conjugacy classes form another point (lighter red) that cannot be separated from the identity.
\end{enumerate}
The largest Hausdorff quotient of $\PU(2,1)/\text{conj}$ is obtained by removing all the parabolic conjugacy classes (all the red components on Figure~\ref{fig:conjugacy-classes}). It is a $2$-dimensional real semi-algebraic variety immersed in 3-dimensional space (it can be embedded in 4-dimensional space).

\begin{landscape}
\begin{table}[ht]
    \centering
    \caption{Conjugacy classes in $\PU(2,1)$.}\label{tab:conjugacy-classes}
    {\tabulinesep=.6mm
    \begin{tabular}{c||c|c||c|c||c|c}
        \toprule
        $\PU(2,1)$ & representative & parameters & invariants in $\CH^2$ & fixed data in $\partial\CH^2$ & centralizer & isomorphism type \\
        \midrule
        
        \makecell{regular\\ elliptic} & $\begin{pmatrix}
            e^{i\alpha} & 0 & 0\\ 0 & e^{i\beta} & 0\\ 0 & 0 & 1
        \end{pmatrix}$ & $\alpha\neq\beta\in (0,2\pi)$ & \makecell{unique fixed point $X$,\\ two invariant $\C$-lines\\$L_1\perp L_2$ through $X$} & none & \makecell{$\stab{L_1}\cap\stab{L_2}$} & $\mathrm{P}(\U(1)\times \U(1)\times \U(1))$ \\ \hline
        
        \makecell{complex\\ reflections} & $\begin{pmatrix}
            1 & 0 & 0\\ 0 & e^{i\alpha} & 0\\ 0 & 0 & 1
        \end{pmatrix}$ & $\alpha\in (0,2\pi)$ & \makecell{unique fixed $\C$-line $L$,\\ all $\C$-lines $L'\perp L$\\ are invariant} & fixed circle $\partial L$ & $\stab L$ & $\mathrm{P}(\U(1)\times \U(1,1))$ \\\hline
        
        \makecell{scalar\\ rotations} & $\begin{pmatrix}
            e^{i\alpha} & 0 & 0\\ 0 & e^{i\alpha} & 0\\ 0 & 0 & 1
        \end{pmatrix}$ & $\alpha\in (0,2\pi)$ & \makecell{unique fixed point $X$,\\ all $\C$-lines through $X$ \\are invariant} & none & $\stab X$ & $\mathrm{P}(\U(2)\times \U(1))$ \\ \midrule
        
        \makecell{regular\\ parabolic} & $\begin{pmatrix}
            1 & -1 & \frac{-1}{2}\\ 0 & 1& 1\\ 0 & 0 & 1
        \end{pmatrix}$ & none & \makecell{invariant horosphere \\foliation $H$ based at $Y$} & unique fixed point $Y$ & \makecell{Heisenberg translations in\\
        $\stab H$ with horizontal\\ translation in same direction $P$} & \makecell{$\begin{pmatrix}
            1 & -u & \frac{-u^2+iv}{2}\\ 0 & 1 & u\\
            0 & 0 & 1
        \end{pmatrix}$\\ $(u,v)\in \R\times\R$} \\ \hline
        
        \makecell{vertical\\ translation} & $\begin{pmatrix}
            1 & 0 & \varepsilon i\\ 0 & 1& 0\\ 0 & 0 & 1
        \end{pmatrix}$ & $\varepsilon=\pm1$ & \makecell{invariant horosphere \\foliation $H$ based at $Y$,\\all $\C$-lines $L'$ through \\$Y$ are invariant} & unique fixed point $Y$ & $\stab H$ & $\mathcal{N}_Y\rtimes \mathcal{CR}_{L'}$ \\ \hline
        
        \makecell{screw\\ parabolic} & $\begin{pmatrix}
            1 & 0 & \varepsilon i\\ 0 & e^{i\theta}& 0\\ 0 & 0 & 1
        \end{pmatrix}$ & \makecell{$\varepsilon=\pm1$\\ $\theta\in (0,2\pi)$} & \makecell{invariant horosphere \\foliation $H$ based at $Y$,\\one invariant $\C$-line $L$\\ through $Y$} & unique fixed point $Y$ & $\stab L\cap \stab H$ & \makecell{$\mathcal{VT}_Y\times \mathcal{CR}_L$\\
        \\$\begin{pmatrix}
            1 & 0 & \frac{iv}{2}\\ 0 & e^{i\phi}& 0\\ 0 & 0 & 1
        \end{pmatrix}$\\ $(\phi, v)\in [0,2\pi)\times \R$} \\
         \midrule
         
         \makecell{hyperbolic} & $\begin{pmatrix}
            r & 0 & 0\\ 0 & e^{i\theta}& 0\\ 0 & 0 & r^{-1}
        \end{pmatrix}$ & \makecell{$r>1$\\ $\theta\in [0,2\pi)$} & \makecell{one invariant $\C$-line $L$\\ through $Y_1,Y_2$} & two fixed points $Y_1,Y_2$ & $\stab Y_1\cap \stab Y_2$ & \makecell{$\begin{pmatrix}
            u & 0 & 0\\ 0 & e^{i\phi}& 0\\ 0 & 0 & u^{-1}
        \end{pmatrix}$\\ $(\phi, u)\in [0,2\pi)\times \R$} \\
        
        \bottomrule
    \end{tabular}
    }
\end{table}
\end{landscape}

\begin{landscape}
\begin{table}[ht]
    \caption{Intersections of centralizers in $\PU(2,1)$.}
    \label{tab:intersections-centralizers}
{%
\tiny
\setlength{\tabcolsep}{3pt}
\renewcommand{\arraystretch}{1.5}
\begin{tabular}{c||c|c|c|c|c|c|c}
\toprule
{\large $Z\cap Z'$} & \begin{tabular}[b]{@{}c@{}} regular elliptic (RE)\\ $X'$, $L_1'$, $L_2'$\end{tabular} &
\begin{tabular}[b]{@{}c@{}} complex reflection (CR) \\ $L'$\end{tabular} &
\begin{tabular}[b]{@{}c@{}}scalar rotation (SR)\\ $X'$\end{tabular} &
\begin{tabular}[b]{@{}c@{}} regular parabolic (RP)\\ $Y'$, $P'$\end{tabular} &
\begin{tabular}[b]{@{}c@{}}vertical translation (VT)\\ $Y'$\end{tabular} &
\begin{tabular}[b]{@{}c@{}} screw parabolic (SP)\\ $Y'$, $L'$\end{tabular} &
\begin{tabular}[b]{@{}c@{}} hyperbolic (H)\\ $Y_1'$, $Y_2'$, $L'$\end{tabular} \\
\midrule\midrule

\begin{tabular}[t]{@{}c@{}}RE\\ $X$, $L_1$, $L_2$\end{tabular} &
\begin{tabular}[t]{@{}l@{\ $\to$\ }l@{}}
$X=X'$, $\{L_i\}=\{L_j'\}$ & $Z=Z'$ \\
$X=X'$, $\{L_i\}\neq \{L_j'\}$ & $\mathcal{SR}_X$ \\
$X\neq X'\in L_i= L_j'$ & $\mathcal{CR}_{XX'}$ \\
otherwise & $1$
\end{tabular} &
\begin{tabular}[t]{@{}l@{\ $\to$\ }l@{}}
$L'=L_i$ & $Z\subset Z'$ \\
$X\in L'$, $L'\neq L_i$ & $\mathcal{SR}_X$ \\
$X\notin L'$, $L'\perp L_i$ & $\mathcal{CR}_{L_i}$ \\
otherwise & $1$
\end{tabular} &
\begin{tabular}[t]{@{}l@{\ $\to$\ }l@{}}
$X'=X$ & $Z\subset Z'$ \\
$X'\in L_i$ & $\mathcal{CR}_{L_i}$ \\
otherwise & $1$
\end{tabular} &
\begin{tabular}[t]{@{}l@{\ $\to$\ }l@{}}
always & $1$
\end{tabular} &
\begin{tabular}[t]{@{}l@{\ $\to$\ }l@{}}
$Y'\in\partial L_i$ & $\mathcal{CR}_{L_i}$ \\
otherwise & $1$
\end{tabular} &
\begin{tabular}[t]{@{}l@{\ $\to$\ }l@{}}
$L'=L_i$ & $\mathcal{CR}_{L'}$ \\
otherwise & $1$
\end{tabular} &
\begin{tabular}[t]{@{}l@{\ $\to$\ }l@{}}
$L'=L_i$ & $\mathcal{CR}_{L'}$ \\
otherwise & $1$
\end{tabular} \\
\hline
\begin{tabular}[t]{@{}c@{}}CR\\ $L$\end{tabular} &
&
\begin{tabular}[t]{@{}l@{\ $\to$\ }l@{}}
$L=L'$ & $Z=Z'$ \\
$L\cap_\infty L'$ & $\mathcal{VT}_{Y}$ \\
$L\cap L'$, $L\perp L'$ & $\mathcal{CR}_{L}\mathcal{CR}_{L'}$ \\
$L\cap L'$, $L\not\perp L'$ & $\mathcal{SR}_{X}$ \\
$L\parallel L'$ & $\mathcal{CR}_{M}$
\end{tabular} &
\begin{tabular}[t]{@{}l@{\ $\to$\ }l@{}}
$X'\in L$ & $\mathcal{CR}_{L}\mathcal{CR}_{M}$ \\
$X'\notin L$ & $\mathcal{CR}_{M}$
\end{tabular} &
\begin{tabular}[t]{@{}l@{\ $\to$\ }l@{}}
$Y'\in\partial L$ & $\mathcal{VT}_{Y'}$ \\
$Y'\notin\partial L$ & $1$
\end{tabular} &
\begin{tabular}[t]{@{}l@{\ $\to$\ }l@{}}
$Y'\in\partial L$ & $\mathcal{VT}_{Y'}\mathcal{CR}_{L}$ \\
$Y'\notin\partial L$ & $\mathcal{CR}_{M}$
\end{tabular} &
\begin{tabular}[t]{@{}l@{\ $\to$\ }l@{}}
$L=L'$ & $Z'\subset Z$ \\
$L\neq L'$, $Y'\in\partial L$ & $\mathcal{VT}_{Y'}$ \\
$L\perp L'$ & $\mathcal{CR}_{L'}$ \\
otherwise & $1$
\end{tabular} &
\begin{tabular}[t]{@{}l@{\ $\to$\ }l@{}}
$L=L'$ & $Z'\subset Z$ \\
$L\perp L'$ & $\mathcal{CR}_{L'}$ \\
otherwise & $1$
\end{tabular} \\
\hline
\begin{tabular}[t]{@{}c@{}}SR\\ $X$\end{tabular} &
& &
\begin{tabular}[t]{@{}l@{\ $\to$\ }l@{}}
$X=X'$ & $Z=Z'$ \\
$X\neq X'$ & $\mathcal{CR}_{XX'}$
\end{tabular} &
\begin{tabular}[t]{@{}l@{\ $\to$\ }l@{}}
always & $1$
\end{tabular} &
\begin{tabular}[t]{@{}l@{\ $\to$\ }l@{}}
always & $\mathcal{CR}_{XY'}$
\end{tabular} &
\begin{tabular}[t]{@{}l@{\ $\to$\ }l@{}}
$X\in L'$ & $\mathcal{CR}_{L'}$ \\
$X\notin L'$ & $1$
\end{tabular} &
\begin{tabular}[t]{@{}l@{\ $\to$\ }l@{}}
$X\in L'$ & $\mathcal{CR}_{L'}$ \\
$X\notin L'$ & $1$
\end{tabular} \\
\hline
\begin{tabular}[t]{@{}c@{}}RP\\ $Y$, $P$\end{tabular} &
& & &
\begin{tabular}[t]{@{}l@{\ $\to$\ }l@{}}
$Y=Y'$, $P=P'$ & $Z=Z'$ \\
$Y=Y'$, $P\neq P'$ & $\mathcal{VT}_{Y}$ \\
$Y\neq Y'$ & $1$
\end{tabular} &
\begin{tabular}[t]{@{}l@{\ $\to$\ }l@{}}
$Y=Y'$ & $Z\subset Z'$ \\
$Y\neq Y'$ & $1$
\end{tabular} &
\begin{tabular}[t]{@{}l@{\ $\to$\ }l@{}}
$Y=Y'$ & $\mathcal{VT}_{Y}$ \\
$Y\neq Y'$ & $1$
\end{tabular} &
\begin{tabular}[t]{@{}l@{\ $\to$\ }l@{}}
always & $1$
\end{tabular} \\
\hline
\begin{tabular}[t]{@{}c@{}}VT\\ $Y$\end{tabular} &
& & & &
\begin{tabular}[t]{@{}l@{\ $\to$\ }l@{}}
$Y=Y'$ & $Z=Z'$ \\
$Y\neq Y'$ & $\mathcal{CR}_{YY'}$
\end{tabular} &
\begin{tabular}[t]{@{}l@{\ $\to$\ }l@{}}
$Y=Y'$ & $Z'\subset Z$ \\
$Y\neq Y'$, $Y\in\partial L'$ & $\mathcal{CR}_{L'}$ \\
$Y\notin\partial L'$ & $1$
\end{tabular} &
\begin{tabular}[t]{@{}l@{\ $\to$\ }l@{}}
$Y\in\partial L'$ & $\mathcal{CR}_{L'}$ \\
$Y\notin\partial L'$ & $1$
\end{tabular} \\
\hline
\begin{tabular}[t]{@{}c@{}}SP\\ $Y$, $L$\end{tabular} &
& & & & &
\begin{tabular}[t]{@{}l@{\ $\to$\ }l@{}}
$Y=Y'$, $L=L'$ & $Z=Z'$ \\
$Y=Y'$, $L\neq L'$ & $\mathcal{VT}_{Y}$ \\
$Y\neq Y'$, $L=L'$ & $\mathcal{CR}_{L}$ \\
$Y\neq Y'$, $L\neq L'$ & $1$
\end{tabular} &
\begin{tabular}[t]{@{}l@{\ $\to$\ }l@{}}
$L=L'$ & $\mathcal{CR}_{L}$ \\
$L\neq L'$ & $1$
\end{tabular} \\
\hline
\begin{tabular}[t]{@{}c@{}}H\\ $Y_1$, $Y_2$, $L$\end{tabular} &
& & & & & &
\begin{tabular}[t]{@{}l@{\ $\to$\ }l@{}}
$\{Y_i\}=\{Y_j'\}$ & $Z=Z'$ \\
$L=L'$, $Y_i\neq Y_j'$ & $\mathcal{CR}_{L}$ \\
$L\neq L'$ & $1$
\end{tabular} \\
\bottomrule
\end{tabular}}
\end{table}
\begin{rem}
Table~\ref{tab:intersections-centralizers} indicates all possible intersections of centralizers in $\PU(2,1)$ according to the relative positions of the invariant data. The table is symmetric. Invariant data is recalled in the first row (primed data) and the first column (unprimed data) according to the notation of Table~\ref{tab:conjugacy-classes}. The inner cells are formated as \texttt{configuration}$\to$\texttt{value of intersection}. We use the following notation. Let $X\neq X'\in\CH^2$, $Y\in\partial\CH^2$, and $L, L'$ be two $\C$-lines : (i) $XX'$ denotes the unique $\C$-line through $X$ and $X'$, (ii) $M$ stands for the unique $\C$-line through a given point and perpendicular to a given $\C$-line (which are always clear from context), (iii) $L\cap_\infty L'$ means that the lines are asymptotic on $\partial\CH^2$, (iv) $L\cap L'$ means that the lines intersect once in $\CH^2$, (v) $L\perp L'$ means that the lines intersect perpendicularly in $\CH^2$, (vi) $L\parallel L'$ means that the lines are ultraparallel, (vii) $\mathcal{SR}_X$ stands for the scalar rotations fixing $X$, (viii) $\mathcal{CR}_{L}$ stands for the complex reflections across the $\C$-line $L$, and (ix) $\mathcal{VT}_Y$ stands for the vertical translations fixing $Y$.
\end{rem}
\end{landscape}

\subsection{Toolbox}
We collect some facts on $\CH^2$ and its isometries that will be useful later.

\subsubsection{Intersections of centralizers}
We described all possible intersections of the centralizers of two non-trivial elements in $\PU(2,1)$, according to their conjugacy classes, in Table~\ref{tab:intersections-centralizers}. The following statements are direct consequences of the classification of all centralizers intersections.

\begin{cor}\label{cor:non-abelian-intersection-centralizers}
    If $A$ and $B$ are non-trivial elements of $\PU(2,1)$ such that $Z(A)\cap Z(B)$ is non-abelian, then they must be of one of the following kinds:
    \begin{itemize}
        \item $A$ and $B$ are complex reflections with the same mirror $L$ and $Z(A)\cap Z(B)=\stab(L)$;
        \item $A$ and $B$ are scalar reflections fixing the same point $X$ and $Z(A)\cap Z(B)=\stab(X)$;
        \item $A$ and $B$ are vertical translations with same invariant horosphere foliation $H$ and the intersection is~$Z(A)\cap Z(B)=\stab(H)$.
    \end{itemize}
\end{cor}

\begin{cor}\label{cor:intersection-centralizers-pseudo-reflections-non-trivial}
    For $A$ and $B$ are pseudo-reflections (Definition~\ref{def:pseudo-reflections}) in $\PU(2,1)$, then the intersection~$Z(A)\cap Z(B)$ is always non-trivial.
\end{cor}

\subsubsection{Intersections of stabilizers}

\begin{lem}\label{lem:intersection-sstabilizers-two-C-line-is-abelian}
    If $L_1$ and $L_2$ are two distinct $\C$-lines in $\CH^2$, then $\stab(L_1)\cap \stab(L_2)$ is abelian.
\end{lem}
\begin{proof}
    We consider the relative positions of $L_1$ and $L_2$.
    \begin{itemize}
        \item If $L_1$ and $L_2$ intersect at a point $X$ in $\CH^2$, then either $L_1\perp L_2$ and $\stab(L_1)\cap \stab(L_2)$ is isomorphic to $\mathrm{P}(\U(1)\times \U(1)\times \U(1))$, or $L_1\not\perp L_2$ and $\stab(L_1)\cap \stab(L_2)=\mathcal{SR}_X\cong \U(1)$. 
        \item If $L_1$ and $L_2$ intersect at a point $Y$ in $\partial\CH^2$, then $\stab(L_1)\cap \stab(L_2)=\mathcal{VT}_Y\cong \R$.
        \item If $L_1$ and $L_2$ are ultra-parallel, then $\stab(L_1)\cap \stab(L_2)=\mathcal{CR}_L\cong \U(1)$, where $L$ is the unique $\C$-line perpendicular to both $L_1$ and $L_2$.
    \end{itemize}
    In all cases, the intersection is abelian.
\end{proof}

\begin{lem}\label{lem:intersection-sstabilizers-one-C-line-and-one-bdry-point}
    If $L$ is a $\C$-line in $\CH^2$ and $Y\in \partial \CH^2$, then $\stab(L)\cap \stab(Y)$ can be described as follows.
    \begin{itemize}
        \item If $Y\in \partial L$, then $\stab(L)\cap \stab(Y)= \mathcal{VT}_Y\rtimes (\mathcal{CR}_L\times \R)$, where the second factor correspond to the stabilizer of $Y$ and any other boundary point.
        \item If $Y\notin \partial L$, then $\stab(L)\cap \stab(Y)=\mathcal{CR}_M$, where $M$ is the unique $\C$-line perpendicular to $L$ and passing through $Y$.
    \end{itemize}
\end{lem}
\begin{proof}
Recall that $\stab(L)\cap \stab(Y)$ contains $\stab(L)\cap \stab(H)$, which is the intersection of two centralizers. According to Table~\ref{tab:intersections-centralizers}, $\stab(L)\cap \stab(H)$ is either equal to $\mathcal{VT}_Y\mathcal{CR}_L$ when $Y\in\partial L$ or to $\mathcal{CR}_M$ when $Y\notin \partial L$. In the first case, we need to take into account the dilatations of $H$ that stabilize $L$ in order to get all of $\stab(L)\cap \stab(Y)$, while in the second case no dilatations of $H$ stabilize $L$.
\end{proof}

\subsubsection{Product of isometries}
\begin{lem}\label{lem:product-pseudo-reflections}
    If $A$ and $B$ are two pseudo-reflections in $\PU(2,1)$ (Definition~\ref{def:pseudo-reflections}) such that $AB$ is the identity or a pseudo-reflection, then
    $A$ and $B$ commute.
\end{lem}
\begin{proof}
    Since $A,B$ are pseudo-reflections, there are lifts $\widehat A,\widehat B$ to $\mathrm{U}(2,1)$ such that
    $\widehat A-\id = A_0$ and  $\widehat B-\id=B_0$
    with $\mathrm{rk}(A_0)=\mathrm{rk}(B_0)=1$. Note that $\widehat A$ and $\widehat B$ being unitary, it follows that $\ker(A_0)=\Im(A_0)^\perp$ and $\ker(B_0)=\Im(B_0)^\perp$. We also have
    \begin{equation*}
        \widehat A\widehat B-\id = A_0+B_0+A_0B_0,
    \end{equation*}
    which means that $\widehat A\widehat B-\id$ vanishes on $\ker(A_0)\cap \ker(B_0)$. Now, note that the rank assumption implies that $\ker(A_0)\cap \ker(B_0)$ has dimension at least $1$, so $\widehat A\widehat B-\id$ has rank at most $2$.

    If it has rank $0$, then $AB=1$, so $A$ and $B$ commute.
    
    Assume now that $\mathrm{rk}(\widehat A\widehat B-\id)=1$. If $\ker(A_0)\neq \ker(B_0)$, then for any $v\in \ker(A_0)\setminus \ker(B_0)$ and $w\in \ker(B_0)\setminus \ker(A_0)$, we have that $A_0w$ and $\widehat{A}B_0v$ both belong to $\Im(\widehat A\widehat B-\id)$. Hence, $\widehat{A}B_0v\in\Im(\widehat A\widehat B-\id)= \Im(A_0)$ and since $\widehat{A}^{-1}$ preserves $\Im(A_0)$, we conclude that $\Im(B_0)=\Im(A_0)$. Taking orthogonal complements, we get $\ker(A_0)=\ker(B_0)$, contradiction. We conclude that $\ker(A_0)=\ker(B_0)$ and $\Im(A_0)=\Im(B_0)$. This means that $A$ and $B$ are either complex reflections with same mirror, or scalar rotations around the same point, or vertical translations fixing the same boundary point. In all cases, $A$ and $B$ commute.

    Finally, assume that $\mathrm{rk}(\widehat A\widehat B-\id)=2$. Since $AB$ is a pseudo-reflection, there is $\mu\in\C^\ast$ such that $\mathrm{rk}(\widehat A\widehat B-\mu\id)=1$ and $\mu\neq 1$. Let $v$ be a non-zero vector in $\ker(B_0)\cap \ker (\widehat A\widehat B-\mu\id)$. Note that $\mu v= \widehat A\widehat B v = \widehat A v$ and thus $A_0v=(\mu-1)v\neq 0$. This shows that $\Im(A_0)=\langle v\rangle\subset \ker(B_0)$, and consequently $\Im(A_0)\perp \Im(B_0)$ and $\Im(B_0)\subset \ker(A_0)$. Therefore, $A_0B_0=B_0A_0=0$ and $\widehat A\widehat B =\id + A_0 + B_0 = \widehat B \widehat A$, proving that $A$ and $B$ commute.
\end{proof}
\subsubsection{Maximal subgroups}
Throughout this section, $X$ will denote a point of $\CH^2$, $L$ will denote a $\C$-line in $\CH^2$, $H$ will denote a horosphere foliation based at $Y\in \partial\CH^2$, and $R$ will denote a maximal totally real totally geodesic 2-dimensional subspace of $\CH^2$. We recall the following:
\begin{itemize}
    \item $\stab(X)\cong \mathrm{P}(\U(2)\times\U(1))$ is the set of isometries fixing $X$;
    \item $\stab(L)\cong\mathrm{P}(\U(1)\times\U(1,1))$ is the set of isometries stabilizing $L$,
    \item $\stab(H)\cong \mathcal{N}_Y\rtimes \mathcal{CR}_{L'}$, where $L'$ is any $\C$-line passing through $Y$, is the set of isometries stabilizing each leaf of $H$ (in particular, it does not contain any hyperbolic isometry);
    \item $\stab(Y)\cong (\R_+\times \mathcal{CR}_{L'})\ltimes \mathcal{N}_Y$ is the set of all Heisenberg similarities fixing $Y$, i.e.~the set of all isometries fixing $Y$;
    \item $\stab(R)\cong\mathrm{PO}(2,1)$ is the set of isometries stabilizing $R$.
\end{itemize}

\begin{lem}\label{lem:maximal-subgroups-of-PU(2,1)-with-non-trivial-centralizer}
    If $\Lambda$ is a maximal subgroup of $\PU(2,1)$ with non-trivial centralizer, then $\Lambda$ is of one of the following forms:
    \begin{enumerate}
        \item 
        $\Lambda = \stab(X)$ with centralizer $\mathcal{SR}_X\cong \U(1)$ consisting of all scalar rotations fixing $X$.
        \item 
        $\Lambda= \stab(L)$ with centralizer $\mathcal{CR}_L\cong \U(1)$ consisting of all complex reflections with mirror~$L$.
        \item $\Lambda=\stab(H)$ with centralizer $\mathcal{VT}_Y\cong \R$
        consisting of all vertical translations fixing $Y$.
    \end{enumerate}
In particular, $Z(\Lambda)$ is always abelian.
\end{lem}
\begin{proof}
    For any $g\in Z(\Lambda)$, we have
    $\Lambda\subset Z(g)$ and $g\in Z(Z(g))$. So, if $g\neq 1$, then $Z(Z(g))\neq \{1\}$ and by maximality of $\Lambda$, we must have $\Lambda = Z(g)$. This means that $\Lambda$ must be one of the centralizers appearing in Table~\ref{tab:conjugacy-classes} and by maximality it must be of the form $\stab(X)$, $\stab(L)$, or $\stab(H)$, as desired.
\end{proof}

\begin{defn}\label{defn:irreducible}
A subgroup $\Lambda\subset \PU(2,1)$ is \emph{reducible} if its projective action on $\mathbb{P}(\C^3)=\mathbb{CP}^2$ fixes a point and it is \emph{irreducible} otherwise.
\end{defn} 

\begin{lem}\label{lem:maximal-reducible-subgroups-of-PU(2,1)}
    If $\Lambda$ is a maximal reducible subgroup of $\PU(2,1)$, then $\Lambda$ is of one of the following forms:
    \begin{enumerate}
        \item $\Lambda = \stab(X)$.
        \item $\Lambda = \stab(L)$.
        \item $\Lambda = \stab(Y)$.
    \end{enumerate}
\end{lem}
\begin{proof}
    The fixed point of $\Lambda$ in $\mathbb{P}(\C^3)$ corresponds to a vector in $\C^{2,1}$. If it has negative norm, then $\Lambda$ fixes the corresponding point in $\CH^2$. Similarly, if it is a null vector, then $\Lambda$ fixes the corresponding point in $\partial\CH^2$. Finally, if it has positive norm, then $\Lambda$ stabilizes the complex $\C$-line that corresponds to its orthogonal complement in $\C^{2,1}$.
\end{proof}

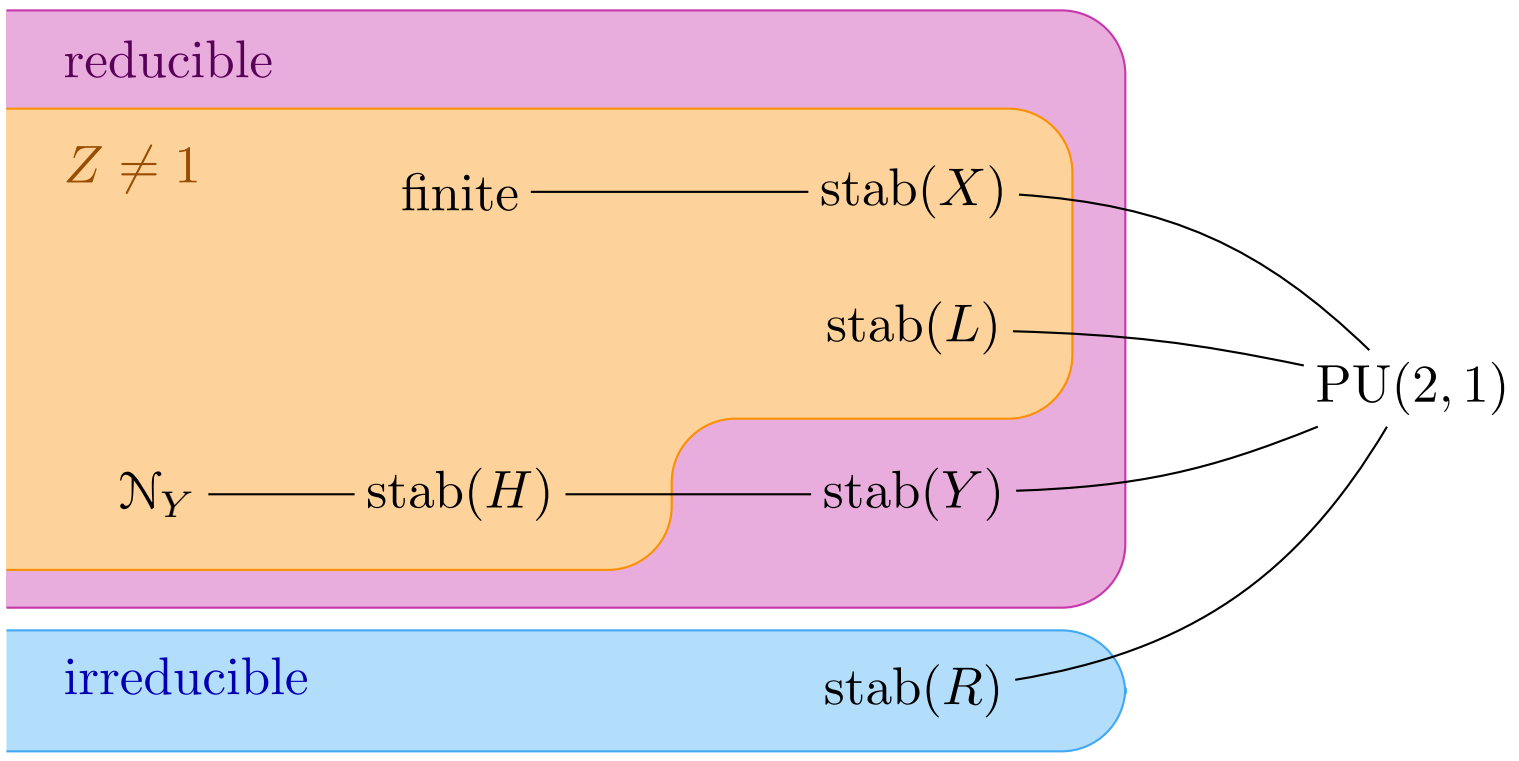
\begin{figure}[ht]
    \centering
    \begin{tikzpicture}[every node/.style={inner sep=2pt}]
    \begin{scope}
    \clip (0,-3) rectangle (11,3.5);
    \fill[lightmauve, rounded corners=12pt] (-1,-1.25) rectangle (7.4,2.7);
    \fill[lightapricot, rounded corners=12pt]
      (-1,-1) -- (-1,2.05) -- (7.05,2.05) -- (7.05,0) -- (4.4,0) -- (4.4,-1) -- cycle;
    \fill[lightsky, rounded corners=12pt]
    (0,-2.2) -- (7.4,-2.2) -- (7.4,-1.4) -- (0,-1.4);
    \end{scope}

    \draw[mauve, rounded corners=12pt]
    (0,2.7) -- (7.4,2.7) -- (7.4,-1.25) -- (0,-1.25);
    \draw[apricot, rounded corners=12pt]
    (0,2.05) -- (7.05,2.05) -- (7.05,0) -- (4.4,0) -- (4.4,-1) -- (0,-1);
    \draw[sky, rounded corners=12pt]
    (0,-2.2) -- (7.4,-2.2) -- (7.4,-1.4) -- (0,-1.4);

    \node[anchor=west, violet!70!black] at (0.3,2.38) {reducible};
    \node[anchor=west, orange!60!black] at (0.3,1.65) {$Z\neq 1$};
    \node[anchor=west, blue!70!black] at (0.3,-1.7) {irreducible};

    \node (F) at (3,1.5)    {finite};
    \node (X) at (6,1.5)    {$\stab(X)$};
    \node (L) at (6,0.6)    {$\stab(L)$};
    \node (H) at (3,-0.5)   {$\stab(H)$};
    \node (N) at (1,-0.5)   {$\mathcal{N}_Y$};
    \node (Y) at (6,-0.5)   {$\stab(Y)$};
    \node (R) at (6,-1.8)   {$\stab(R)$};
    \node (G) at (9.3,0.2)  {$\mathrm{PU}(2,1)$};

    \draw (F) -- (X);
    \draw (H) -- (Y);
    \draw (N) -- (H);
    \draw (X) edge[-, bend left=20] (G);
    \draw (L) edge[-, bend left=5] (G);
    \draw (Y) edge[-, bend right=10] (G);
    \draw (R) edge[-, bend right=25] (G);
\end{tikzpicture}
    \caption{Some remarkable Zariski closed subgroups of $\PU(2,1)$ (including all maximal proper ones). A line indicates an inclusion of the left group into the right group.}
    \label{fig:Zariski-closed-subgroups}
\end{figure}

\begin{lem}\label{lem:maximal-Zariski-closed-subgroups-of-PU(2,1)}
    If $\Lambda$ is a maximal, closed proper subgroup of $\PU(2,1)$, then $\Lambda$ is of one of the following forms:
    \begin{enumerate}
        \item $\Lambda$ is discrete.
        \item $\Lambda=\stab(X)$.
        \item $\Lambda = \stab(L)$.
        \item $\Lambda = \stab(Y)$.
        \item $\Lambda = \stab(R)$.
    \end{enumerate}
    If we further assume $\Lambda$ to be Zariski closed, then only (2)--(5) apply.
\end{lem}
\begin{proof}
    This is a particular case of~\cite[Theorem~4.4.2]{chen-greenberg}. 
\end{proof}
Remarkable Zariski closed subgroups of $\PU(2,1)$ are illustrated on Figure~\ref{fig:Zariski-closed-subgroups}. An immediate consequence of Lemma~\ref{lem:maximal-Zariski-closed-subgroups-of-PU(2,1)} is the following result.
\begin{cor}\label{cor:Zariski-dense-implies-discrete-or-dense}
    A Zariski dense subgroup of $\PU(2,1)$ is either dense or discrete.
\end{cor}

\section{Mapping class group fixed points in genus zero}\label{sec:fixed-points-genus-zero}
\subsection{Background}
Let $S_{0,n}$ denote an oriented and connected sphere with $n\geq 3$ points removed. Its \emph{pure mapping class group}, denoted by $\PMod_{0,n}$, consists of all isotopy classes of orientation-preserving homeomorphisms of $S_{0,n}$ that fix each puncture individually. The Dehn--Nielsen--Baer Theorem identifies it with a subgroup $\Out^\star(\pi_1S_{0,n})$ of the group $\Out(\pi_1S_{0,n})$ of all outer automorphisms of $\pi_1S_{0,n}$ (see e.g.~\cite[Theorem~8.8]{farb-margalit}). More precisely, $\Out^\star(\pi_1S_{0,n})$ consists of all outer automorphisms that preserve the conjugacy classes of peripheral loops around each puncture of $S_{0,n}$. This identification defines an action of $\PMod_{0,n}$ on conjugacy classes of representations $\rho\colon\pi_1S_{0,n}\to\PU(2,1)$ by pre-composition, known as \emph{mapping class group action on character variety}. We are interested in understanding the fixed points of this action.
\begin{notation}
    Given a group $\Lambda$, we write
    \[
    \mathcal{X}(\Lambda)=\Hom(\Lambda,\PU(2,1))/\PU(2,1)
    \]
    for the set of conjugacy classes of representations $\Lambda\to\PU(2,1)$, endowed with the quotient topology. When $\Lambda=\pi_1 S_{0,n}$, we abbreviate $\mathcal{X}_{0,n}=\mathcal{X}(\pi_1S_{0,n})$ and call it \emph{character variety} (even though we only consider the topological quotient). The set of fixed points for the action of $\PMod_{0,n}$ on $\mathcal{X}_{0,n}$ will be denoted by
    \[
    \FixedPts_{0,n} \subset \mathcal{X}_{0,n}.
    \]
\end{notation}
\begin{terminology}
    When we say that a representation $\rho\colon \pi_1S_{0,n}\to\PU(2,1)$ is a \emph{fixed point}, we mean that its conjugacy class $[\rho]$ is fixed by the action of $\PMod_{0,n}$. In particular, unless explicitly stated, we never refer to the action of $\Aut^\star(\pi_1S_{0,n})$ on representations, but only to the action of $\Out^\star(\pi_1S_{0,n})$ on conjugacy classes of representations.
\end{terminology}

The group $\PMod_{0,n}$ is naturally isomorphic to the pure braid group on $n-1$ strands $PB_{n-1}$ modulo its centre (see e.g.~\cite[Section~9.3]{farb-margalit}):
\[
\PMod_{0,n}\cong PB_{n-1}/Z(B_{n-1}).
\]
It also fits into a short exact sequence known as the \emph{Birman exact sequence}:
\begin{equation}\label{eq:Birman-short-exact-sequence}
1\longrightarrow \pi_1 S_{0,n}\longrightarrow \PMod_{0,n+1}\longrightarrow \PMod_{0,n}\longrightarrow 1.
\end{equation}
The first arrow is the \emph{point-pushing map} and the second arrow forgets the extra puncture, see~\cite[Theorem~4.6 and the discussion thereafter]{farb-margalit} for more details. The Birman exact sequence can be rewritten as the purely algebraic short exact sequence
\[
1\longrightarrow \pi_1 S_{0,n}\longrightarrow \Aut^\star(\pi_1S_{0,n})\longrightarrow\Out^\star(\pi_1S_{0,n})\longrightarrow 1.
\]
In particular, there is a natural isomorphism $\PMod_{0,n+1}\cong \Aut^\star(\pi_1 S_{0,n})$. A remarkable feature of the Birman exact sequence in genus $0$ is that it splits. We include an explanation for the reader's convenience.
\begin{lem}\label{lem:Birman-splits}
    There exists a section $s\colon \PMod_{0,n}\to \PMod_{0,n+1}$ that splits the Birman exact sequence~\eqref{eq:Birman-short-exact-sequence} for every $n\geq 3$. Moreover, $s$ can be chosen such that its image, seen as a subgroup of $\Aut^\ast(\pi_1 S_{0,n})$, fixes a peripheral loop in $\pi_1 S_{0,n}$. 
\end{lem}
\begin{proof}
    The Birman exact sequence can also be seen as coming from the long exact sequence in homotopy associated to the fibration
    \[
    S_{0,n}\to \mathcal{M}_{0,n+1}\to \mathcal{M}_{0,n},
    \]
    where $\mathcal{M}_{0,n}$ is the moduli space of $n$-punctured spheres. Identifying $\mathcal{M}_{0,n}$ with the configuration space of $n-3$ distinct points on $\mathbb{CP}^1\setminus\{0,1,\infty\}$, we can define a continuous section $\mathcal{M}_{0,n}\to \mathcal{M}_{0,n+1}$ by
    \[
    (z_1,\ldots,z_{n-3})\mapsto (z_1,\ldots,z_{n-3},z),\quad \text{where } z=2+\sum_{i=1}^{n-3}|z_i|. 
    \]
    Note that $z$ is real, distinct from $0,1$, and distinct from each $z_i$. The induced map on fundamental groups gives the desired splitting section.

    Fix a base point $(z_1,\ldots,z_{n-3},z)$ in $\mathcal{M}_{0,n+1}$ and consider the loop $\gamma(t)=(z_1,\ldots,z_{n-3},z e^{2\pi i t})$. It is a peripheral loop around the puncture at infinity of the fibre $\CP^1\setminus\{0,1,\infty,z_1,\ldots,z_{n-3}\}$. On the other hand, a loop $\alpha(t)$ in $\mathcal{M}_{0,n}$ based at $(z_1,\ldots, z_{n-3})$ is mapped by $s$ to the loop $\overline{\alpha}(t)=(\alpha(t), z(\alpha(t))$ in $\mathcal{M}_{0,n+1}$. The loops $\gamma$ and $\overline{\alpha}$ commute in $\pi_1\mathcal{M}_{0,n+1}$, giving the desired conclusion.
\end{proof}

\subsection{Mapping class group representations}
Mapping class group fixed points and representations of mapping class groups are closely related. Most of this correspondence is well-known; we include it for the reader's convenience before proving more specific statements for $\PU(2,1)$ representations. 

For now, let $G$ be an arbitrary centre-free group.

\begin{lem}\label{lem:restriction-of-PMod-representations-are-fixed-points}
If $\rho\colon\pi_1S_{0,n}\to G$ is the restriction of a representation $\widehat{\rho}\colon\PMod_{0,n+1}\to G$, then $\rho$ is a fixed point.
\end{lem}
\begin{proof}
Recall that $\PMod_{0,n+1}\cong\Aut^\star(\pi_1 S_{0,n})$. To show that $\rho$ is a fixed point, we have to prove that $\rho\circ f$ is conjugate to $\rho$ for every $f\in \Aut^\star(\pi_1 S_{0,n})$. Via the point-pushing map $\pi_1S_{0,n}\to \Aut^\star(\pi_1 S_{0,n})$, we map $\gamma$ to the inner automorphism $\imath_\gamma(x)=\gamma x\gamma^{-1}$. Now, note that
\[
(\rho\circ f)(\gamma)=\rho(f(\gamma))=\widehat{\rho}(\imath_{f(\gamma)})=\widehat{\rho}(f\circ \imath_\gamma\circ f^{-1})=\widehat{\rho}(f)\widehat{\rho}(\imath_\gamma)\widehat{\rho}(f)^{-1}=\widehat{\rho}(f)\rho(\gamma)\widehat{\rho}(f)^{-1}.
\]
As $\gamma$ was arbitrary, we conclude that $\rho\circ f$ is conjugate to $\rho$ by $\widehat{\rho}(f)$.
\end{proof}

Lemma~\ref{lem:restriction-of-PMod-representations-are-fixed-points} defines a map
\[
\res\colon\mathcal{X}(\PMod_{0,n+1})\to \FixedPts_{0,n},
\]
which also makes sense for representations into any group. This map is not expected to be surjective for arbitrary $G$. However, some representations always lie in its image.

\begin{defn}\label{def:regular-representations}
    A representation $\rho\colon\Lambda\to G$ is \emph{regular} if the centralizer of $\rho(\Lambda)$ in $G$ is trivial. 
\end{defn} 
For instance, when $G$ is (semi-)algebraic, all Zariski dense representations are regular.
\begin{notation}
    The centralizer of the image of a representation $\rho\colon\Lambda\to G$ will be denoted by $Z_\rho$, and its normalizer by $N_\rho$. 
    
    We will also write $\mathcal{X}^\reg(\Lambda)\subset \mathcal{X}(\Lambda)$, and similarly $\mathcal{X}_{0,n}^\reg\subset \mathcal{X}_{0,n}$ and $\FixedPts_{0,n}^\reg\subset\FixedPts_{0,n}$, for the subsets of regular representations.
\end{notation}
\begin{lem}\label{lem:regular-and-abelian-fixed-points-are-restrictions}
If $\rho\colon\pi_1S_{0,n}\to G$ is a fixed point and $\rho$ is regular, then it is the restriction of a unique representation $\widehat{\rho}\colon\PMod_{0,n+1}\to G$. Moreover, $\widehat{\rho}$ is regular.

On the other hand, if $\rho\colon\pi_1S_{0,n}\to G$ is abelian, then $\rho$ is a fixed point. Moreover, $\res^{-1}([\rho])$ is non-empty and can be identified with  $\Hom(\PMod_{0,n}, Z_\rho)/Z_\rho$.
\end{lem}
\begin{proof}
First, assume that $\rho\colon\pi_1S_{0,n}\to G$ is regular and a fixed point. For every $f\in \Aut^\star(\pi_1 S_{0,n})$, we define $\widehat{\rho}(f)$ as the unique element $g_f$ of $G$ such that $\rho\circ f=g_f\rho g_f^{-1}$. Note that $g_f$ exists because $\rho$ is a fixed point and $g_f$ is unique because $\rho$ is regular. It is not difficult to verify that $\widehat{\rho}$ defines a group morphism and extends $\rho$. If $\widehat{\rho}_1$ and $\widehat{\rho}_2$ were two lifts of $\rho$, then the previous computation would give, for every $f\in  \Aut^\star(\pi_1 S_{0,n})$,
\[
\rho\circ f = \widehat{\rho}_1(f)\cdot\rho\cdot \widehat{\rho}_1(f)^{-1}= \widehat{\rho}_2(f)\cdot\rho\cdot \widehat{\rho}_2(f)^{-1},
\]
which implies that $\widehat{\rho}_2(f)^{-1}\widehat{\rho}_1(f)\in Z_\rho=\{1\}$ since $\rho$ is regular. Hence, $\widehat{\rho}_1=\widehat{\rho}_2$.

Now, assume that $\rho\colon\pi_1S_{0,n}\to G$ is abelian. Since $\Aut^\star(\pi_1 S_{0,n})$ preserves the conjugacy classes of peripheral loops, it fixes $\rho$. Using the section $s\colon \PMod_{0,n}\to\PMod_{0,n+1}\cong \Aut^\star(\pi_1 S_{0,n})$ from Lemma~\ref{lem:Birman-splits}, we may write $\PMod_{0,n+1}\cong \pi_1S_{0,n}\rtimes \PMod_{0,n}$. In particular, a representation $\widehat{\rho}\colon \PMod_{0,n+1}\to G$ is the data of a pair $\rho\colon\pi_1S_{0,n}\to G$ together with a morphism $\sigma\colon \PMod_{0,n}\to G$ such that, for every $\tau\in\PMod_{0,n}$,
    \begin{equation}\label{eq:reformulation-representation-of-semi-direct-produc}
    \rho \circ s(\tau)=\sigma(\tau)\rho\sigma(\tau)^{-1}.
    \end{equation}
When $\rho$ is abelian, $\rho\circ s(\tau)=\rho$ (because $s(\tau)$ preserves the conjugacy classes of peripheral loops) and exactly the morphisms $\sigma\colon\PMod_{0,n}\to Z_\rho$ (which include the trivial one) satisfy this condition.
\end{proof}


We may also describe other fibres of the restriction map, assuming they are non-empty.

\begin{lem}\label{lem:fibers-restriction-map-over-abelian-centralizers}
     If $\rho\colon\pi_1S_{0,n}\to G$ is the restriction of $\widehat{\rho}\colon\PMod_{0,n+1}\to G$ and $Z_\rho$ is abelian, then $\res^{-1}([\rho])$ is isomorphic to the twisted cohomology group $H^1(\PMod_{0,n}; Z_\rho)$ (the $\PMod_{0,n}$-module structure on $Z_\rho$ is detailed in the proof).
\end{lem}
\begin{proof}
    Note that $\rho$ is a fixed point by Lemma~\ref{lem:restriction-of-PMod-representations-are-fixed-points}. So, for every $f\in \Aut^\star (\pi_1 S_{0,n})$, there exists $g_f\in G$ such that $\rho\circ f =g_f \rho g_f^{-1}$. Note that $g_f\in N_\rho$ and it is uniquely defined up to right multiplication by an element of $Z_\rho$. This defines a morphism
    \begin{equation}\label{eq:morphism-eta-bar}
    \overline{\eta}_\rho \colon \Aut^\star (\pi_1 S_{0,n}) \to N_\rho/Z_\rho.
    \end{equation}
    Recall that $Z_\rho$ is a normal subgroup of $N_\rho$. The conjugation action of $N_\rho$ on $Z_\rho$ defines a morphism $\overline{\theta}_\rho\colon N_\rho\to \Aut(Z_\rho)$, which descends to a morphism $\theta_\rho\colon N_\rho/Z_\rho\to \Aut(Z_\rho)$ because $Z_\rho$ is assumed to be abelian.\footnote{In general, if $Z_\rho$ is not abelian, then we get a morphism $\theta_\rho\colon N_\rho/Z_\rho\to \Out(Z_\rho)$ instead.} In other words, we obtain a morphism $\overline{\varphi}_\rho\colon \Aut^\star (\pi_1 S_{0,n})\to \Aut(Z_\rho)$ defined by $\overline{\varphi}_\rho=\theta_\rho\circ \overline{\eta}_\rho$. Note that $\overline{\varphi}_\rho(\imath_\gamma)=\id$ for every inner automorphism $\imath_\gamma$ (although $\overline\eta_\rho(\rho(\gamma))\neq \id$ in general), so that there is an induced morphism
    \[
    \varphi_\rho \colon \underbrace{\Out^\star(\pi_1S_{0,n})}_{\cong \PMod_{0,n}}\to \Aut(Z_\rho).
    \]
    The morphism $\varphi_\rho$ defines a $\PMod_{0,n}$-module structure on $Z_\rho$ which we use to define the cohomology group $H^1(\PMod_{0,n};Z_\rho)$.

    Now, if $\widehat{\rho}_1$ and $\widehat{\rho}_2$ are two lifts in $\res^{-1}([\rho])$, then, as in the proof of Lemma~\ref{lem:regular-and-abelian-fixed-points-are-restrictions}, we have $\widehat{\rho}_2(f)^{-1}\widehat{\rho}_1(f)\in Z_\rho$ for every $f\in\Aut^\star(\pi_1S_{0,n})$. In other words, there is a map $\overline{\sigma}\colon \Aut^\star(\pi_1S_{0,n})\to Z_\rho$ defined by $\overline{\sigma}(f)=\widehat{\rho}_2(f)^{-1}\widehat{\rho}_1(f)$. A simple computation shows that
    \[
    \overline{\sigma}(fg)=\widehat{\rho}_2(g)^{-1}\widehat{\rho}_2(f)^{-1}\widehat{\rho}_1(f)\widehat{\rho}_1(g)=\overline{\sigma}(g)\cdot \overline{\varphi}_{\rho}(g)\big(\overline{\sigma}(f) \big). 
    \]
    Hence, $\overline{\sigma}$ defines a cocycle. Furthermore, for any inner automorphism $\imath_\gamma$, we have $\widehat{\rho}_1(\imath_\gamma)=\widehat{\rho}_2(\imath_\gamma)=\rho(\gamma)$, so that $\overline{\sigma}(\imath_\gamma)=1$. Therefore, $\overline{\sigma}$ descends to a cocycle
    \[
    \sigma\colon \Out^\star(\pi_1S_{0,n})\to Z_\rho,
    \]
    that defines a cohomology class in $H^1(\PMod_{0,n};Z_\rho)$. A careful analysis of all quotients involved in this construction shows that $H^1(\PMod_{0,n};Z_\rho)$ is the right parameter space for $\res^{-1}([\rho])$.
\end{proof}

We now specialize to the case $G=\PU(2,1)$ and prove that every fixed point is the restriction of a representation of mapping class group. To the best of our knowledge, this result is new, although non-regular representations may have simply been left aside in previous works.
\begin{thm}\label{thm:all-fixed-points-lift-PU(2,1)}
    When $G=\PU(2,1)$, the map $\res\colon \mathcal{X}(\PMod_{0,n+1})\to \FixedPts_{0,n}$ is surjective. 
\end{thm}
\begin{proof}
        We start by refining the observations we made before. In the proof of Lemma~\ref{lem:regular-and-abelian-fixed-points-are-restrictions}, when considering abelian representations, we used the section $s\colon\PMod_{0,n}\to \PMod_{0,n+1}$ from Lemma~\ref{lem:Birman-splits} to show that a fixed point $\rho\colon \pi_1S_{0,n}\to\PU(2,1)$ is the restriction of a morphism $\widehat{\rho}\colon \PMod_{0,n+1}\to\PU(2,1)$ if and only if there exists a morphism $\sigma\colon\PMod_{0,n}\to \PU(2,1)$ such that, for every $\tau\in\PMod_{0,n}$,
        \begin{equation}\label{eq:condition-lift}
        \rho\circ s(\tau)=\sigma(\tau)\rho \sigma(\tau)^{-1}.
        \end{equation}
        While it is not obvious that such a morphism $\sigma$ exists a priori, what we always have is a morphism
        \[
        \eta_\rho \colon \PMod_{0,n}\to N_\rho/Z_\rho,
        \]
        obtained by composing the morphism $\overline{\eta}_\rho$ from~\eqref{eq:morphism-eta-bar} with the section $s$. If we can lift $\eta_\rho$ to a morphism $\sigma\colon \PMod_{0,n}\to N_\rho$, then it will automatically satisfy~\eqref{eq:condition-lift}.

        The target of the map $\eta_\rho$ can be further restricted. According to Lemma~\ref{lem:Birman-splits}, the image of the section $s$ fixes a peripheral curve $\gamma\in\pi_1 S_{0,n}$. Evaluating~\eqref{eq:condition-lift} at $\gamma$, we see that $\eta_\rho(\tau)$ belongs to $A_\gamma=N_\rho\cap Z(\rho(\gamma))$ for every $\tau\in\PMod_{0,n}$. In other words, we really want to lift
        \[
        \eta_\rho \colon \PMod_{0,n}\to A_\gamma/Z_\rho
        \]
        to a morphism $\sigma\colon \PMod_{0,n}\to A_\gamma$.
        Note that, up to relabelling the punctures, we may always pick $\gamma$ as our favourite peripheral curve in $\pi_1 S_{0,n}$.

        \begin{claim}\label{claim:A_gamma-abelian-is-enough}
            If $A_\gamma$ is abelian, then $\eta_\rho$ always lifts to the desired morphism $\sigma$.
        \end{claim}
        \begin{proofclaim}
            If $A_\gamma$ is abelian, then $\eta_\rho$ factorizes through the abelianization $\PMod_{0,n}^\mathrm{ab}$ of $\PMod_{0,n}$. Since $\PMod_{0,n}^\mathrm{ab}$ is a free abelian group~(see e.g.~\cite[p.~252]{farb-margalit}), there exists a morphism $\PMod_{0,n}^\mathrm{ab}\to A_\gamma$ that makes the following diagram commute and induces the morphism $\sigma$.
            \begin{center}
                \begin{tikzcd}
                    \PMod_{0,n} \arrow[r, "\eta_\rho"] \arrow[d] \arrow[dr, dashed, "\sigma"] & A_\gamma/Z_\rho \\
                    \PMod_{0,n}^\mathrm{ab} \arrow[r, dashed, "\exists"] & A_\gamma \arrow[u]
                \end{tikzcd}
            \end{center}
        \end{proofclaim}

        Let now $\rho\colon\pi_1 S_{0,n}\to \PU(2,1)$ be a fixed point. If $\rho$ is abelian or if $Z_\rho=\{1\}$, then $\rho$ is the restriction of a mapping class group representation by Lemma~\ref{lem:regular-and-abelian-fixed-points-are-restrictions}. We now consider the case where $\rho$ is not abelian and $Z_\rho\neq \{1\}$. By the classification of maximal subgroups of $\PU(2,1)$ with non-trivial centralizer (Lemma~\ref{lem:maximal-subgroups-of-PU(2,1)-with-non-trivial-centralizer}), $\Gamma=\Im(\rho)$ is contained in $\stab(X)$, $\stab(L)$, or $\stab(H)$, where $X\in\CH^2$ is a point, $L\subset \CH^2$ is a $\C$-line, and $H$ is a horosphere foliation based at $Y\in\partial \CH^2$. We consider each case separately and show that $\eta_\rho$ always lifts.
        \begin{enumerate}
            \item $\Gamma\subset\stab(X)$: First note that $\Gamma$ fixes no other points $X'\in\CH^2$, for otherwise it would fix the whole $\C$-line through $X$ and $X'$, meaning that $\Gamma\subset \mathcal{CR}_{XX'}$, which is impossible because $\rho$ is not abelian. Next, note that $N_\rho\subset \stab(X)$, since if $g\in N_\rho$, then $(g^{-1}\Gamma g)X=\Gamma X=X$ forcing $gX=X$. Since $\rho$ is not abelian, at least one peripheral curve $\gamma\in\pi_1 S_{0,n}$ satisfies $\rho(\gamma)\notin \mathcal{SR}_X$. However, recall from Table~\ref{tab:conjugacy-classes} that $\stab(X)=Z(h)$, where $h$ is any non-trivial element of $\mathcal{SR}_{X}$. By Corollary~\ref{cor:non-abelian-intersection-centralizers}, it follows that $Z(h)\cap Z(\rho(\gamma))$ is abelian, hence $A_\gamma=N_\rho\cap Z(\rho(\gamma))\subset \stab(X)\cap Z(\rho(\gamma))$ is abelian too and Claim~\ref{claim:A_gamma-abelian-is-enough} applies.
            \item $\Gamma\subset\stab(L)$: First note that $\Gamma$ does not stabilize any other $\C$-line as otherwise it would be abelian by Lemma~\ref{lem:intersection-sstabilizers-two-C-line-is-abelian}. As in the previous case, it follows that $N_\rho\subset \stab(L)$. Since $\rho$ is not abelian, there exists a peripheral curve $\gamma\in\pi_1 S_{0,n}$ such that $\rho(\gamma)\notin \mathcal{CR}_L$. However, recall from Table~\ref{tab:conjugacy-classes} that $\stab(L)=Z(h)$ for any non-trivial $h\in \mathcal{CR}_L$. By Corollary~\ref{cor:non-abelian-intersection-centralizers}, it follows that $Z(h)\cap Z(\rho(\gamma))$ is abelian, and thus $A_\gamma$ is also abelian and Claim~\ref{claim:A_gamma-abelian-is-enough} applies.
            \item $\Gamma\subset \stab(H)$: First note that $\Gamma$ does not fix another point $Y'\subset \partial\CH^2$, for otherwise $\Gamma$ would be abelian. As in the first case, it follows that $N_\rho\subset \stab(Y)$. Since $Z(\stab(H))=\mathcal{VT}_Y$, it follows that $\mathcal{VT}_Y\subset Z_\rho$. Furthermore, for every non-trivial $z\in Z_\rho$, we have $\Gamma\subset Z(z)$. Since $Z_\rho\subset \stab(Y)$, the following cases may occur:
            \begin{itemize}
                \item $z$ is hyperbolic, regular parabolic, or screw parabolic: then $Z(z)$ is abelian (according to Table~\ref{tab:conjugacy-classes}), which is impossible because $\Gamma$ is not abelian.
                \item $z$ is a complex reflection: in that case, $Z(z)=\stab(L)$ and thus $\Gamma\subset \stab(L)$, so it reduces to the previous case.
                \item $z$ is a vertical translation fixing $Y$.
            \end{itemize}
            Excluding the first two eventualities, we conclude that $Z_\rho =\mathcal{VT}_Y$.
            
            Since $\rho$ is not abelian, there exists a peripheral curve $\gamma\in\pi_1S_{0,n}$ such that $\rho(\gamma)\notin \mathcal{VT}_Y$. If $\rho(\gamma)$ is regular parabolic or screw parabolic, then $Z(\rho(\gamma))$ is abelian (according to Table~\ref{tab:conjugacy-classes}), so $A_\gamma$ is abelian and Claim~\ref{claim:A_gamma-abelian-is-enough} applies. If not, then $\rho(\gamma)$ is a non-trivial complex reflection with mirror $L\ni Y$ and $Z(\rho(\gamma))=\stab(L)$. In that case, we have $A_\gamma\subset \stab(Y)\cap\stab(L) = \mathcal{VT}_Y\rtimes (\mathcal{CR}_L\times \R)$ (by Lemma~\ref{lem:intersection-sstabilizers-one-C-line-and-one-bdry-point}) and, since $Z_\rho = \mathcal{VT}_Y\subset A_\gamma$, we deduce that $A_\gamma/Z_\rho \cong A_\gamma\cap (\mathcal{CR}_L\times \R_{>0})$. This last identification means that $\eta_\rho$ lifts, as desired.\qedhere
        \end{enumerate}
\end{proof}

\subsection{Elementary examples of fixed points}
In this section, $\PU(2,1)$ may be replaced by any group $G$. However, it is important that the surface is a punctured sphere (and not a higher genus surface).
\begin{defn}\label{defn:pants-like}
    We say that a representation $\rho\colon\pi_1S_{0,n}\to G$ is \emph{pants-like} if $\rho$ maps the peripheral loops around at least $n-3$ punctures of $S_{0,n}$ inside $Z_\rho$, where $Z_\rho$ denotes the centralizer of the image of $\rho$ in $G$.
\end{defn}

\begin{lem}
Conjugacy classes of pants-like representations are fixed by $\PMod(S_{0,n})$.
\end{lem}
\begin{proof}
We fix an auxiliary presentation $\pi_1S_{0,n}=\langle c_1,\ldots,c_n : c_1\cdots c_n=1\rangle$. For any set $\mathcal I$ of at least two and at most $n-2$ consecutive integers in $\{1,\ldots, n-1\}$, we write $c_\mathcal{I}=\prod_{i\in\mathcal{I}} c_i$, where indices are taken in increasing order. We also denote by $f_\mathcal{I}$ the automorphism of $\pi_1\Sigma$ defined by
\[
f_{\mathcal{I}} (c_i) =\begin{cases}
   c_\mathcal{I}c_ic_{\mathcal{I}}^{-1} & \text{if } i\in \mathcal{I},\\
    c_i & \text{if } i\notin \mathcal{I}.
\end{cases}
\]
Let $\mathcal{I}^c=\{1,\ldots,n\}\setminus \mathcal{I}$ and write $f_{\mathcal{I}^c}=c_{\mathcal{I}}^{-1}f_\mathcal{I}c_{\mathcal{I}}$. Note that $c_\mathcal{I}^{-1}$ can be expressed as a product of $c_i^{-1}$ for $i\in \mathcal{I}^c$.

Now, if $\rho\colon\pi_1S_{0,n}\to G$ is pants-like, then it maps at least $n-3$ of the generators $c_1,\ldots,c_n$ inside $Z_\rho$. Without loss of generality, assume that $\rho(c_4),\ldots,\rho(c_n)\in Z_\rho$. By the pigeonhole principle, one of the two sets $\mathcal{I}\cap\{1,2,3\}$ and $\mathcal{I}^c\cap\{1,2,3\}$ contains at most one element. Assume first that $\mathcal{I}\cap\{1,2,3\}$ contains at most one element. If $\mathcal{I}\cap\{1,2,3\}=\emptyset$, then we immediately see that $\rho\circ f_\mathcal{I}=\rho$. Similarly, if $\mathcal{I}\cap\{1,2,3\}=\{j\}$, then clearly $\rho\circ f_\mathcal{I}(c_i)=c_i$ for $i\neq j$ and $\rho\circ f_{\mathcal{I}}(c_j)=c_\mathcal{I}c_jc_\mathcal{I}^{-1}=c_jc_jc_j^{-1}=c_j$, implying $\rho\circ f_\mathcal{I}=\rho$ as well. If instead $\mathcal{I}^c\cap\{1,2,3\}$ contains at most one element, then a similar argument shows that $\rho\circ f_{\mathcal{I}^c}=\rho$. Since $\PMod(S_{0,n})$ is generated by the outer automorphisms associated to all $f_\mathcal{I}$ (for $\mathcal{I}$ as above), we conclude that $\PMod(S_{0,n})$ fixes $[\rho]$.
\end{proof}

\begin{rem}\label{rem:pants-like-rank-2}
When $G=\mathrm{GL}_2\C$, then pants-like representations can only be of two forms: either they are abelian, or the peripheral loops around at least $n-3$ punctures of $S_{0,n}$ are mapped into $Z(\mathrm{GL}_2\C)$. Moreover, these representations account for all mapping class group fixed points in rank $2$~\cite[Lemma~3.3]{bronstein-maret}. However, when $G=\mathrm{GL}_n\C$ with $n\geq 3$, or when $G=\PU(2,1)$, then there are pants-like representations that are not abelian and which do not map any peripheral curve in the centre of $G$.
\end{rem}

\subsection{Braid group representations}\label{sec:braid-group-representations}
Representations of the pure braid group $B_n$ into $\PU(2,1)$ that vanish on $Z(B_n)$ give rise to mapping class group fixed points by Lemma~\ref{lem:restriction-of-PMod-representations-are-fixed-points}. In this section we recall the classical examples of such representations given by Burau representations. We start with one particularly striking example that initially sparked our curiosity to investigate fixed points for $\PU(2,1)$ representations.

\subsubsection{The Casella fixed points}
 Interestingly, the first non-elementary fixed point we encountered comes from a representation $\varphi\colon \pi_1 S_{1,1}\to \PU(2,1)$ of the once-punctured torus $S_{1,1}$. It was originally studied by Casella and exhibited as the one with all eight Fock--Goncharov coordinates equal to the third root of unity $\omega=e^{i2\pi/3}$ in the $\mathrm{PGL}_3\C$ character variety of $S_{1,1}$. According to Casella, $\varphi$ and its conjugate $\overline{\varphi}$ are the only two non-elementary fixed points in the rank-3 decorated character variety of the once-punctured torus~\cite[Lemma~13 and the discussion thereafter]{casella-phd}. The explicit form of $\varphi$ can be found in~\cite[Section~6.2 and pages 11--13]{casella-article}. Writing $\pi_1S_{1,1}=\langle a,b\rangle$, we have
\begin{equation}\label{eq:casella-orbit-torus}
\varphi(a)=\begin{pmatrix}
    1-\omega & 0 & -\omega^2\\
    0 & \omega^2 & 0\\
    -\omega^2 & 0 & 0    
\end{pmatrix}, \quad
\varphi(b) = \begin{pmatrix}
    -\omega^2 & -\omega^2 & -1\\
    \omega^2 & 0 & -\omega\\
    -\omega & 1 & -\omega^2
\end{pmatrix}.
\end{equation}
Both $\varphi(a)$ and $\varphi(b)$ preserve the anti-diagonal Hermitian form $J'$ from~\eqref{eq:hermitian-form-antidiagonal}. This means that $\varphi$ is valued in $\PU(J';\Z[\omega])$~\cite[Lemma~15]{casella-phd}. In particular, $\varphi$ is a discrete representation. Casella also proved that $\varphi$ is irreducible (actually even Zariski dense) and has infinite image~\cite[Lemma~14]{casella-phd}.

The following procedure turns $\varphi$ into a representation of a 4-punctured sphere instead of a once-punctured torus. First, since the conjugacy class of $\varphi\colon \pi_1 S_{1,1}\to\PU(2,1)$ is fixed by $\PMod(S_{1,1})$ and $\varphi$ is irreducible (hence regular), the representation $\varphi$ factorizes through the point-pushing morphism $\pi_1 S_{1,1}\hookrightarrow \PMod(S_{1,2})$ by an analogue of Lemma~\ref{lem:regular-and-abelian-fixed-points-are-restrictions}.
\[
\begin{tikzcd}
\pi_1S_{1,1} \arrow[r, "\varphi"] \arrow[d, hook] & \PU(2,1)\\
\PMod(S_{1,2}) \arrow[ur, dashed, start anchor=east, "\exists\, \widehat{\varphi}"]
\end{tikzcd}
\]
The group $\PMod(S_{1,2})$ is generated by three Dehn twists $\tau_1,\tau_2,\tau_3$. The embedding $\PMod(S_{1,2})\subset \Aut(\pi_1S_{1,1})$ maps $\tau_1,\tau_2,\tau_3$ to the following automorphisms of $\pi_1 S_{1,1}=\langle a,b\rangle$:
\begin{equation}\label{eq:embedding-PMod_1,2-in-Aut(pi_1S_1,1)}
\tau_1\mapsto \begin{cases}
     a \mapsto a  \\
     b \mapsto ba \end{cases}, \quad
\tau_2\mapsto \begin{cases}
     a \mapsto ab^{-1}\\
     b \mapsto b \end{cases}, \quad
\tau_3\mapsto \begin{cases}
     a \mapsto a  \\
     b \mapsto ab. \end{cases}
\end{equation}
The representation $\widehat\varphi$ maps $\tau_i$ to the unique element $g_i$ of $\PU(2,1)$ such that $\varphi\circ\tau_i= g_i\varphi g_i^{-1}$. The elements $g_i$ exist because $[\varphi]$ is fixed by $\PMod(S_{1,1})$, and they are unique because $\varphi$ is irreducible and $\PU(2,1)$ is centre-free. One way to determine each $g_i$, is to start by diagonalizing $\varphi(a)$ or $\varphi(b)$ and then solve for $g_i$. A direct computation shows that
\[
g_1 = \begin{pmatrix}
    1 & 0 & 0\\
    0 & -\omega^2 & 0\\
    0 & 0 & 1
\end{pmatrix}, \quad g_2 = \begin{pmatrix}
    \omega & \omega & \omega^2\\
    0 & -1 & 1\\
    0 & 0 & \omega
\end{pmatrix}, \quad g_3 = \begin{pmatrix}
    \omega + 2 & 0 & \omega\\
    0 & 1 & 0\\
    \omega & 0 & 0
\end{pmatrix}.
\]

There is a somewhat exceptional isomorphism between $\PMod(S_{1,2})$ and the braid group on four strands $B_4$ quotient by its centre:
\begin{equation}\label{eq:PMod_1,2=B_4/Z(B_4)}
    \PMod(S_{1,2})\cong B_4/Z(B_4).
\end{equation}
This isomorphism is explained by the fact that the generators $\tau_1,\tau_2,\tau_3$ satisfy precisely the same relations as Artin's generators $\sigma_1,\sigma_2,\sigma_3$ of $B_4$, see~e.g.~\cite[discussion after Theorem~9.2]{farb-margalit}.

On the other hand, if we consider the pure braid group $PB_4$, then there is a canonical isomoprhism between $PB_4/Z(B_4)$ and $\PMod(S_{0,5})$. Recall that the fundamental group of a $4$-punctured sphere $\pi_1 S_{0,4}$ embeds in $\PMod(S_{0,5})$ by the point-pushing morphism. More concretely, if we fix a presentation of $\pi_1 S_{0,4}$ of the following kind:
\[
\pi_1 S_{0,4} = \langle c_1,c_2,c_3, c_4 : \prod c_i =1 \rangle,
\]
then the composition of the point-pushing morphism $\pi_1 S_{0,4}\to\PMod(S_{0,5})$ and the inclusion $\PMod(S_{0,5})\subset B_4/Z(B_4)$ has the explicit form
\begin{align}
    \pi_1 S_{0,4}&\to B_4/Z(B_4)\label{eq:morphism-pi_1S_0,4-to-B_4/Z(B_4)}\\
    c_1 &\mapsto \big(\sigma_2^{-1}\sigma_3^{-1}\big)^3\nonumber\\
    c_2 &\mapsto \big(\sigma_1^{-1}\sigma_2^{-1}\sigma_1\sigma_3^{-1}\big)^3\nonumber\\
    c_3 &\mapsto \big(\sigma_1^{-1}\sigma_2^{-1}\sigma_3^{-1}\sigma_2\big)^3\nonumber\\
    c_4 &\mapsto \big(\sigma_1^{-1}\sigma_2^{-1}\big)^3.\nonumber
\end{align}
The composition of the map $\pi_1S_{0,4}\to B_4/Z(B_4)$ with $\widehat{\varphi}\colon \PMod(S_{1,2})\to\PU(2,1)$ leads to a representation $\rho\colon\pi_1S_{0,4}\to \PU(2,1)$, explicitly given by
\begin{align}
\rho(c_1)= \begin{pmatrix}
    -\omega-1 & -\omega+1 & \omega+2\\
    -\omega+1 & 2\omega+2 & 2\omega+1\\
    \omega+2 & 2\omega+1 & 2\omega-1\\
\end{pmatrix}, \quad&
\rho(c_2)= \begin{pmatrix}
    \omega & -2\omega-1 & -2\omega-1\\
    \omega-1 & -2\omega & -2\omega-1\\
    -2\omega-1 & \omega+2 & \omega+3\\
\end{pmatrix}, \label{eq:Casella-orbit-antidiagonal} \\
\rho(c_3)= \begin{pmatrix}
    -2 & 0 & -6\omega -3\\
    0 & 1 & 0\\
    -2\omega-1 & 0 & 4\\
\end{pmatrix}, \quad&
\rho(c_4)= \begin{pmatrix}
    1 & 0 & -2\omega-1\\
    0 & 1 & 0\\
    0 & 0 & 1\\
\end{pmatrix}.\nonumber
\end{align}

\begin{defn}\label{defn:Casella-fixed-point}
    Any representation $\pi_1 S_{0,4}\to \PU(2,1)$ conjugate to $\rho$ or to its complex conjugate $\overline{\rho}$ is called a \emph{Casella fixed point}.
\end{defn}

\begin{rem}
Recall that complex conjugation is the only non-trivial outer automorphism of $\PU(2,1)$ (it is conjugate to the inverse transpose by definition of $\PU(2,1)$). On the other hand, $\PMod(S_{1,2})$ also has only one non-trivial outer automorphism $\imath$ given by $\imath(\tau_i)= \tau_i^{-1}$. We claim that $\widehat{\varphi}\circ\imath$ is conjugate to the complex conjugate of $\widehat{\varphi}$. Since $(\widehat{\varphi}\circ\imath)(\tau_i)=\widehat{\varphi}(\tau_i)^{-1}=g_i^{-1}$, it is sufficient to prove that $g_1^{-1},g_2^{-1},g_3^{-1}$ can be simultaneously conjugated into $\overline{g_1},\overline{g_2},\overline{g_3}$. (Recall that $J'g_i^{-1}J'=\overline{g_i}^t$, but since $g_2$ is not symmetric, the matrix $J'$ does not do job here.) This is achieved by the matrix
\[
\begin{pmatrix}
    1 & 0 & 2\omega-1\\
    0 & \omega & 0\\
    0& 0 &1
\end{pmatrix}.
\]
\end{rem}

\begin{lem}\label{lem:properties-Casella-fixed-point}
The Casella fixed point $\rho$ has the following properties:
\begin{enumerate}
    \item The image of $\rho$ is contained in $\PU(J';\Z[\omega])$. In particular, $\rho$ is discrete.
    \item The conjugacy class $[\rho]$ is fixed by $\PMod(S_{0,4})$.
    \item All non-peripheral simple closed curves in $S_{0,4}$ are mapped by $\rho$ to regular elliptic elements.
    \item Each $\rho(c_i)$ is a vertical translation fixing a different point on $\partial\CH^2$. No three of the fixed points are contained in the boundary of the same $\C$-line.
    \item The image of $\rho$ is a Zariski dense subgroup of $\PU(2,1)$ (in particular, infinite).   
\end{enumerate}
\end{lem}
\begin{proof}
The first property is inherited from $\varphi$. It can also be checked directly from~\eqref{eq:Casella-orbit-antidiagonal}. Since $\rho$ was constructed as the restriction of a representation of $\PMod(S_{0,5})$, it follows that $[\rho]$ is fixed by $\PMod(S_{0,4})$ by Lemma~\ref{lem:restriction-of-PMod-representations-are-fixed-points}, proving the second property. To prove the third property, it is sufficient to check that $\rho(c_1c_2)$, $\rho(c_2c_3)$, and $\rho(c_1c_3)$ are regular elliptic, because $[\rho]$ is fixed by $\PMod(S_{0,4})$. This can be established by direct computations.

Now, let us consider the fourth property. The trace of each $\rho(c_i)$ belongs to $\{3,3\omega,3\omega^2\}$. Furthermore, each $\rho(c_i)$ is represented by a matrix of determinant $1$ in~\eqref{eq:Casella-orbit-antidiagonal} with minimal polynomial of degree $2$. This shows that each $\rho(c_i)$ is a vertical translation. In order to understand the fixed points of $\rho(c_1),\ldots,\rho(c_4)$, it is convenient to conjugate $\rho$ so that it preserves the standard Hermitian form $J$ from~\eqref{eq:hermitian-form-diagonal}, rather than the anti-diagonal form $J'$, and so that the unique fixed point of $\rho(c_1c_2)$ is the origin in the ball model. Doing so, we obtain the representation $\rho'\colon\pi_1S_{0,4}\to \PU(2,1)$ given by
\begin{align}
\rho'(c_1)= \begin{pmatrix}
    1 & 0 & 0\\
    0 & -2\omega & -\sqrt{3}\\
    0 & -\sqrt{3} & -2\omega^2\\
\end{pmatrix}, \quad&
\rho'(c_2)= \begin{pmatrix}
    1 & 0 & 0\\
    0 & -2\omega & \sqrt{3}\omega\\
    0 & \sqrt{3}\omega^2 & -2\omega^2\\
\end{pmatrix}, \label{eq:Casella-fixed-point-diagonal}\\
\rho'(c_3)= \begin{pmatrix}
    -2\omega & 0 & 2\omega + 1\\
    0 & 1 & 0\\
    -2\omega-1 & 0 & -2\omega^2\\
\end{pmatrix}, \quad&
\rho'(c_4)= \begin{pmatrix}
    -2\omega & 0 & \omega + 2\\
    0 & 1 & 0\\
    -\omega+1 & 0 & -2\omega^2\\
\end{pmatrix}.\nonumber
\end{align}
It can be checked directly that the fixed points of the images of the generators are the following points in $\partial \CH^2$ (illustrated on Figure~\ref{fig:casella-fixed-point}):
\[
Y_1=(0,i),\quad Y_2=\big(0,e^{\frac{i\pi}{6}}\big),\quad Y_3=(1,0), \quad Y_4=(-\omega,0).
\]
The fifth property follows directly from the fourth one and the fact that the four fixed points are not contained in the boundary of the same $\C$-line.  
\end{proof}
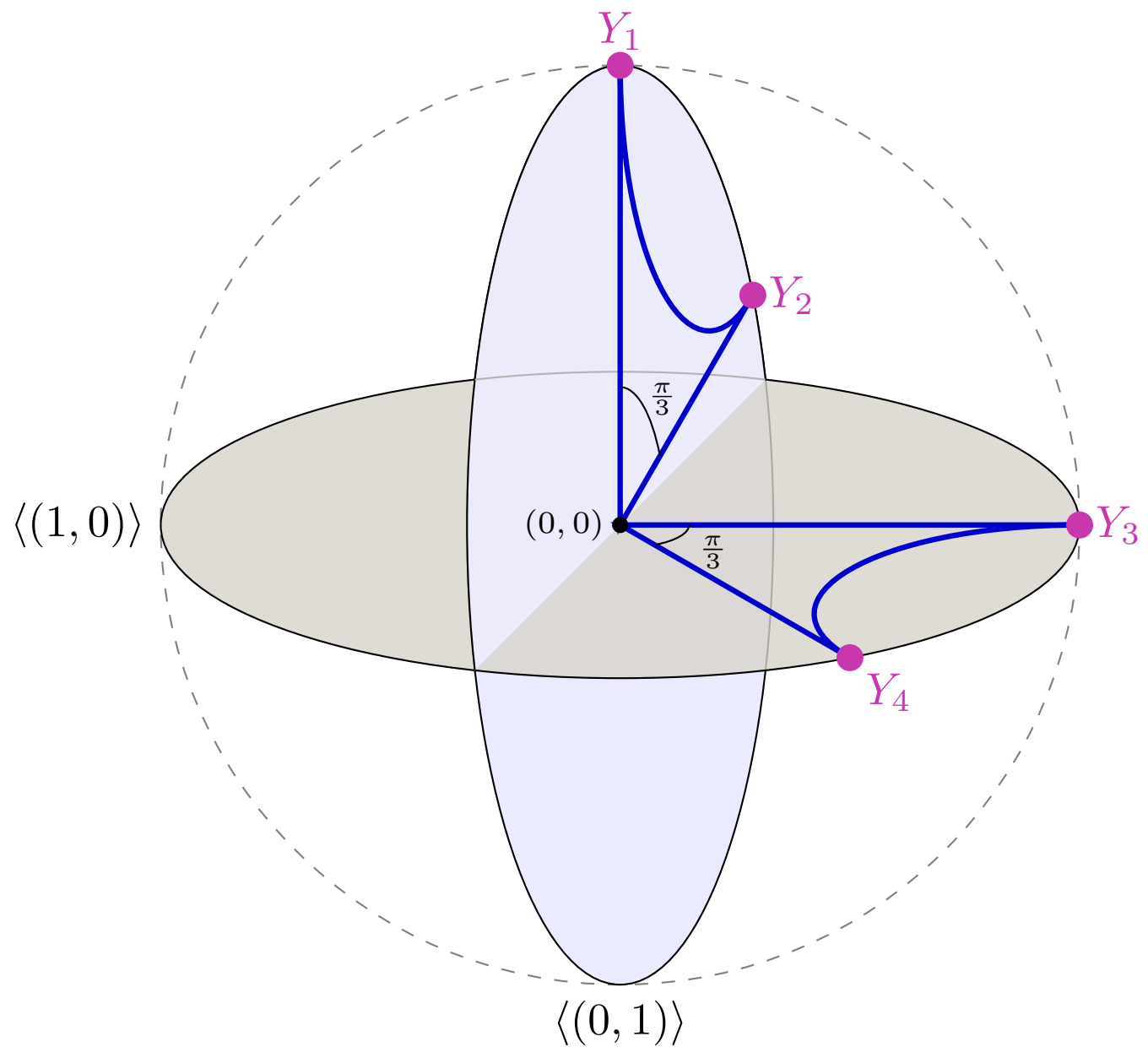
\begin{figure}[ht]
    \centering
    \begin{tikzpicture}[scale=1.2]
  \draw[gray, dashed] (0,0) circle (3);

  \draw[fill=blue!8] (0,0) ellipse [x radius=1, y radius=3];
  \draw[fill=smoked, fill opacity=0.8] (0,0) ellipse [x radius=3, y radius=1];

  \begin{scope}
    \clip (-4,-4) -- (4,4) -- (-4,4) -- cycle;
    \draw[fill=blue!80!black!8] (0,0) ellipse [x radius=1, y radius=3];
    \clip (0,0) ellipse [x radius=1, y radius=3];
    \draw[gray!60] (0,0) ellipse [x radius=3, y radius=1];
  \end{scope}

  \draw[blue!80!black, line width=1.2pt]
    (0,0) -- (90:1 and 3)
    .. controls (90:0.56 and 1.67) and (30:0.56 and 1.67) .. (30:1 and 3)
    -- cycle;

  \draw[blue!80!black, line width=1.2pt]
     (0,0) -- (0:3 and 1)
     .. controls (0:1.67 and 0.56) and (-60:1.67 and 0.56) .. (-60:3 and 1)
     -- cycle;

  \draw (30:0.3 and 0.9) arc[start angle=30, end angle=90, x radius=0.3, y radius=0.9];
  \node[font=\tiny] at (60:0.55 and 0.95) {$\frac{\pi}{3}$};
  \draw (-60:0.45 and 0.15) arc[start angle=-60, end angle=0, x radius=0.45, y radius=0.15];
   \node[font=\scriptsize] at (-30:.7 and 0.35) {$\frac{\pi}{3}$};

  \fill (0,0) circle (1.5pt);
  \node[left, font=\scriptsize] at (0,0) {$(0,0)$};

  \fill[mauve] (90:1 and 3)  circle (2.5pt) node[above] {$Y_1$};
  \fill[mauve] (30:1 and 3)  circle (2.5pt) node[right] {$Y_2$};
  \fill[mauve] (0:3 and 1)   circle (2.5pt) node[right] {$Y_3$};
  \fill[mauve] (-60:3 and 1) circle (2.5pt) node[below right] {$Y_4$};

  \node[left]  at (180:3 and 0) {$\langle(1,0)\rangle$};
  \node[below] at (270:1 and 3) {$\langle(0,1)\rangle$};
\end{tikzpicture}
    \caption{An illustration of the Casella fixed point~\eqref{eq:Casella-fixed-point-diagonal} with the four fixed points $Y_i$ of $\rho'(c_i)$ and the two $\C$-lines invariant by $\rho'(c_1),\rho'(c_2)$ (in blue), respectively $\rho'(c_3),\rho'(c_4)$ (in gray). The two $\C$-lines intersect at $(0,0)$ which is the unique fixed point of the regular elliptic $\rho'(c_1c_2)=\rho'(c_3c_4)^{-1}$.}
    \label{fig:casella-fixed-point}
\end{figure}

\subsubsection{The Burau fixed points}\label{sec:Burau-fixed-points-sphere}
Casella's orbit belongs to a larger family of fixed points coming from the Burau representations of braid groups.

The \emph{Burau representation} of the braid group $B_n$ is a morphism $B_n\to \mathrm{GL}_n\big(\Z[t,t^{-1}]\big)$, which can be reduced to a representation 
\[
\beta_n\colon B_n\to \mathrm{GL}_{n-1}\big(\Z[t,t^{-1}]\big).
\]
It turns out that $\beta_n(Z(B_n))$ is made of scalar matrices, so that $\beta_n$ descends to a representation
\[
\beta_n\colon B_n/Z(B_n)\to \mathrm{PGL}_{n-1}\big(\Z[t,t^{-1}]\big).
\]
When $n=4$, the image of Artin's generators $\sigma_1,\sigma_2,\sigma_3$ are the following three matrices (see e.g.~\cite[(3-26)]{birman}):
\begin{equation}\label{eq:Burau-matrices-4}
\beta_4(\sigma_1)=\begin{pmatrix}
 -t & 1 & 0\\
 0 & 1 & 0\\
 0 & 0 & 1
\end{pmatrix}, \quad
\beta_4(\sigma_2)=\begin{pmatrix}
 1 & 0 & 0\\
 t & -t & 1\\
 0 & 0 & 1
\end{pmatrix}, \quad
\beta_4(\sigma_3)=\begin{pmatrix}
 1 & 0 & 0\\
 0 & 1 & 0\\
 0 & t & -t
\end{pmatrix}.
\end{equation}
Squier proved that $\beta_n(B_n)$ preserves a Hermitian form on $\C^{n-1}$ after specifying $t$ as a complex number of modulus $1$~\cite{squier}. When $n=4$, the form is given by
\[
J_s=\begin{pmatrix}
    s+s^{-1} & -s^{-1} & 0\\
    -s & s+s^{-1} & -s^{-1}\\
    0 & -s & s+s^{-1}
\end{pmatrix},
\]
where $s^2=t$. The eigenvalues of $J_0$ are $s+s^{-1}-\sqrt{2}$, $s+s^{-1}$, and $s+s^{-1}+\sqrt{2}$. So, if $s=e^{i\theta}$, then $J_s$ is of the following type according to the value of $\theta$ (see also Figure~\ref{fig:signature-Burau-form} and~\cite[Proposition~5.2]{FGG26} for a different parametrization):
\smallskip
\begin{center}
\begin{tabular}{c||c|c|c}
\toprule
         $\theta$ & $\big(0, \frac{\pi}{4}\big)\cup \big (\frac{3\pi}{4},\pi\big)$ & $\big(\frac{\pi}{4},\frac{3\pi}{4}\big)\setminus\big\{\frac{\pi}{2}\big\}$ & $\big\{\frac{k\pi}{4}\big\}$ \\\midrule
         $\mathrm{sign}(J_s)$ & $(3,0)$ & $(2,1)$ & \text{singular}\\ 
\bottomrule
\end{tabular}
\end{center}
\smallskip
Note that it is convenient to replace $J_s$ by $-J_s=J_{-s}$ when $\theta>\pi/2$ in order to have more positive than negative directions.

There is an equivalent form of the Burau representation that can also be found in the literature (this is the form currently implemented in SageMath). It is given by
\begin{equation}\label{eq:Burau-reduced-n=4-sage}
\beta_4'(\sigma_1)=\begin{pmatrix}
 -s^2 & s & 0\\
 0 & 1 & 0\\
 0 & 0 & 1
\end{pmatrix}, \quad
\beta_4'(\sigma_2)=\begin{pmatrix}
 1 & 0 & 0\\
 s & -s^2 & s\\
 0 & 0 & 1
\end{pmatrix}, \quad
\beta_4'(\sigma_3)=\begin{pmatrix}
 1 & 0 & 0\\
 0 & 1 & 0\\
 0 & s & -s^2
\end{pmatrix}.
\end{equation}
Note that $\beta_4'=P\beta_4P^{-1}$, where $P=\mathrm{diag}(s,1,s^{-1})$. The Hermitian form preserved by $\beta_4'$ is
\[
J_s'=\begin{pmatrix}
    s+s^{-1} & -1 & 0\\
    -1 & s+s^{-1} & -1\\
    0 & -1 & s+s^{-1}
\end{pmatrix}.
\]

Consider the one-parameter family of Burau representations
\[
\beta_{4,t}'\colon B_4/Z(B_4)\to \PU(2,1)
\]
indexed by $t=e^{2i\theta}$ with $\theta\in (\pi/4, 3\pi/4)\setminus\{\pi/2\}$. If we pre-compose it with the morphism $\pi_1S_{0,4}\to B_4/Z(B_4)$ from~\eqref{eq:morphism-pi_1S_0,4-to-B_4/Z(B_4)}, we obtain a one-parameter family of representations
\[
\rho_t\colon \pi_1 S_{0,4}\to \PU(2,1),
\]
explicitly given by
\begin{align}
\rho_t(c_1)= \begin{pmatrix}
    1 & 0 & 0\\
    s-s^5 & s^6 & 0\\
    s^2-s^4 & 0 & s^6\\
\end{pmatrix}, \quad&
\rho_t(c_2)= \begin{pmatrix}
    s^4 & s^7-s^5 & 0\\
    s^3-s^{-1} & s^6-s^4+1 & 0\\
    s^2-1 & s-s^3 & s^6\\
\end{pmatrix}, \label{eq:Burau-representation} \\
\rho_t(c_3)= \begin{pmatrix}
    s^6 & -s^5+s^3 & s^6-s^4\\
    0 & s^2 & s^7-s^3\\
    0 & s-s^{-1} & s^6-s^2+1\\
\end{pmatrix}, \quad&
\rho_t(c_4)= \begin{pmatrix}
    s^6 & 0 & s^2-s^4\\
    0 & s^6 & s-s^5\\
    0 & 0 & 1\\
\end{pmatrix}.\nonumber
\end{align}
One can verify that $\rho_t(c_1)\rho_t(c_2)\rho_t(c_3)\rho_t(c_4)=\diag(s^{16},s^{16},s^{16})$ is indeed trivial in $\PU(2,1)$.

Since every $\rho_t$ is the restriction of a morphism of $\PMod(S_{0,5})$, the conjugacy classes $[\rho_t]$ are all fixed by $\PMod(S_{0,4})$ by Lemma~\ref{lem:restriction-of-PMod-representations-are-fixed-points}.

\begin{defn}\label{def:Burau-fixed-points}
The conjugacy classes $[\rho_t]$ are called \emph{Burau fixed points}.
\end{defn}

\begin{lem}\label{lem:properties-Burau-fixed-points}
    The Burau fixed points have the following properties. Recall that they are indexed by $t=e^{2i\theta}$ for $\theta\in (\pi/4, 3\pi/4)\setminus\{\pi/2\}$. 
    \begin{enumerate}
        \item Replacing $\theta$ by $\pi-\theta$ complex conjugates $\rho_t$, i.e.~$\rho_{\overline{t}}=\overline{\rho_t}$.
        \item In all cases, $\rho_t$ maps every non-peripheral simple closed curve on $S_{0,4}$ to regular elliptic elements.
        \item If $\theta\in\{\pi/3, 2\pi/3\}$, then $\rho_t$ is one of the two Casella representations.
        \item If $\theta\in (\pi/4,\pi/3)\cup (2\pi/3,3\pi,4)$, then $\rho_t(c_i)$ are scalar rotations.
        \item If $\theta\in (\pi/3,2\pi/3)\setminus \{\pi/2\}$, then $\rho_t(c_i)$ are complex reflections.
        
    \end{enumerate}
\end{lem}
\begin{figure}[ht]
    \centering
    \begin{tikzpicture}
    \draw (-3.8,0) -- (3.8,0);

    \draw[smoked, line width=2pt] (0:3) arc (0:45:3);
    \draw[smoked, line width=2pt] (135:3) arc (135:180:3);

    \draw[blue!60!black, line width=2pt] (45:3) arc (45:60:3);
    \draw[blue!60!black, line width=2pt] (120:3) arc (120:135:3);

    \draw[sky, line width=2pt] (60:3) arc (60:120:3);

    \fill[mauve] (60:3) circle (3pt);
    \fill[mauve] (120:3) circle (3pt);

    \foreach \a in {45,90,135}
    \draw[fill=white] (\a:3) circle (3pt);

    \foreach \a/\l in {45/{\pi/4}, 60/{\pi/3}, 90/{\pi/2}, 120/{2\pi/3}, 135/{3\pi/4}}
    \node at (\a:2.5) {\small $\l$};

    \node[right, smoked] at (20:3.1) {$\PU(3)$};
    \node[left, smoked]  at (160:3.1) {$\PU(3)$};
    \node[above right, blue!80!black] at (67:3.1) {$\PU(2,1)$};
    \node[above left, blue!80!black] at (112:3.1) {$\PU(2,1)$};
    \end{tikzpicture}
    \caption{In dark blue: the Burau fixed points where peripheral monodromies are scalar rotations. In light blue : the Burau fixed points where peripheral monodromies are complex reflections. In mauve : the transition regime given by the Casella representations (peripheral monodromies are vertical translations). The white points are the values of degeneracy of the Hermitian form.}
    \label{fig:signature-Burau-form}
\end{figure}
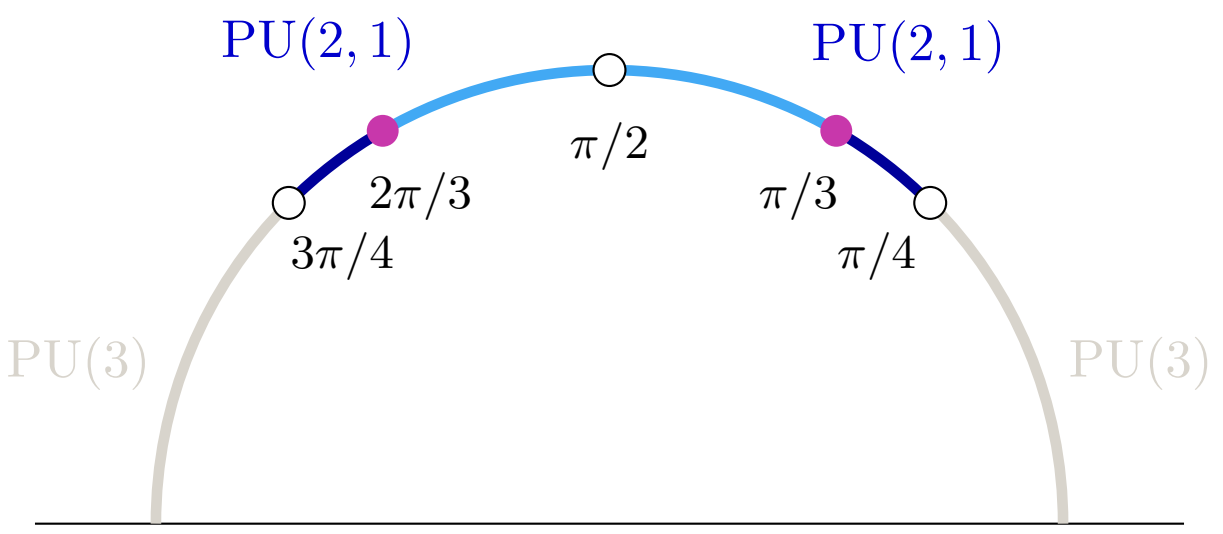
\begin{proof}
    If we take the four matrices from~\eqref{eq:Burau-representation} and  multiply them by $s^{-4}$, we lift $\rho_t$ into a representation $\widehat{\rho}_t\colon\pi_1S_{0,4}\to\SU(2,1)$. The eigenvalues of $\widehat{\rho}_t(c_i)$ are $\{s^{-4}, s^2, s^2\}$. Hence, the trace of $\widehat{\rho}_t(c_i)$ equals $2s^2+s^{-4}=2e^{2i\theta}+e^{-4i\theta}$, which corresponds to a parametrization of Goldman's deltoid~\cite[Section~6.3.2 and Figure~6.1]{goldman-book}, see also Figure~\ref{fig:Goldman's-deltoid}.
    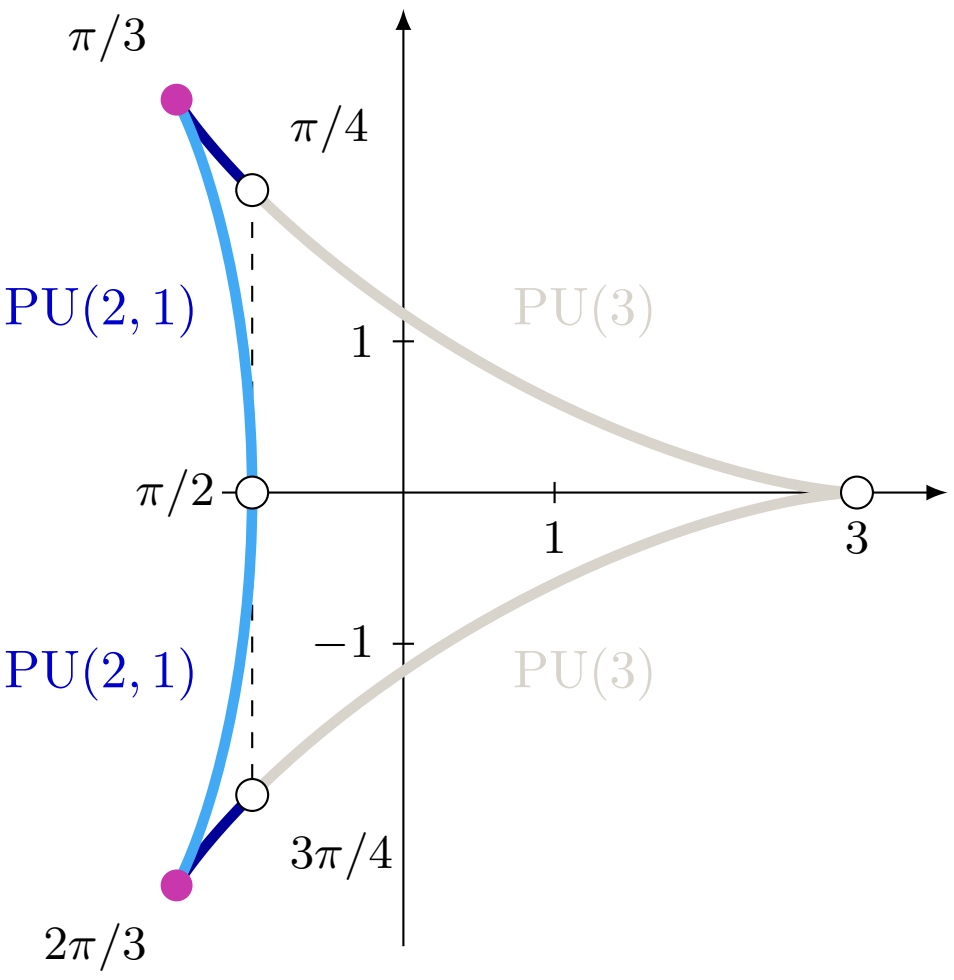
\begin{figure}[ht]
    \centering
        \begin{tikzpicture}
        \draw[->] (-1.2,0) -- (3.6,0);
        \draw[->] (0,-3) -- (0,3.2);
        \foreach \x in {1,3} \draw (\x,2pt) -- (\x,-2pt) node[below] {\small $\x$};
        \foreach \y in {-1,1} \draw (2pt,\y) -- (-2pt,\y) node[left] {\small $\y$};
        \draw[dashed] (-1,2) -- (-1,-2);
        
        \foreach \a/\b/\c in {0/45/smoked, 135/180/smoked, 45/60/blue!60!black, 120/135/blue!60!black, 60/120/sky} \draw[\c, line width=2pt, domain=\a:\b, samples=40, variable=\x] plot ({2*cos(2*\x)+cos(4*\x)}, {2*sin(2*\x)-sin(4*\x)});

        \foreach \a/\l/\d in {45/{\pi/4}/45, 60/{\pi/3}/120, 90/{\pi/2}/180, 120/{2\pi/3}/240, 135/{3\pi/4}/-45}
        \path ({2*cos(2*\a)+cos(4*\a)}, {2*sin(2*\a)-sin(4*\a)}) node[label={\d:\small $\l$}] {};

        \node[smoked] at (1.2,1.2) {$\PU(3)$};
        \node[smoked] at (1.2,-1.2) {$\PU(3)$};
        \node[blue!80!black] at (-2,1.2) {$\PU(2,1)$};
        \node[blue!80!black] at (-2,-1.2) {$\PU(2,1)$};
        
        \foreach \p in {(120:3), (240:3)} \fill[mauve] \p circle (3pt);
        \foreach \p in {(-1,2), (-1,0), (-1,-2), (3,0)} \draw[fill=white] \p circle (3pt);
        \end{tikzpicture}
    \caption{The peripheral traces of $\widehat{\rho}_t$ drawn on Goldamn's deltoid according to the value of $\theta\in [0,\pi)$, with the regime $\theta\in (\pi/4,3\pi/4)\setminus\{\pi/2\}$ corresponding to $\PU(2,1)$ representations.}
    \label{fig:Goldman's-deltoid}
    \end{figure}

    Property (1) is immediate. To prove (2), not that, since $\rho_t$ is a fixed point, it is sufficient to check that $\widehat{\rho}_t(c_1c_2), \widehat{\rho}_t(c_2c_3), \widehat{\rho}_t(c_1c_3)$ are regular elliptic elements. A direct computation shows that their traces equal $1+s^{-4}+s^4=1+2\cos(4\theta)$, which, by Goldman's criterion~\cite[Theorem~6.2.4]{goldman-book}, lies in the interior of the deltoid for all $\theta\neq k\pi/4$. This proves that the images of non-peripheral simple closed curves are regular elliptic when $\theta\neq k\pi/4$.

    We move to (3). When $\theta=2\pi/3$, then $t=\omega^2$ and $s=\omega$ are primitive cubic roots of unity. To see the relation with Casella's fixed point, we conjugate $\rho_{\omega^2}$ so that it preserves the diagonal Hermitian form $J$ from~\eqref{eq:hermitian-form-diagonal}. Consider the following matrices:
    \[
    P=\begin{pmatrix}
        1 & 0 & -1 \\ 0 & 0 & 1 \\ 0 & 1 & -1
    \end{pmatrix}, \quad Q=\begin{pmatrix}
        1 & 0 & 0\\
        0 & i\omega^2 & 0\\
        0 & 0 & 1
    \end{pmatrix}.
    \]
    The representation $P^{-1}\rho_{\omega^2}P$ takes image in $\PU(J)$. If we further conjugate by $Q$, then the representation $Q^{-1}P^{-1}\rho_{\omega^2}PQ$ can be checked to precisely coincide with the complex conjugate of the Casella representation given by~\eqref{eq:Casella-fixed-point-diagonal} (up to a cyclic permutation of the punctures).

    To prove (4) and (5), it is not sufficient to know the peripheral traces of $\widehat{\rho}_t$ as they do not determine a single conjugacy class (traces on the deltoid away from its vertices correspond to complex reflections, vertical translations, or screw parabolic). Instead, we have to determine if the images of the generators fix a whole $\C$-line, or just a point in $\CH^2$. For instance, consider
    \[
    \widehat{\rho}_t(c_4)=\begin{pmatrix}
        s^2 & 0 & s^{-2}-1\\
        0 & s^2 & s^{-3}-s\\
        0 & 0 & s^{-4}
    \end{pmatrix} \in\SU(J_s').
    \]
    The eigenspace $E_{s^2}\subset \C^3$ is simply $\C e_1\oplus \C e_2$, while the eigenspace $E_{s^{-4}}$ is equal to $\C v$ for $v=(s^{-2}-1, s^{-3}-1, s^{-4}-s^2)$. We compute that $\langle v,v\rangle_{J_s'}=4\sin(4\theta)\sin(3\theta)$. Moreover, the restriction of $J_s'$ to $E_{s^2}$ is
    \[
    J_s'\vert_{E_{s^2}} = \begin{pmatrix}
        s^{-1}+s & -1\\ -1 & s^{-1}+s
    \end{pmatrix}.
    \]
    It has trace $4\cos(\theta)$ and determinant $2\cos(2\theta)+1$. In conclusion, we have the following sign table.
    \smallskip
    \begin{center}
    \begin{tabular}{c||c|c|c|c|c|c|c}
    \toprule
         $\theta$ & $\big(\frac{\pi}{4},\frac{\pi}{3}\big)$ & $\frac{\pi}{3}$ & $\big(\frac{\pi}{3},\frac{\pi}{2}\big)$ & $\frac{\pi}{2}$ & $\big(\frac{\pi}{2},\frac{2\pi}{3}\big)$ & $\frac{2\pi}{3}$ & $\big (\frac{2\pi}{3},\frac{3\pi}{4}\big)$ \\\midrule
         $\mathrm{sign}(J_s'\vert_{E_{s^2}})$ & $(2,0)$ & singular & $(1,1)$ & $(1,1)$ & $(1,1)$ & singular & $(0,2)$\\ 
         $\mathrm{sign}(\langle v,v\rangle_{J_s'})$ & $-$ & $0$ & $+$ & $0$ & $-$ & $0$ & $+$\\
    \bottomrule
    \end{tabular}
    \end{center}
    \smallskip
    Note that, in order to have two positive directions and one negative one, we should replace $J_s'$ by $-J_s'$ in the regime $\theta\in (\pi/2,3\pi/4)$, After this adjustment, we note that $E_{s^2}$ projects to a fixed $\C$-line when $\theta\in (\pi/3,\pi/2)\cup (\pi/2,2\pi/3)$ and $\widehat{\rho}_t(c_4)$ is a complex reflection, while $E_{s^{-4}}$ corresponds to a fixed point in $\CH^2$ when $\theta\in (\pi/4,\pi/3)\cup (2\pi/3,3\pi/4)$ and $\widehat{\rho}_t(c_4)$ is a scalar rotation. Similar computations lead to the same conclusions for the other $\widehat{\rho}_t(c_i)$.  
\end{proof}

\begin{rem}
    When $\theta=\pi/2$ and $t=-1$, then $\rho_{-1}(c_i)$ are involutions of the form
    \[
    \begin{pmatrix}
        1 & \ast & \ast\\
        0 & -1 & 0\\
        0 & 0 & -1
    \end{pmatrix}
    \]
    and the image of $\rho_{-1}$ is virtually abelian (index 2). When $\theta\in\{\pi/4 , 3\pi/4\}$ and $t=\pm \sqrt{2}/2$, then $\rho_{\pm\frac{\sqrt{2}}{2}}(c_i)$ are order $4$ elements of the form
    \[
    \begin{pmatrix}
        -1 & \ast & \ast\\
        0 & i & 0\\
        0 & 0 & i
    \end{pmatrix}
    \]
    and the image of $\rho_{\pm\frac{\sqrt{2}}{2}}$ is virtually abelian (index 4).
\end{rem}

\begin{lem}\label{lem:Burau-fixed-points-Zariski-dense}
    The Burau fixed points $\rho_t$ are Zariski dense representations for all values of $\theta\in (\pi/4,3\pi/4)\setminus \{\pi/2\}$.
\end{lem}
\begin{proof}
    It has been established by Formanek that the Burau representation $\beta_{4,t}'\colon B_4/Z(B_4)\to\mathrm{PGL}_3\C$ is irreducible if and only if $t^3+t^2+t+1\neq 0$, or equivalently $\theta\notin\{\pi/4,\pi/2,3\pi/4\}$~\cite[Lemma~6]{formanek}. When $\theta\in (\pi/4,3\pi/4)\setminus \{\pi/2\}$, the image of $\beta_{4,t}'$, say $\Lambda$, belongs to $\PU(2,1)$ and irreducibility implies that it fixes no point in $\CH^2\cup\partial\CH^2$ and preserves no $\C$-line.  By Lemma~\ref{lem:maximal-Zariski-closed-subgroups-of-PU(2,1)}, the Zariski closure $\overline{\Lambda}$ of $\Lambda$ is either equal to $\PU(2,1)$ or it is contained in a conjugate of $\mathrm{PO}(2,1)$. However, complex reflections and scalar rotations are contained (up to conjugation) in $\mathrm{PO}(2,1)$ if and only if they are involutions, and vertical translations are never contained in $\mathrm{PO}(2,1)$, so we may exclude the case thanks to item (4) of Lemma~\ref{lem:properties-Casella-fixed-point} and items (4)--(5) from Lemma~\ref{lem:properties-Burau-fixed-points} (recalling that $\theta\in (\pi/4,3\pi/4)\setminus \{\pi/2\})$. Hence, $\overline{\Lambda}=\PU(2,1)$.
    
    Next, we prove that when $\beta_{4,t}'$ is Zariski dense, then its restriction $\rho_t\colon\pi_1S_{0,4}\to\PU(2,1)$ remains Zariski dense. The point is that $\pi_1S_{0,4}$ is a normal subgroup of $B_4/Z(B_4)$, implying that the image of $\rho_t$, say $\Gamma$, is a normal subgroup of $\Lambda$. So, $\Lambda$ is a subgroup of the normalizer $N(\Gamma)$ of $\Gamma$. Hence, the Zariski closure of $N(\Gamma)$, which coincides with the normalizer of the Zariski closure $\overline{\Gamma}$ of $\Gamma$, is equal to $\PU(2,1)$. Since $\PU(2,1)$ is a simple group and $N(\overline{\Gamma})=\PU(2,1)$, the group $\overline{\Gamma}$ is either discrete, hence finite, or equal to $\PU(2,1)$. In the former case, $\overline{\Gamma}$ is a finite normal subgroup of $\overline{\Lambda}=\PU(2,1)$. Since $\PU(2,1)$ is connected and centre-free, $\overline{\Gamma}$ must be trivial, which is never the case when $t\in (0,\pi)$. We conclude that $\overline{\Gamma}=\PU(2,1)$ as desired.
\end{proof}

The values of $t$ for which Burau representations have discrete image have been worked out by McMullen~\cite[Corollary~7.2]{mcmullen}. In the case of four strands and in the regime $\theta\in (\pi/4,\pi/2)$, the values of $t=e^{2 i\theta}$ are precisely
\begin{equation}\label{eq:values-of-theta-for-discrete-Burau}
\theta\in \left\{\frac{3\pi}{10},\frac{\pi}{3},\frac{5\pi}{14},\frac{3\pi}{8},\frac{2\pi}{5},\frac{5\pi}{12},\frac{4\pi}{9}\right\}.
\end{equation}
\begin{cor}
    In the regime $\theta\in (\pi/4,3\pi/4)\setminus\{\pi/2\}$, the Burau fixed points $\rho_t$ have discrete image if and only if $\theta$ or $\pi-\theta$ belongs to the list~\eqref{eq:values-of-theta-for-discrete-Burau}, otherwise they have dense image.
\end{cor}
\begin{proof}
    By McMullen's result, it is enough to verify that when the image $\Gamma$ of $\rho_t$ is discrete, then so is the image $\Lambda$ of $\beta_{4,t}'$. Since $\Gamma$ is discrete, and also Zariski dense by Lemma~\ref{lem:Burau-fixed-points-Zariski-dense}, its normalizer $N(\Gamma)$ in $\PU(2,1)$ is discrete. However, $\Lambda$ normalizes $\Gamma$, hence it is discrete too. Furthermore, by Corollary~\ref{cor:Zariski-dense-implies-discrete-or-dense}, if $\Gamma$ is not discrete, then it must be dense because it is Zariski dense.
\end{proof}

\subsection{Pure braid group representations}\label{sec:pure-braid-group-representation} Pure braid groups admit more irreducible representations than braid groups. Classical examples include the 4-parameter family of Gassner representations, which we now introduce using different perspectives.

\subsubsection{Gassner fixed points}
Let
\[\Omega_0 = \left \{\left (z_i\right )_{1\leqslant i\leqslant 5}\in \U(1)^5\mid z_1\cdots z_5 = 1\right \}.\]
For $\alpha = (z_i)_{1\leqslant i\leqslant 5}$ in $\Omega_0$, let $\rho_\alpha\colon \pi_1S_{0,4}\to \mathrm{PGL}_3\C$ be the representation defined by:

\begin{align*}
\
\rho_\alpha(c_1) = \begin{pmatrix}
    z_2z_3z_4 & 0 & 0\\
    z_3z_4(1-z_2) & 1 & 0\\
    z_4(1-z_3) & 0 & 1
\end{pmatrix},\quad
\rho_\alpha(c_2) = \begin{pmatrix}
    1 & 1-z_1 & 0\\
    0 & z_1 z_3z_4 & 0\\
    0 & z_1z_4(1-z_3) & 1
\end{pmatrix}
\end{align*}
\begin{align*}
\rho_\alpha(c_3) = \begin{pmatrix}
    1 & 0 & z_2(1-z_1)\\
    0 & 1 & 1-z_2\\
    0 & 0 & z_1z_2z_4
\end{pmatrix},
\quad
\rho_\alpha(c_4) = \begin{pmatrix}
    1 + z_2z_3(z_1 - 1) & z_2z_3(z_1 - 1) & z_2z_3(z_1 - 1)\\
    z_3(z_2 - 1) & 1 + z_3(z_2 - 1) & z_3(z_2 - 1)\\
    z_3 -1 & z_3 - 1 & z_3
\end{pmatrix}.
\end{align*}
This indeed defines a representation $\pi_1S_{0,4}\to\mathrm{PGL}_3\C$, since one can check that \[\rho_\alpha(c_1)\rho_\alpha(c_2)\rho_{\alpha}(c_3)\rho_\alpha(c_4) = z_5^{-1}\cdot1.\]

\begin{defn}\label{def:gassner_fixed_point}
    A \emph{Gassner fixed} point is a representation conjugate to $\rho_\alpha\colon \pi_1S_{0,4}\to \mathrm{PGL}_3\C$ for some $\alpha\in \Omega_0$.
\end{defn}

We will recall later why conjugacy classes of Gassner fixed points are indeed fixed by $\PMod(S_{0,n})$ (Lemma~\ref{lem:gassner-fixed-points-are-restrictions-hence-fixed-points}) and for which parameters $\alpha$ they preserve a Hermitian form of signature $(2,1)$ (Proposition~\ref{prop:gassner-parameters-signature-(2,1)}).

\begin{lem}\label{lem:Gassner-peripheral-monodromies-are-pseudo-refelections}
    For every $\alpha\in\Omega_0$, each peripheral curve of $S_{0,4}$ is mapped by $\rho_\alpha$ to a pseudo-reflection.
\end{lem}
\begin{proof}
Indeed, for every $1\leqslant i\leqslant 4$, we can write $\rho_\alpha(c_i) = 1 + A_iB_i$, where $(B_1, B_2, B_3)$ is the canonical basis of $\C^3$, as row vectors, $B_4 = \begin{pmatrix}
    1 & 1& 1
\end{pmatrix}$, and
\begin{align}
    A_1 &= \begin{pmatrix}
        z_2z_3z_4 - 1 \\ z_3z_4(1-z_2) \\z_4(1-z_3)    \end{pmatrix} ,\quad A_2 = \begin{pmatrix}
            1 - z_1\\
            z_1z_3z_4 - 1\\
            z_1z_4(1-z_3)
        \end{pmatrix},\label{eq:vectors-Ai-gassner}\\
        A_3 &= \begin{pmatrix}
            z_2(1 - z_1)\\
            1-z_2\\
            z_1z_2z_4 - 1
        \end{pmatrix},\quad 
        A_4 = \begin{pmatrix}
            z_2z_3(z_1-1)\\
            z_3(z_2-1)\\
            z_3 - 1
        \end{pmatrix}.\nonumber\qedhere
\end{align}
\end{proof}
Using Formanek's criterion \cite[Theorem 5]{formanek}, we can easily decide when $\rho_\alpha$ is irreducible.
\begin{lem}\label{lem:gassner_reducible}
    Let $\alpha\in \Omega_0$. The representation $\rho_\alpha$ is reducible if and only if there exists $1\leqslant i\leqslant 5$ such that $z_i = 1$.
\end{lem}
\begin{proof}
    Let us form the matrix $M = (B_i A_j)_{1\leqslant i,j\leqslant 3}$, using the vectors from~\eqref{eq:vectors-Ai-gassner}. In other words, $M$ is the matrix whose columns are $A_1, A_2, A_3$.
    We check that the determinant of $M$ is $(z_4-1)(z_1z_2z_3z_4 - 1)$.
    Let $\mathcal G$ be the oriented graph with vertices $\{1, 2, 3\}$ and an arrow $i\to j$ if and only if $M_{ij}\neq 0$.
    If $z_i = 1$ for $i\leqslant 3$, then there is no path starting from $i$.
    If however $z_i\neq 1$ for all $i\leqslant 3$, then $\mathcal G$ is complete.
    It follows from \cite[Theorem 5]    {formanek} that $\rho_\alpha$ is irreducible if and only if $\mathcal G$ is strongly connected and $M$ is invertible.
    This happens exactly when $z_i\neq 1$ for every $1\leqslant i\leqslant 5$.
\end{proof}
Let $\Omega_1\subset\Omega_0$ be the subset of all $\alpha$ such that $\rho_\alpha$ is irreducible.

\subsubsection{Variation of local systems}
We now show that the representations $\rho_\alpha$ arise in the context of variations of local systems studied by Dettweiler--Wewers in~\cite{DW06}.
In particular we show that they are fixed points, by extending them to $\PMod_{0, 5}$.
Let $S = S_{0,5}$, $x\in S$ and write $\pi = \pi_1(S, x)=\langle c_1,\ldots, c_5:\prod c_i=1\rangle$.
An element $f\in \Mod(S)$ induces a permutation of the punctures and thus acts on $\Omega_1$ by permuting the $z_i$. 
We identify an element $\alpha = (z_i)_{1\leq i\leq 5}\in \Omega_1$ with the unique homomorphism $\pi\to \C$ sending $c_i$ to $z_i$.
For $\alpha\in \Omega_1$, consider the vector space \[V_\alpha = H^1\left (\pi, \mathbb C_{\alpha}\right )\]
defined as the three dimensional vector space $Z_\alpha / B_\alpha$, where 
\[
Z_\alpha = \{t\colon \pi\to \mathbb C\mid t(ab) = t(a) + \alpha(a) t(b)\} \quad\text{and}\quad B_\alpha = \C (1-\alpha).
\]
An element $f\in \Mod(S, x)$ induces an automorphism of $\pi$ and a linear map $\Phi(f, \alpha)\colon V_\alpha\to V_{f\cdot \alpha}$. 
The image of the point-pushing map $\pi\to \Mod(S, x)$ acts by inner automorphisms on $\pi$, and one checks that they induce scalar multiplications on $V_\alpha$.
Therefore $\Phi(f, \alpha)$ is well-defined for $f\in \Mod(S)$ up to scalar multiplication.
We follow the convention of \cite{DW06} and let linear maps act on the right.
The maps $\Phi(f, \alpha)$ satisfy $\Phi(\mathrm{id}_S, \alpha) = 1$ and
\begin{equation}\label{eq:cocycle_phi}
\Phi(f, \alpha)\cdot\Phi(g, f\cdot \alpha) = \Phi(g\circ f, \alpha).    
\end{equation}
The subgroup of $\Mod(S)$ fixing the fifth puncture is isomorphic to $B_4 / Z(B_4)$.
Let us fix the bases $\left ([e^\alpha_1], [e^\alpha_2], [e_\alpha^3]\right )$ of $V_\alpha$ defined by $e^\alpha_i(c_i) = z_i - 1$ and $e^\alpha_i(c_j) = 0$ for $1\leq i\neq j\leq 3$.
Using \cite[Section 2.2]{DW06}, we compute the image of Artin's generators of $B_4 / Z(B_4)$ in these bases:
\[
    \Phi(\sigma_1, \alpha) = \begin{pmatrix}
    0 & z_2 & 0\\
    1 & 1-z_2 & 0\\
    0 & 0 & 1
\end{pmatrix}, \,
\Phi(\sigma_2, \alpha) = \begin{pmatrix}
    1 & 0 & 0\\
    0 & 0 & z_3\\
    0 & 1 & 1-z_3
\end{pmatrix}, \,
\Phi(\sigma_3, \alpha) = \begin{pmatrix}
    1 & 0 & 0\\
    0 & 1 & 0\\
    -z_4 & -z_4 & -z_4
\end{pmatrix}.
\]
Since pure braids act trivially on $\Omega_0$ and $\PMod(S) = PB_4 / Z(B_4)$ we get representations 
\[
\Phi_\alpha := \Phi(\cdot, \alpha)\colon\PMod(S)\to \mathrm{PGL}(V_\alpha).
\]
\begin{rem}
    One can use \eqref{eq:cocycle_phi} to compute $\Phi(f, \alpha)\colon V_\alpha\to V_{f\cdot \alpha}$ for any $f$ fixing the fifth puncture from the matrices $\Phi(\sigma_i, \beta)$.
\end{rem}


\begin{lem}\label{lem:gassner-fixed-points-are-restrictions-hence-fixed-points}
    The restriction of $\Phi_\alpha$ to $\pi_1 S_{0,4}$ is $\rho_\alpha$. Consequently, the representations $\rho_\alpha$ are fixed points.
\end{lem}
\begin{proof}
    We only show that $\Phi_\alpha(c_1) = \rho_\alpha(c_1)$, and leave the other cases to the reader.
    Recall from \eqref{eq:morphism-pi_1S_0,4-to-B_4/Z(B_4)} that $c_1$ is sent by the point-pushing map to $x^{-3}$, where $x = \sigma_3\sigma_2$ in $B_4 / Z(B_4)$.
    Using \eqref{eq:cocycle_phi}, 
    \begin{align*}
     \Phi(x, \alpha) &= \Phi(\sigma_2, \alpha)\cdot\Phi(\sigma_3, \sigma_2\cdot \alpha)\\
     &= \begin{pmatrix}
    1 & 0 & 0\\
    0 & 0 & z_3\\
    0 & 1 & 1-z_3
\end{pmatrix} \begin{pmatrix}
    1 & 0 & 0\\
    0 & 1 & 0\\
    -z_4 & -z_4 & -z_4
\end{pmatrix}.\\
&= \begin{pmatrix}
    1 & 0 & 0\\
    -z_3z_4 & -z_3z_4 & -z_3z_4\\
    z_4(z_3 - 1) & 1 + z_4(z_3 - 1) & z_4(z_3 - 1)
\end{pmatrix}.
    \end{align*}
    Therefore,
    \begin{align*}
        \Phi(x^2, \alpha) &= \Phi(x, \alpha)\cdot\Phi(x, x\cdot \alpha)\\
        &= \begin{pmatrix}
    1 & 0 & 0\\
    -z_3z_4 & -z_3z_4 & -z_3z_4\\
    z_4(z_3 - 1) & 1 + z_4(z_3 - 1) & z_4(z_3 - 1)
\end{pmatrix}\begin{pmatrix}
    1 & 0 & 0\\
    -z_2z_4 & -z_2z_4 & -z_2z_4\\
    z_2(z_4 - 1) & 1 + z_2(z_4 - 1) & z_2(z_4 - 1)
\end{pmatrix}\\
&= \begin{pmatrix}
    1 & 0 & 0 \\
    z_3z_4(z_2 - 1) & z_3z_4(z_2 - 1) & z_2z_3z_4\\
    -z_4(z_2z_3 - z_3 + 1) & -z_4(z_2z_3 - z_3 + 1) & -z_2z_3z_4
\end{pmatrix}.
    \end{align*}
    \begin{align*}
        \Phi(x^3, \alpha) &= \Phi(x^2, \alpha)\cdot\Phi(x, x^2\cdot \alpha)\\
        &= \begin{pmatrix}
    1 & 0 & 0 \\
    z_3z_4(z_2 - 1) & z_3z_4(z_2 - 1) & z_2z_3z_4\\
    -z_4(z_2z_3 - z_3 + 1) & -z_4(z_2z_3 - z_3 + 1) & -z_2z_3z_4
\end{pmatrix}\begin{pmatrix}
    1 & 0 & 0\\
    -z_2z_3 & -z_2z_3 & -z_2z_3\\
    z_3(z_2 - 1) & 1 + z_3(z_2 - 1) & z_3(z_2 - 1)
\end{pmatrix}\\
&= \begin{pmatrix}
    1 & 0 & 0\\
    z_3z_4(z_2-1) & z_2z_3z_4 & 0\\
    z_4(z_3-1) & 0 & z_2z_3z_4
\end{pmatrix}
    \end{align*}
    Finally, we get $\Phi_\alpha(x^{-3}) = \Phi_\alpha(x^3)^{-1} = z_1z_5\rho_\alpha(c_1)$.
\end{proof}

\subsubsection{Extending Burau representations}\label{sec:Burau-is-Gassner}
For $\alpha \in \Omega_0$ of the form $\alpha = (t, t, t, t, z)$ with $z=t^{-4}$, the matrices $\Phi(f, \alpha)$ form a representation $B_4\to \mathrm{PGL}(V_\alpha)$.
If 
\[
P = \begin{pmatrix}
    z^2 & 0 & 0\\
    -z
    ^2 & z & 0\\
    0 & -z & 1
\end{pmatrix},
\]
then we check that
\[\Phi(\sigma_i, \alpha)P = P \beta_4 (\sigma_i),\]
where $\beta_4(\sigma_i)$ can be found in~\eqref{eq:Burau-matrices-4}.
Therefore, the Burau fixed points can be seen as a special case of Gassner fixed points where $z_1 = z_2 = z_3 = z_4$, with the Casella fixed points corresponding to the case $z_1 = z_2 = z_3 = z_4 = e^{\pm 2i\pi/3}$.

\subsubsection{Monodromy of hypergeometric functions and shapes of polyhedra}
Deligne--Mostow computed in \cite{DM86} the signature of some hermitian form preserved by $\Phi_\alpha$. In particular, it follows from \cite[Section 2]{DM86} (see also \cite[Theorem 2.3]{DW06}) that $\Phi_\alpha$ preserves a Hermitian form $h_\alpha$ whose signature can be computed from $\alpha$.
\begin{prop}\label{prop:gassner-parameters-signature-(2,1)}
    Let $(z_i)_{1\leqslant i\leqslant 5}\in \Omega_1$.
    Write $z_i = e^{2i\pi \mu_i}$ with $\mu_i\in (0, 1)$.
    There exists a non-degenerate Hermitian form $h_\alpha$ that is preserved by $\Phi_\alpha$ whose signature is \[\left(\sum_{i=1}^5\mu_i - 1, 4 - \sum_{i=1}^5 \mu_i\right).\]
\end{prop}

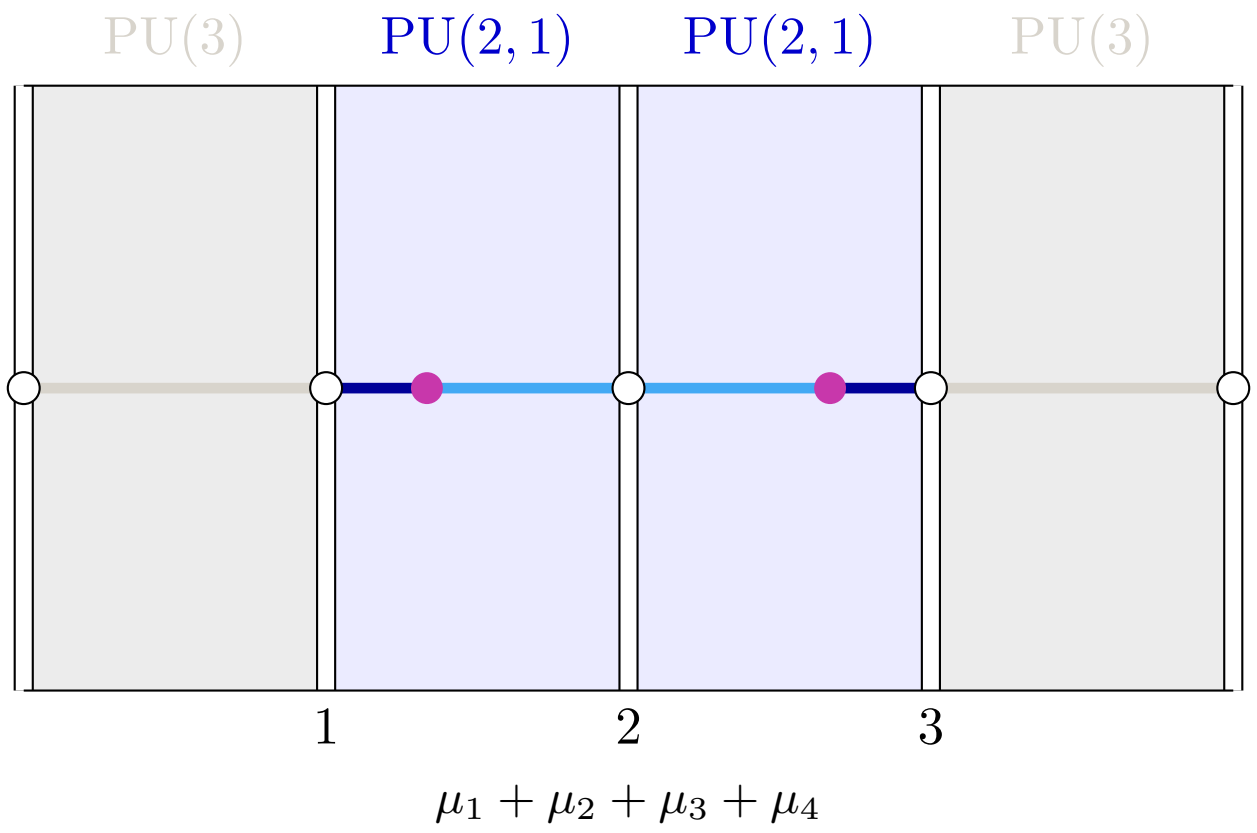
\begin{figure}[h]
    \centering
    \begin{tikzpicture}[x=2cm]
  \fill[gray!15] (0,0) rectangle (1,4);
  \fill[gray!15] (3,0) rectangle (4,4);
  \fill[blue!8]  (1,0) rectangle (3,4);

  \foreach \x in {0,1,2,3,4}
    \draw[double, double distance=3pt] (\x,0) -- (\x,4);

  \draw (0,4) -- (4,4);
  \draw (0,0) -- (4,0);

  \draw[smoked, line width=2pt] (0,2) -- (1,2);
  \draw[smoked, line width=2pt] (3,2) -- (4,2);
  \draw[blue!60!black, line width=2pt] (1,2) -- (4/3,2);
  \draw[blue!60!black, line width=2pt] (8/3,2) -- (3,2);
  \draw[sky, line width=2pt] (4/3,2) -- (8/3,2);

  \fill[mauve] (4/3,2) circle (3pt);
  \fill[mauve] (8/3,2) circle (3pt);

  \foreach \x in {0,1,2,3,4}
    \draw[fill=white] (\x,2) circle (3pt);

  \node[smoked, above] at (0.5,4) {$\PU(3)$};
  \node[blue!80!black, above] at (1.5,4) {$\PU(2,1)$};
  \node[blue!80!black, above] at (2.5,4) {$\PU(2,1)$};
  \node[smoked, above] at (3.5,4) {$\PU(3)$};

  \node[below] at (1,0) {$1$};
  \node[below] at (2,0) {$2$};
  \node[below] at (3,0) {$3$};
  \node[below] at (2,-.5) {\small $\mu_1+\mu_2+\mu_3+\mu_4$};
\end{tikzpicture}
    \caption{The shaded blue region in the middle is a representation of the (4-dimensional) domain $\Omega$ where $\mu_1+\mu_2+\mu_3+\mu_4\in (1,3)\setminus\{2\}$ corresponds to the values of $\alpha$ for which the Gassner fixed point $\rho_\alpha$ preserves a Hermitian form of signature $(2,1)$. The horizontal line in the middle represents the locus of Burau fixed points (with the same color code as on Figure~\ref{fig:signature-Burau-form}), with the two Casella fixed points in mauve.}
    \label{fig:Gassner-range}
\end{figure}

Denote by $\Omega\subset \Omega_1$ the subset of all $\alpha$ for which $\rho_\alpha$ is irreducible and $h_\alpha$ is of signature $(2,1)$. We illustrate $\Omega$ on Figure~\ref{fig:Gassner-range}.

In \cite{DM86}, Deligne--Mostow studied the monodromy of hypergeometric functions and use them to construct a family of complex hyperbolic structures on $\mathcal M_{0,n}$.
The associated holonomy representations $\PMod(S_{0,n})\to \mathrm{PU}(n-3, 1)$ coincide with $\Phi_\alpha$, when $n= 5$ and $\alpha\in \Omega$.
Thurston gave another geometric interpretation of these complex hyperbolic structures, as local period coordinates on the space of $\mathcal M_{0,n}$, identified with the space of polyhedra with $n$ cone angles in $(0, 2\pi)$ \cite{Thurston2022}.

\subsubsection{Gassner representations}
The reduced Gassner representations are representations \[PB_n\to \mathrm{GL}_{n-1}(\mathbb Z[t_1^{\pm 1}, \ldots, t_n^{\pm 1}])\] obtained from applying Fox calculus to the action of braid groups on free groups.
They are for example presented in \cite[Section 3.3]{birman}.
In \cite[Theorem C]{KM01}, Kapovich--Millson showed that the representations $\Phi_\alpha$ are specializations of Gassner representations at roots of unity.

Abdulrahim described in \cite{abdulrahim} the roots of unity for which the specialisation of Gassner representations is reducible. 
It implies that $\Phi_\alpha$ is reducible if and only if $z_5 = 1$.
 From Lemma~\ref{lem:gassner_reducible}, there exists $\alpha\in \Omega_0$ such that $\Phi_\alpha$ is irreducible while $\rho_\alpha$ is not.
 We recover the well-known fact that Burau representations restricted to pure braids are particular cases of Gassner representations, see e.g.~\cite[Section 3.3]{birman}

\begin{lem}\label{lem:gassner-fixed-points-interior-scc-are-RE-or-SP}
    Let $\alpha\in\Omega$ and $\rho=\rho_\alpha\colon \pi_1 S_{0, 4} \to \PU(2,1)$ be a Gassner fixed point. Let $\gamma\in\pi_1S_{0,4}$ represent a simple closed curve. If $\gamma$ is peripheral, then $\rho(\gamma)$ is a pseudo-reflection. Otherwise, $\rho(\gamma)$ is either regular elliptic or screw parabolic.
\end{lem}
\begin{proof}
We have seen in Lemma~\ref{lem:Gassner-peripheral-monodromies-are-pseudo-refelections} that the peripheral curves are pseudo-reflections.
The curve $\gamma = c_1c_2$ satisfies 
\[\rho(\gamma) =\begin{pmatrix}
z_2z_3z_4 & z_2z_3z_4(1-z_1) & 0\\
z_3z_4(1-z_2) & z_3z_4(1-z_2+z_1z_2) & 0\\
z_4(1-z_3) & z_4(1-z_3) & 1
\end{pmatrix}. \]
Its eigenvalues are $1, z_3z_4$, and $z_5^{-1}$.
If $z_3z_4z_5\neq 1$ then $\rho_\alpha(\gamma)$ is regular elliptic.
Otherwise $z_1z_2 = 1$ and $\rho(\gamma) - z_3z_4\cdot1$ has its first $2\times 2$ minor \[\begin{pmatrix}
    z_3z_4(z_2 -1 ) & z_2z_3z_4(1-z_1)\\
    z_3z_4(1-z_2) & z_3z_4(1-z_2)
\end{pmatrix} = z_3z_4(z_2-1)\begin{pmatrix}
    1 & 1\\
    -1 & -1
\end{pmatrix}.
\]
Therefore, $\mathrm{rank}(\rho(c_1c_2) - z_3z_4\cdot 1) = 2$ and the only possibility is that $\rho(c_1c_2)$ is screw parabolic.
Similarly, we check that $\rho(c_2c_3)$ and $\rho(c_2c_4)$ are either regular elliptic or screw parabolic.
Since $\rho$ is a fixed point and $\PMod(S_{0,4})$ maps any non-peripheral simple closed curve to either $c_1c_2$, $c_2c_3$, or $c_2c_4$, we are done.
\end{proof}

\begin{lem}\label{lem:gassner-fixed-points-reducible-on-all-pairs-of-pants}
    Let $\rho\colon \pi_1 S_{0,4}\to \PU(2,1)$ be a Gassner representation.
    The restriction of $\rho$ to any pair of pants is reducible.
\end{lem}
\begin{proof}
    Let $P_{12}, P_{23}, P_{24}$ be the pair of pants on $S_{0,4}$ whose boundary curves are homotopic respectively to $(c_1, c_2, c_1c_2), (c_2, c_3, c_2c_3), (c_2, c_4, c_2c_4)$.
    For any pair of pants $P\subset S_{0,4}$, there exists $f\in \PMod(S_{0,4})$ that sends $P$ to one of $P_{12}, P_{23}, P_{24}$.
    Since $\rho$ is a fixed point, it suffices to show that $\rho$ is reducible when restricted to these pair of pants. This is the case since the restrictions of $\rho$ to these pairs of pants preserve the point in $\CP^2$ represented by, respectively,

    \[
    V_{12} = \begin{pmatrix}
        0 \\ 0 \\1    \end{pmatrix}, \quad 
        V_{23} = \begin{pmatrix}
            1\\  0\\   0
        \end{pmatrix},\quad
        V_{24} = \begin{pmatrix}
            1\\
            0\\
            -1
        \end{pmatrix}.\qedhere
\]
\end{proof}

\section{Reducible fixed points in genus zero}\label{sec:reducible-fixed-points-genus-0}
\subsection{Overview}
The goal of this section is to classify all representations $\rho\colon\pi_1S_{0,n}\to\PU(2,1)$ with reducible image (in the sense of Definition~\ref{defn:irreducible}) whose conjugacy classes are mapping class group fixed points (in short: the reducible fixed points). We will prove the following.
\begin{thm}\label{thm:reducible-fixed-points-are-alomst-abelian}
    If $\rho\colon\pi_1S_{0,n}\to\PU(2,1)$ is a reducible fixed point, then $\rho$ is pants-like in the sense of Definition~\ref{defn:pants-like}.
\end{thm}
The proof of Theorem~\ref{thm:reducible-fixed-points-are-alomst-abelian} unfolds over the next few sections and is concluded in Section~\ref{sec:proof-thm-reducible-implies-almost-abelian}. We will typically treat the case $n=4$ first and then generalize to all $n\geq 5$. Recall that, by Lemma~\ref{lem:maximal-reducible-subgroups-of-PU(2,1)}, $\rho$ is reducible if and only if its image is contained in $\stab(X)$, $\stab(L)$, or $\stab(Y)$ for some point $X\in\CH^2$, $\C$-line $L\subset\CH^2$, and boundary point $Y\in\partial\CH^2$. Also recall that we denote by $Z_\rho$ the centralizer of the image of a representation $\rho\colon\pi_1S_{0,4}\to\PU(2,1)$ and we say that $\rho$ is a fixed point to mean that $[\rho]$ is a fixed point of the mapping class group action.

\subsection{A useful lemma}
We begin with a general result about mapping class group fixed points.
\begin{lem}\label{lem:general-lemma-fixed-points}
Let $G$ be a centre-free group, $S=S_{g,n}$ be an oriented surface of genus $g\geq 0$ with $n\geq 0$ punctures, and $\rho\colon \pi_1 S\to G$ be a fixed point. Then, for every $1\neq \gamma\in\pi_1S$ representing a simple closed curve that separates $S$ into two sub-surfaces $S_1\sqcup S_2$, with the restrictions of $\rho$ to $\pi_1S_i$ denoted by $\rho_i=\rho\vert_{\pi_1S_i}$, one of the following is true:
\begin{enumerate}
    \item
    $Z_{\rho_2}=\{1\}$ and $\rho(\gamma)\in Z_{\rho_1}$;
    \item
    $Z_{\rho_1}=\{1\}$ and $\rho(\gamma)\in Z_{\rho_2}$;
    \item
    $Z_{\rho_1}\neq \{1\}$ and $Z_{\rho_2}\neq \{1\}$.
\end{enumerate}
Moreover, when $S=S_{0,4}$, the condition $\rho(\gamma)\in Z_{\rho_i}$ is equivalent to $\rho_i$ being abelian.
\end{lem}
\begin{proof}
    Let $\tau$ denote the Dehn twist along $\gamma$. Recall that $\tau$ is represented by the automorphism $T$ of the amalgamated product $\pi_1S = \pi_1S_1 \ast_{\langle\gamma\rangle} \pi_1 S_2$ that acts trivially on $\pi_1S_1$ and acts by conjugation by $\gamma$ on $\pi_1S_2$, see e.g.~\cite[Section~3.2]{goldman-ergodic}. In particular, $\rho_1\circ T=\rho_1$ and $\rho_2\circ T= \rho(\gamma)\rho_2\rho(\gamma)^{-1}$.
    
    Since $\rho$ is a fixed point, the representations $\rho\circ T$ and $\rho$ are conjugate. Therefore, there exists $A\in Z_{\rho_1}$ such that
    $A\rho(\gamma)\in Z_{\rho_2}$. If $Z_{\rho_1}=\{1\}$, then $A=1$ and $\rho(\gamma)\in Z_{\rho_2}$. Else, if $Z_{\rho_2}=\{1\}$, then $\rho(\gamma)=A^{-1}\in Z_{\rho_1}$. Otherwise, $Z_{\rho_1}\neq\{1\}$ and $Z_{\rho_2}\neq\{1\}$.

    It remains to consider the case $S=S_{0,4}$. Each simple closed curve $\gamma$ separates $S_{0,4}$ into two pairs of pants $S_1$ and $S_2$ with fundamental groups $\pi_1S_i=\langle \alpha_i,\beta_i,\gamma: \alpha_i\beta_i\gamma=1\rangle$. Therefore, the condition $\rho(\gamma)\in Z_{\rho_i}$ means that $[\rho_i(\gamma),\rho_i(\alpha_i)]=[\rho_i(\gamma),\rho_i(\beta_i)]=1$. Since $\rho_i(\beta_i)=\rho_i(\alpha_i)^{-1}\rho_i(\gamma)^{-1}$, we also have $[\rho_i(\alpha_i),\rho_i(\beta_i)]=1$, implying that $\rho_i$ is abelian. 
\end{proof}

\subsection{Stabilizing a point or a \texorpdfstring{$\C$}-line}
We start by considering the cases where $\Im(\rho)\subset \stab(X)$ or $\Im(\rho)\subset \stab(L)$.
\begin{prop}\label{prop:reducible-fixed-point-in-stab(X)-or-stab(L)-are-almost-abelian}
    If $\rho\colon \pi_1 S_{0,n}\to \PU(2,1)$ is a reducible fixed point and $\Im(\rho)\subset \Lambda$ for $\Lambda=\stab(X)$ or $\Lambda=\stab(L)$, then $\rho$ is pants-like.
\end{prop}
\begin{proof}
    Let $\rho\vert^\Lambda\colon \pi_1S_{0,n}\to\Lambda$ be the co-restriction of $\rho$ to $\Lambda$. Since $\rho$ is a fixed point, for every $f\in\Aut^\ast (\pi_1S_{0,n})$, there is $g\in \PU(2,1)$ such that $\rho\circ f=g\rho g^{-1}$. In particular, if $\Lambda = \stab(X)$, then $\Im(\rho)$ also fixes $g^{-1}X$. So, either $X$ is the only fixed point of $\Im(\rho)$ and $g\in \Lambda$, or $\Im(\rho)$ fixes another point $X'$ and $\Im(\rho)\subset \mathcal{CR}_{XX'}$, hence $\rho$ is abelian. Similarly, if $\Lambda=\stab(L)$, then either $g\in \Lambda$ or $\Im(\rho)$ stabilizes two $\C$-lines and $\rho$ is abelian (Lemma~\ref{lem:intersection-sstabilizers-two-C-line-is-abelian}). This shows that the co-restriction $\rho\vert^\Lambda$ is also a fixed point.

    Recall that $\stab(X)\cong \mathrm{P}(\U(2)\times\U(1))\cong \U(2)$ and $\stab(L)\cong \mathrm{P}(\U(1)\times\U(1,1))\cong \U(1,1)$. Up to multiplying $\rho\vert^\Lambda$ by a character $\chi\colon\pi_1S_{0,n}\to \U(1)$, we may assume that it takes images in $\SU(2)$ or $\SU(1,1)$, or, to put it simply, in $\mathrm{SL}_2\C$. Such fixed points have been classified: they are either abelian or map all but at most three peripheral curves to $\pm 1$ (see e.g.~\cite[Lemma~3.3]{bronstein-maret}). We conclude that $\rho$ is pants-like.
\end{proof}

\subsection{Stabilizing a boundary point}
We will now consider the case where $\Im(\rho)\subset \stab(Y)$ for a boundary point $Y\in\partial\CH^2$. Let $H$ denote the horosphere foliation based at $Y$. We start by considering $4$-punctured spheres.

\begin{lem}\label{lem:reducible-fixed-point-stab(Y)-reduce-to-stab(H)}
    If $\rho\colon\pi_1 S_{0,4}\to\PU(2,1)$ is a reducible fixed point with $\Im(\rho)\subset \stab(Y)$, then either $Z_\rho\neq \{1\}$ and $\Im(\rho)\subset \stab(H)$, or $\rho$ is pants-like.
\end{lem}
\begin{proof}
    Assume $\rho$ is not pants-like. In particular, it does not map any peripheral curve to $1$. In order to conclude that $\Im(\rho)\subset \stab(H)$, it is enough to prove that $\rho$ does not map any peripheral curve to a hyperbolic element. We fix a standard presentation of $\pi_1S_{0,4}$ with generators $c_1,c_2,c_3,c_4$ and write $C_i=\rho(c_i)$.
    
    For the sake of contradiction, assume that $\rho$ maps a peripheral curve to a hyperbolic element and without loss of generality assume that $C_1$ is hyperbolic. In particular, $C_1$ preserves a $\C$-line $L$ with $Y\in\partial L$. Using Table~\ref{tab:intersections-centralizers}, we observe that, for $i\neq 1$, three possibilities may occur:
    \begin{enumerate}
        \item $[C_1,C_i]=1$;
        \item $Z(\langle C_1, C_i\rangle)=\mathcal{CR}_L$;
        \item $Z(\langle C_1, C_i\rangle)=\{1\}$.   
    \end{enumerate}
    We now consider the nine possible cases resulting from choosing $i=2$ and $i=4$ above. The symmetry in the indices allows us to only consider six cases.
    \begin{enumerate}
        \item[(1-1)]
        $[C_1,C_2]=[C_1,C_4]=1$: This forces $C_1\in Z_\rho$, contradiction.
        \item[(1-2)]
        $[C_1,C_2]=1$ and $Z(\langle C_1,C_4\rangle)=\mathcal{CR}_L$:  In this case, $\Im(\rho)\subset \stab(L)$ and so $\rho$ is pants-like by Proposition~\ref{prop:reducible-fixed-point-in-stab(X)-or-stab(L)-are-almost-abelian}, contradiction.
        \item[(1-3)]
        $[C_1,C_2]=1$ and $Z(\langle C_1, C_4\rangle)=\{1\}$: The centralizer condition implies $C_4C_1\neq 1$ and $C_4\notin \stab(L)$. The commutator condition implies $C_2\in \stab(L)$. Also, by Lemma~\ref{lem:general-lemma-fixed-points}, we have $C_4C_1\in Z(\langle C_2, C_3\rangle)$, thus $C_2\in \stab(L)\cap Z(C_4C_1)$. Since $C_4C_1\in \stab(Y)\setminus \stab(L)$, we can use Table~\ref{tab:intersections-centralizers} to see that $\stab(L)\cap Z(C_4C_1)=\mathcal{VT}_Y$. So, $C_2$ is a vertical translation at $Y$, contradicting $[C_1,C_2]=1$.
        \item[(2-2)]
        $Z(\langle C_1,C_2\rangle)=Z(\langle C_1,C_4\rangle)=\mathcal{CR}_L$: This forces $\Im(\rho)\subset \stab(L)$ and so $\rho$ is pants-like by Proposition~\ref{prop:reducible-fixed-point-in-stab(X)-or-stab(L)-are-almost-abelian}, contradiction.
        \item[(2-3)]
        $Z(\langle C_1, C_2\rangle)=\mathcal{CR}_L$ and $Z(\langle C_1, C_4\rangle)=\{1\}$: In this case, $C_2\in\stab(L)$ and the same reasoning as in (1-3) shows that $C_2$ must be a vertical translation at $Y$ because $Z(\langle C_1, C_4\rangle)=\{1\}$. Hence, $C_1C_2$ is a hyperbolic element stabilizing $L$. Since $\rho$ is a fixed point, there exists $A\in Z(\langle C_1,C_2\rangle)=\mathcal{CR}_L$ such that $AC_1C_2\in Z(\langle C_3, C_4\rangle)$ (see the proof of Lemma~\ref{lem:general-lemma-fixed-points}). Note that $AC_1C_2$ is hyperbolic and stabilizes $L$. It commutes with $C_4$, so $C_4$ stabilizes $L$, contradicting $Z(\langle C_1, C_4\rangle)=\{1\}$.
        \item[(3-3)]
        $Z(\langle C_1,C_2\rangle)=Z(\langle C_1, C_4\rangle)=\{1\}$: From the statement in Lemma~\ref{lem:general-lemma-fixed-points} about 4-punctured spheres, we conclude that $[C_3,C_4]=1$ and $[C_2,C_3]=1$, so $C_3\in Z_\rho$, contradiction. \qedhere
    \end{enumerate}
\end{proof}

\begin{lem}\label{lem:reducible-fixed-point-stab(H)-are-totally-hyperbolic}
    If $\rho\colon \pi_1S_{0,4}\to\PU(2,1)$ is a reducible fixed point with $\Im(\rho)\subset \stab(H)$, then either $\rho$ is pants-like or all peripheral monodromies are regular parabolic and $\Im(\rho)\subset\mathcal{N}_Y$.
\end{lem}
\begin{proof}
    Assume that $\rho$ is not pants-like. Since $Z(\stab(H))=\mathcal{VT}_Y$, we conclude that the images of all peripheral curves are non-trivial and cannot be vertical translations. So, peripheral monodromies are either regular parabolic, screw parabolic, or complex reflections.

    With the same notation as in Lemma~\ref{lem:reducible-fixed-point-stab(Y)-reduce-to-stab(H)}, assume for the sake of contradiction that $C_1$ is screw parabolic stabilizing a $\C$-line $L$ or a complex reflection with mirror $L$, with $Y\in \partial L$. Note that, according to Table~\ref{tab:intersections-centralizers}, every other $C_i$ satisfies one of the following:
    \begin{enumerate}
        \item $[C_1, C_i]=1$;
        \item $Z(\langle C_1, C_i\rangle) =\mathcal{VT}_Y$.
    \end{enumerate}
    We proceed with the same case distinction as in the proof of Lemma~\ref{lem:reducible-fixed-point-stab(Y)-reduce-to-stab(H)} with $i=2$ and $i=4$.
    \begin{enumerate}
        \item[(1-1)] $[C_1,C_2]=[C_1,C_4]=1$: This forces $C_1\in Z_\rho$, contradiction.
        \item[(1-2)] $[C_1,C_2]=1$ and $Z(\langle C_1,C_4\rangle)=\mathcal{VT}_Y$: Here, $C_2$ is either screw parabolic and stabilizes $L$, or is a complex reflection with mirror $L$. Moreover, since $\rho$ is a fixed point, there exists $A\in Z(\langle C_1,C_4\rangle)=\mathcal{VT}_Y$ such that $AC_4C_1\in Z(\langle C_2,C_3\rangle)$. Since $AC_4C_1$ commutes with $C_2$, it must stabilize $L$, hence $C_4$ stabilizes $L$, contradicting $Z(\langle C_1,C_4\rangle)=\mathcal{VT}_Y$.
        \item[(2-2)] $Z(\langle C_1,C_2\rangle)=Z(\langle C_1,C_4\rangle)=\mathcal{VT}_Y$: As before, there exists $A\in Z(\langle C_1,C_4\rangle)=\mathcal{VT}_Y$ such that $AC_4C_1\in Z(\langle C_2,C_3\rangle)$. Since $Z(\langle C_1,C_4\rangle)=\mathcal{VT}_Y$, $C_4$ does not stabilize $L$, thus $AC_4C_1$ neither. In particular, $AC_4C_1$ fixes $Y$ but is not a vertical translation. Using Table~\ref{tab:intersections-centralizers}, this forces $[C_2,C_3]=1$. Similarly, $[C_3,C_4]=1$, and thus $C_3\in Z_\rho$, contradiction.
    \end{enumerate}
    This means that the only possibility is for all peripheral monodromies of $\rho$ to be regular parabolic. The product of regular parabolic elements fixing the same boundary point $Y$ is an element of $\mathcal{N}_Y$, hence $\Im(\rho)\subset \mathcal{N}_Y$.
    \end{proof}

\begin{prop}\label{prop:reducible-fixed-point-in-stab(y)-are-almost-abelian}
    If $\rho\colon \pi_1S_{0,4}\to\PU(2,1)$ is a reducible fixed point with $\Im(\rho)\subset \stab(Y)$, then $\rho$ is pants-like.
\end{prop}
\begin{proof}
    Assume for the sake of contradiction that $\rho$ is not pants-like. Then, by Lemmas~\ref{lem:reducible-fixed-point-stab(Y)-reduce-to-stab(H)} and~\ref{lem:reducible-fixed-point-stab(H)-are-totally-hyperbolic}, it follows that all peripheral monodromies of $\rho$ are regular parabolic and $\Im(\rho)\subset \mathcal{N}_Y$.

    We denote the peripheral monodromies of $\rho$ by $C_i=A(\zeta_i,v_i)$, using the notation~\eqref{eq:Heisenberg-translations}.
    Recall the product rule: $A(\zeta,v)A(\zeta',v')=A(\zeta+\zeta',v+v'-2\Im(\overline{\zeta}\zeta'))$. Also remember that $A(\zeta,v)$ and $A(\zeta',v')$ commute if and only if $\overline{\zeta}\zeta'\in\R$, and the product $A(\zeta,v)A(\zeta',v')$
    is a vertical translation or the identity
    if and only if $\zeta+\zeta'=0$.

    Since $\rho$ is a fixed point, there are $A\in Z(\langle C_1,C_2\rangle)$ such that $AC_1C_2\in Z(\langle C_3, C_4\rangle)$ and $B\in Z(\langle C_4, C_1\rangle)$ such that $BC_4C_1\in Z(\langle C_2, C_3\rangle)$.
    This forces $A,B\in \mathcal{N}_Y$ and if we write $A=A(\zeta_a,v_a)$ and $B=A(\zeta_b,v_b)$, then
    \[
    \overline{\zeta_a}\zeta_1, \overline{\zeta_a}\zeta_2, \overline{\zeta_b}\zeta_1,
    \overline{\zeta_b}\zeta_4\in\R.
    \]
    We consider the following cases.
    \begin{enumerate}
        \item $\zeta_a, \zeta_b\neq 0$ : This implies $\overline{\zeta_1}\zeta_2\in\R$ and $\overline{\zeta_1}\zeta_4\in\R$. So. $C_1$ commutes with $C_2$ and $C_4$, and thus $C_1\in Z_\rho$, contradiction.
        \item $\zeta_a=0$ and $\zeta_b\neq 0$: This means that $\overline{\zeta_1}\zeta_4\in\R$, so $[C_1,C_4]=1$, and $A$ is trivial or a vertical translation. Since $\rho$ is not pants-like, $C_1$ does not commute with $C_2$, and so $C_1C_2$ is non-trivial and not a vertical translation. Therefore, $AC_1C_2$ is non-trivial and not a vertical translation. However, $AC_1C_2\in Z(\langle C_3, C_4\rangle)$, so $[C_3,C_4]=1$, and $C_4\in Z_\rho$, contradiction.
        \item $\zeta_a\neq 0$ and $\zeta_b= 0$: This case is analogous to (2).
        \item $\zeta_a=\zeta_b=0$: This means that $A$ and $B$ are either trivial or vertical translations. If $C_1C_2$ is either trivial or a vertical translation, then $[C_1,C_2]=1$. If not, then as in (2) we conclude that $[C_3,C_4]=1$. Similarly, arguing about $C_4C_1$, we conclude that either $[C_4,C_1]=1$ or $[C_2,C_3]=1$. In all configurations, there is a peripheral monodromy $C_i$ that commutes with the other three, contradiction. \qedhere
    \end{enumerate}
\end{proof}

\subsection{Proof of Theorem~\ref{thm:reducible-fixed-points-are-alomst-abelian}}\label{sec:proof-thm-reducible-implies-almost-abelian}
\begin{proof}
    When $\Im(\rho)\subset \stab(X)$ or $\Im(\rho)\subset \stab(L)$, then the conclusion follows from Proposition~\ref{prop:reducible-fixed-point-in-stab(X)-or-stab(L)-are-almost-abelian}. If $\Im(\rho)\subset \stab(Y)$ and $n=4$, then the conclusion follows from Proposition~\ref{prop:reducible-fixed-point-in-stab(y)-are-almost-abelian}. It remains to prove that, for all $n\geq 5$, if $\rho\colon\pi_1S_{0,n}\to \PU(2,1)$ is a reducible fixed point with $\Im(\rho)\subset \stab(Y)$, then $\rho$ is pants-like.

    We will prove this statement by induction on $n$, using the case $n=4$ as base case. Let $n\geq 5$ and assume that all fixed points $\rho'\colon \pi_1S_{0,n'}\to\PU(2,1)$ with $\Im(\rho')\subset \stab(Y)$ are pants-like for $n'<n$. Let now $\rho\colon\pi_1S_{0,n}\to\PU(2,1)$ be a fixed point with $\Im(\rho)\subset\stab(Y)$ and peripheral monodromies $C_1,\ldots,C_n$. We will assume that $\Im(\rho)\not\subset \stab(L)$ for any $\C$-line $L$ (else the other case applies) and $C_i\neq 1$ for all $i$ (else the induction hypothesis applies immediately). For the sake of contradiction, we will also assume that $\rho$ is not pants-like.

    First consider the case where a peripheral monodromy $C_i$ of $\rho$ is hyperbolic and stabilizes a unique $\C$-line $L$. Since $\Im(\rho)\not\subset\stab(L)$, there is $j\neq i$ such that $C_j$ does not stabilize $L$. This implies that $Z_\rho=\{1\}$. Up to relabelling the punctures, we may assume that $i=1$ and $j=2$. For every $k\in\{1,\ldots,n-1\}$, consider the sub-sphere $S^{(k)}$ of $S_{0,n}$ obtained from $S_{0,n}$ by removing a disk that contains the punctures $k$ and $k+1$ (and no other punctures). The restriction of $\rho$ to $\pi_1S^{(k)}$ is denoted by $\rho^{(k)}$. It is still a fixed point with image contained in $\stab(Y)$ and has peripheral monodromy $C_1, C_2, \ldots, C_{k-1}, C_kC_{k+1}, C_{k+2},\ldots, C_n$. By construction, if $k\geq 3$, then $C_1$ and $C_2$ are still peripheral monodromies of $\rho^{(k)}$ and hence $Z_{\rho^{(k)}}=\{1\}$. The induction hypothesis implies that $\rho^{(k)}$ is pants-like, meaning that at least $n-4$ peripheral monodromies among $\{C_j : j\neq 1,2,k,k+1\}\cup \{C_kC_{k+1}\}$ belong to $Z_{\rho^{(k)}}$, hence are trivial. If $n\geq 6$, this implies that $C_\ell=1$ for some $\ell$, contradiction. If $n=5$ and $C_kC_{k+1}=1$ for $k\in\{3,4\}$, then $C_3=C_4^{-1}=C_5$ and $C_4=C_1C_2$. Hence, $C_4$ does not stabilize $L$ and $Z_{\rho^{(2)}}=\{1\}$. Applying the induction hypothesis to $\rho^{(2)}$ implies that $C_2C_3=1$ (or $C_\ell=1$ for some $\ell$) and thus $C_1=1$, contradiction.

    Now, consider the case where none of the peripheral monodromies of $\rho$ are hyperbolic.  In that case, $\Im(\rho)\subset \stab(H)$ and $\mathcal{VT}_Y\subset Z_\rho$. The only other elements that $Z_\rho$ may contain are complex reflections, screw parabolic and regular parabolic elements. If it contains a complex reflection with mirror $L$ or a screw parabolic element stabilizing $L$, then $\Im(\rho)\subset \stab(L)$, and if it contains a regular parabolic element, then $\rho$ is abelian, contradiction. Therefore, we have $Z_\rho = \mathcal{VT}_Y$. If some monodromy $C_k$ of $\rho$ is a vertical translation, then consider the restriction $\rho^{(k)}$ of $\rho$ to $S^{(k)}$. Note that $Z_{\rho^{(k)}}=\mathcal{VT}_Y$. So, by the induction hypothesis, at least $n-4$ peripheral monodromies among $\{C_j: j\neq k, k+1\}\cup \{C_kC_{k+1}\}$ are vertical translations. Since $C_{k+1}\neq 1$, it is a vertical translation if and only if $C_kC_{k+1}$ is. Therefore, there are always at least $n-3$ peripheral monodromies of $\rho$ that are vertical translations, implying that $\rho$ is pants-like, contradiction.
    
    We may thus assume that no peripheral monodromy of $\rho$ is a vertical translation. Since $\rho$ is not abelian, there are two peripheral monodromies $C_i$ and $C_j$ that do not commute. In particular, $Z(\langle C_i, C_j\rangle)=\mathcal{VT}_Y$. As before, we may assume that $i=1$ and $j=2$ and we consider the restrictions $\rho^{(k)}$ of $\rho$ to the sub-spheres $S^{(k)}$ for $k\geq 3$. Note that $Z_{\rho^{(k)}}=\mathcal{VT}_Y$. The induction hypothesis applied to $\rho^{(k)}$ says that at least $n-4$ peripheral monodromies among $\{C_j : j\neq k, k+1\}\cup \{C_kC_{k+1}\}$ are vertical translations. If $n\geq 6$, this implies that $C_\ell$ is a vertical translation for some $\ell$, contradiction. If $n=5$ and $C_kC_{k+1}$ are vertical translations for $k\in\{3,4\}$, then $C_4=C_1C_2V$ for some vertical translation $V$. Therefore, $[C_4,C_1]=[C_1C_2,C_1]=C_1[C_2,C_1]C_1^{-1}\neq 1$ and $Z_{\rho^{(2)}}=\mathcal{VT}_Y$. Applying the induction hypothesis to $\rho^{(2)}$ implies that $C_2C_3$ is a vertical translation (or $C_\ell$ is for some $\ell$) and thus $C_1=(C_4C_5)^{-1}(C_2C_3)^{-1}$ is a vertical translation, contradiction.
\end{proof}
\section{Reducing to pseudo-reflections for 4-punctured spheres}\label{sec:reducing-to-pseudo-reflections}
\subsection{Overview}
The goal of this section is to reduced the classification of irreducible fixed points to representations with peripheral monodromies consisting only of pseudo-reflections (Definition~\ref{def:pseudo-reflections}), i.e.~complex reflections, scalar rotations, or vertical translations. The following terminology will be useful.
\begin{defn}\label{def:corked}
    A representation $\rho\colon\pi_1S_{g,n}\to\PU(2,1)$ is \emph{corked} if all its peripheral monodromies are non-trivial.
\end{defn}
The main result of this section is the following.
\begin{thm}\label{thm:excluding-regular-elliptic-regular-parabolic-screw-parabolic-regular-hyperbolic}
    Let $\rho\colon\pi_1S_{0,4}\to\PU(2,1)$ be a corked fixed point. Assume that $\rho$ is regular (i.e.~$Z_\rho=\{1\}$). Then all peripheral curves are mapped to pseudo-reflections.
\end{thm}
We divide the proof into several steps, increasing in technicity:
first we prove that no peripheral curve can be sent to a hyperbolic, regular parabolic, or screw parabolic element.
Then, we prove that if there is a regular elliptic peripheral monodromy, then one pair of pants in $S_{0,4}$ must have abelian image. Finally, we prove that, in this case, the representation must have non-trivial centre. 

We recommend keeping Tables~\ref{tab:conjugacy-classes} and~\ref{tab:intersections-centralizers} at hand.

\subsection{Hyperbolic, regular parabolic or screw parabolic peripheral monodromies}

It will be useful to fix a standard presentation of $\pi_1S_{0,4}$:
\[
\pi_1S_{0,4}=\langle c_1,c_2,c_3,c_4 : \prod  c_i=1\rangle.
\]

Throughout this section, we assume that $\rho\colon\pi_1S_{0,4}\to\PU(2,1)$ is a fixed point and that the kernel of $\rho$ contains no peripheral curves. We will write $C_i=\rho(c_i)$. A useful fact is the following.

\begin{fact}\label{fact:at-most-one-pair-of-pants-with-trivial-centralizer}
    At most one of the following subgroups of $\PU(2,1)$ may have a trivial centralizer:
    \[
    \langle C_1,C_2\rangle, \quad \langle C_1,C_3\rangle, \quad \langle C_1,C_4\rangle.
    \]
\end{fact}
\begin{proof}
    For the sake of contradiction, assume that at least two of them have trivial centralizer. For instance, if $Z(\langle C_1,C_2\rangle)=Z(\langle C_1,C_4\rangle)=\{1\}$, then, by Lemma~\ref{lem:general-lemma-fixed-points}, the complementary groups $\langle C_3,C_4\rangle$ and $\langle C_2,C_3\rangle$ are abelian. So, $C_3$ commutes with $C_2$ and $C_4$, and thus also with $C_1=C_4^{-1}C_3^{-1}C_2^{-1}$. Hence, $C_3\in Z_{\rho}\subset Z(\langle C_1,C_2\rangle)=\{1\}$, contradicting the assumption $C_3\neq 1$. Similarly, if $Z(\langle C_1,C_2\rangle)=Z(\langle C_1,C_3\rangle)=\{1\}$, then, by Lemma~\ref{lem:general-lemma-fixed-points}, the complementary groups $\langle C_3,C_4\rangle$ and $\langle C_1C_2C_1^{-1},C_4\rangle$ are abelian. In this case, $C_4$ commutes with $C_3$ and $C_1C_2C_1^{-1}=C_4^{-1}C_3^{-1}C_1^{-1}$. So, $C_4$ commutes with $C_1$ and as before $C_4\in Z_\rho\subset Z(\langle C_1,C_2\rangle)=\{1\}$, contradiction. The last case is similar.
\end{proof}

\begin{lem}\label{lem:excluding-regular-hyperbolic}
    Assume that $\rho$ maps a peripheral curve to a hyperbolic element. Then $\Im(\rho)\subset \stab (L)$ for a $\C$-line $L\subset\CH^2$.
\end{lem}
\begin{proof}
    Assume that $C_1$ is hyperbolic, hence stabilizing a unique $\C$-line $L$. By Fact~\ref{fact:at-most-one-pair-of-pants-with-trivial-centralizer}, at least two among $\langle C_1,C_2\rangle$, $\langle C_1,C_3\rangle$, and $\langle C_1,C_4\rangle$ have non-trivial centralizer. Using Table~\ref{tab:intersections-centralizers}, we note that $\langle C_1,C_i\rangle$ has non-trivial centralizer if and only if $C_i$ stabilizes $L$. So, at least two among $C_2,C_3,C_4$ must stabilize $L$ and since $C_1C_2C_3C_4=1$, all four do and $\Im(\rho)\subset \stab(L)$.
\end{proof}
\begin{lem}\label{lem:excluding-regular-parabolic}
    Assume that $\rho$ maps a peripheral curve to a regular parabolic element. Then $\Im(\rho)\subset \stab(H)$ for a horosphere foliation $H$.
\end{lem}
\begin{proof}
    Let $C_1$ be the regular parabolic image of a peripheral curve, fixing a point $Y\in\partial\CH^2$ and stabilizing the horosphere foliation $H$ based at $Y$. By Fact~\ref{fact:at-most-one-pair-of-pants-with-trivial-centralizer}, at least two among $\langle C_1,C_2\rangle$, $\langle C_1,C_3\rangle$, and $\langle C_1,C_4\rangle$ have non-trivial centralizer. Using Table~\ref{tab:intersections-centralizers}, observe that $\langle C_1,C_i\rangle $ has non-trivial centralizer if and only if $C_i$ fixes $Y$ and $C_i$ is not hyperbolic, i.e.~$C_i$ stabilizes $H$. So, at least two among $C_2,C_3,C_4$ must stabilize $H$, and hence all of them do. So, $\Im(\rho)\subset\stab (H)$.
\end{proof}
\begin{lem}\label{lem:excluding-screw-parabolic}
Assume that $\rho$ maps a peripheral curve to a screw parabolic element. Then either $\Im(\rho)\subset \stab(L)$ or $\Im(\rho)\subset \stab(H)$, where $L$ is a $\C$-line and $H$ a horosphere foliation.
\end{lem}
\begin{proof}
Let $C_1$ be the screw parabolic image of a peripheral curve, stabilizing a unique $\C$-line $L$ and a horosphere foliation $H$ based at $Y\in\partial\CH^2$. By Fact~\ref{fact:at-most-one-pair-of-pants-with-trivial-centralizer}, at least two among $\langle C_1,C_2\rangle$, $\langle C_1,C_3\rangle$, and $\langle C_1,C_4\rangle$ have non-trivial centralizer. From Table~\ref{tab:intersections-centralizers}, note that $\langle C_1,C_i\rangle$ has non-trivial centralizer if and only if
$C_i$ stabilizes $L$ or $C_i$ is not hyperbolic and fixes $Y$, i.e.~$C_i$ stabilizes $H$. So, at least two among $C_2,C_3,C_4$ must stabilize either $L$ or $H$. If two stabilize $L$, or if two stabilize $H$, then we are done. It remains to consider the mixed case where, say, $C_2$ stabilizes $L$, $C_3$ stabilizes $H$, and $Z(\langle C_1, C_4\rangle )=\{1\}$. Note that Lemma~\ref{lem:general-lemma-fixed-points} implies that $\langle C_2,C_3\rangle$ is abelian. Assuming that $C_3$ does not stabilize $L$, then $C_3$ is one of the following isometries:
\begin{itemize}
    \item a regular parabolic element and $C_2\in Z(C_3)\subset \stab(H)$;
    \item a screw parabolic element and $C_2\in Z(C_3)\subset \stab(H)$;
    \item a complex reflection through any mirror $L'\neq L$ containing $Y$ and $C_2\in Z(C_3)\cap \stab(L)=\stab(L')\cap \stab(L)=\mathcal{VT}_Y\subset \stab(H)$ (see Lemma~\ref{lem:intersection-sstabilizers-two-C-line-is-abelian}).
\end{itemize}
So, in all cases $C_2$ stabilizes $H$ and we are done.
\end{proof}

\subsection{Regular elliptic peripheral monodromy}
We will prove the following statement.

\begin{lem}\label{lem:excluding-regular-elliptic}
    Assume that $\rho$ maps a peripheral curve to a regular elliptic element. Then, either $\Im(\rho)\subset \stab(X)$ or $\Im(\rho)\subset\stab(L)$, where $X\in \CH^2$ and $L\subset \CH^2$ is a $\C$-line.
\end{lem}
\begin{proof}
We will split the proof into two different claims. Again, we assume that $\rho\colon\pi_1S_{0,4}\to\PU(2,1)$ is a fixed point and that the kernel of $\rho$ contains no peripheral curves; and we write $C_i=\rho(c_i)$. We also assume that $C_1$ is the image of a peripheral curve that is regular elliptic.

If one among $C_2,C_3,C_4$ is hyperbolic or screw parabolic, then the conclusion follows from Lemmas~\ref{lem:excluding-regular-hyperbolic} and~\ref{lem:excluding-screw-parabolic}. Furthermore, none of $C_2,C_3,C_4$ can be regular parabolic by Lemma~\ref{lem:excluding-regular-parabolic}. Let us now assume  that no peripheral curve is of those types. Denote by $X\in\CH^2$ the fixed point of $C_1$ and by $L_1\perp L_2$ the stabilized $\C$-lines. Using Table~\ref{tab:intersections-centralizers}, note that, if $A\neq 1$, then $\langle C_1,A\rangle$ has non-trivial centralizer if and only if $A$ fixes $X$ or preserves one of $L_1$ or $L_2$.

\begin{claim}\label{claim:excluding-regular-elliptic-no-abelian-pants}
    If no pair of pants in $S_{0,4}$ has abelian image under $\rho$, then $\Im(\rho)\subset \stab(X)$ or $\Im(\rho)\subset\stab(L)$.
\end{claim}
\begin{proofclaim}
Since no pair of pants has abelian image, by Lemma~\ref{lem:general-lemma-fixed-points}, all subgroups $\langle C_1,C_i\rangle$ for $i\neq 1$ have non-trivial centralizer. If two among $C_2,C_3,C_4$ both preserve $L_1$, both preserve $L_2$, or both fix $X$, so does the whole image of the representation and we are done.

The remaining cases are the following.
\begin{itemize}
    \item $C_2\in\stab(X)$, $C_3\in\stab (L_1)$, and $C_4\in\stab(L_2)$: Here, $C_4(X)=C_3^{-1}C_2^{-1}C_1^{-1}(X)=C_3^{-1}(X)$. Since $C_4(X)\in L_2$, $C_3^{-1}(X)\in L_1$, and $L_1\cap L_2=\{X\}$, we must have $C_4(X)=C_3^{-1}(X)=X$, implying $\Im(\rho)\subset \stab(X)$.
    \item $C_2\in \stab(L_1)$, $C_3\in\stab(X)$, and $C_4\in\stab(L_2)$: Here, $C_2(X)=C_1^{-1}C_4^{-1}C_3^{-1}(X)=C_1^{-1}C_4^{-1}(X)$. Since both $C_4$ and $C_1$ stabilize $L_2$, we have, as before, $C_2(X)=C_4^{-1}(X)=X$ and $\Im(\rho)\subset \stab(X)$.
    \item All other cases can be treated similarly.\qedhere
\end{itemize}

\end{proofclaim}

\begin{claim}\label{claim:excluding-regular-elliptic-abelian-pants}
    If some pair of pants in $S_{0,4}$ has abelian image under $\rho$, then $\Im(\rho)\subset\stab(X)$ or $\Im(\rho)\subset \stab(L)$.
\end{claim}
\begin{proofclaim}
Up to the good choice of generators of $\pi_1S_{0,4}$, the abelian image of a pair of pants is of the following form: $\langle C_1,C_i\rangle $ for $i\neq 1$ or $\langle C_i,C_j\rangle$ with $i\neq j> 1$, with $C_1$ being regular elliptic.

First, consider the case where $\langle C_1,C_i\rangle$ is abelian for some $i\neq 1$. In that case, $C_i$ also fixes $X$ and stabilizes $L_1$ and $L_2$. By Fact~\ref{fact:at-most-one-pair-of-pants-with-trivial-centralizer}, there is at least one $j\in\{2,3,4\}\setminus \{i\}$ such that $\langle C_1,C_j\rangle$ has non-trivial
centralizer. Using Table~\ref{tab:intersections-centralizers}, $C_j$ must either fix $X$ or stabilize $L_1$ or $L_2$. So, $C_1,C_i,C_j$ either all fix $X$ or all stabilize $L_1$ or all stabilize $L_2$, proving that $\Im(\rho)\subset \stab(X)$ or $\Im(\rho)\subset \stab(L_k)$ for some $k\in\{1,2\}$.

Now, let us assume that $\langle C_i,C_j\rangle$ is abelian for some $i\neq j> 1$ and that all the groups $\langle C_1, C_k\rangle$ for $k\neq 1$ are not abelian. Without loss of generality, we may assume that $\{i,j\}=\{2,3\}$. Note that $\langle C_1,C_2\rangle$ and $\langle C_1,C_3\rangle$ must have non-trivial centralizer, for otherwise either $\langle C_3,C_4\rangle$ or $\langle C_3^{-1}C_2C_3,C_4\rangle = \langle C_2,C_4\rangle$ is abelian. In the first case, we would conclude that $C_3\in Z_\rho$ and, in the second case $C_2\in Z_\rho$, contradicting the assumption that $\langle C_1,C_2\rangle$ and $\langle C_1,C_3\rangle$ are not abelian. 

Using that $\langle C_1,C_2\rangle$ and $\langle C_1,C_3\rangle$ are non-abelian with non-trivial centralizer, Table~\ref{tab:intersections-centralizers} implies that $C_2$ and $C_3$ must be of one of types in Table~\ref{tab:excluding-regular-elliptic-possible-peripheral-monodromies-abelian-pair-of-pants-case}.

\begin{table}[h]
    \centering
\begin{tabular}{c|c|c}
\toprule
Conjugacy type & Fixed data & Invariant objects
\\\midrule
\multirow{3}{*}{RE}  & $X_i=X$ & $L_1'\perp L_2'$ and $\{L_1',L_2'\}\neq \{L_1,L_2\}$\\\cline{2-3}
& $X_i\neq X$, $X_i\in L_1$ &  $L_1'=L_1$, $L_2'\perp L_1$\\\cline{2-3}
& $X_i\neq X$, $X_i\in L_2$ &  $L_1'=L_2$, $L_2'\perp L_2$\\
\midrule
\multirow{3}{*}{CR} & $M_i\ni X$, $M_i\notin \{L_1,L_2\}$ &
all $\C$-lines $\perp M_i$\\\cline{2-3}
& $M_i\perp L_1$, $X\notin M_i$ & all $\C$-lines $\perp M_i$, including $L_1$\\\cline{2-3}
& $M_i\perp L_2$, $X\notin M_i$ & all $\C$-lines $\perp M_i$, including $L_2$\\
\midrule
\multirow{2}{*}{SR} &
$X_i\in L_1$, $X_i\neq X$ & all $\C$-lines through $X_i$, including $L_1$\\\cline{2-3}
&
$X_i\in L_2$, $X_i\neq X$ & all $\C$-lines through $X_i$, including $L_2$\\
\midrule
\multirow{2}{*}{VT} & $Y_i\in\partial L_1$ & all $\C$-lines containing $Y_i$, including $L_1$\\\cline{2-3}
& $Y_i\in\partial L_2$ & all $\C$-lines containing $Y_i$, including $L_2$\\ 
\bottomrule
\end{tabular}
    \caption{The possible types for the peripheral monodromies $C_i$ with $i\in\{2,3\}$.}
    \label{tab:excluding-regular-elliptic-possible-peripheral-monodromies-abelian-pair-of-pants-case}
\end{table}

Note that if $C_2$ or $C_3$ is regular elliptic, then, since $\langle C_2,C_3\rangle$ is abelian, all three of $C_1,C_2,C_3$ will either fix $X$ or stabilize a $C$-line, and so will $\Im(\rho)$.
Now, assuming that $C_2$ and $C_3$ are not regular elliptic, and using that $\langle C_2, C_3\rangle$ is abelian, we are left with the combinations in Table~\ref{tab:excluding-regular-elliptic-possible-peripheral-monodromies-combinations}. 

\begin{table}[ht]
\centering
\begin{tabular}{c||c|c}
\toprule
Conjugacy type & Configuration & $\rho$-invariant data
\\\midrule
\multirow{3}{*}{CR--CR} & $M_2\ni X,\quad M_3\ni X, \quad M_2=M_3$ or $M_2\perp M_3$ & $\stab(X)$\\\cline{2-3}
& $M_2=M_3 \perp L_1, \quad X\notin M_2,M_3$ & $\stab(L_1)$\\\cline{2-3}
& $M_2=M_3 \perp L_2, \quad X\notin M_2,M_3$ & $\stab(L_2)$\\
\midrule

\multirow{2}{*}{CR--SR} & $M_2\perp L_1, \quad X_3\in M_2\cap L_1$ & $\stab(L_1)$\\\cline{2-3}
& $M_2\perp L_2, \quad X_3\in M_2\cap L_2$ & $\stab(L_2)$\\
\midrule

CR--VT & cannot commute ($Y_3\notin \partial M_2$) & \\
\midrule

\multirow{2}{*}{SR--SR} & $X_2=X_3\in L_1\setminus \{X\}$ & $\stab(L_1)$\\\cline{2-3}
& $X_2=X_3\in L_2\setminus \{X\}$ & $\stab(L_2)$\\
\midrule

SR--VT & cannot commute & \\
\midrule

\multirow{2}{*}{VT--VT} & $Y_2=Y_3\in\partial L_1$ & $\stab(L_1)$\\\cline{2-3}
& $Y_2=Y_3\in\partial L_2$ & $\stab(L_2)$\\
\bottomrule
\end{tabular}
\caption{Possible resulting combinations for $(C_2,C_3)$ using that $\langle C_2,C_3\rangle$ is abelian.}
\label{tab:excluding-regular-elliptic-possible-peripheral-monodromies-combinations}
\end{table}

In all cases, we deduce that $\Im(\rho)$ fixes a point or stabilizes a $\C$-line, as desired.
\end{proofclaim}

Claims~\ref{claim:excluding-regular-elliptic-no-abelian-pants} and~\ref{claim:excluding-regular-elliptic-abelian-pants} prove Lemma~\ref{lem:excluding-regular-elliptic}.
\end{proof}

\subsection{Proof of Theorem~\ref{thm:excluding-regular-elliptic-regular-parabolic-screw-parabolic-regular-hyperbolic}}
\begin{proof}
    Let $\rho\colon\pi_1S_{0,4}\to\PU(2,1)$ be a fixed point. We assume that the kernel of $\rho$ contains no peripheral curves and that $Z_\rho=\{1\}$. According to Lemmas~\ref{lem:excluding-regular-hyperbolic}--\ref{lem:excluding-regular-elliptic}, if $\rho$ maps a peripheral curve to a hyperbolic, regular parabolic, screw parabolic, or regular elliptic element, then $\Im(\rho)$ is contained in $\stab(X)$, $\stab(L)$, or $\stab(H)$ for some point $X\in \CH^2$, $\C$-line $L\subset \CH^2$, and horosphere foliation $H\subset \CH^2$. In all those cases, we would have $Z_\rho\neq \{1\}$ (see Lemma~\ref{lem:maximal-subgroups-of-PU(2,1)-with-non-trivial-centralizer}). Therefore, all peripheral monodromies of $\rho$ are complex reflections, scalar rotations, or vertical translations.
\end{proof}

\section{Irreducible fixed points in genus zero}\label{sec:irreducible-fixed-points-genus-0}
\subsection{Overview}
The goal of this section is to classify all representations $\pi_1 S_{0,n}\to \PU(2,1)$ with irreducible image (Definition~\ref{defn:irreducible}) whose conjugacy class is fixed by the mapping class group action (in short: irreducible fixed points). We will prove the following.
\begin{thm}\label{thm:irreducible-fixed-points-genus-zero}
Let $\rho\colon\pi_1 S_{0,n}\to \PU(2,1)$ be an irreducible fixed point with $n\geq 3$. 
\begin{itemize}
    \item If $n=3$, then $\rho$ is any irreducible representation.
    \item If $n=4$ and $\rho$ is corked, then $\rho$ is a Gassner fixed point (Definition~\ref{def:gassner_fixed_point}).
    \item If $n\geq 5$, then either $\rho$ has exactly $n-4$ trivial peripheral monodromies and, after filling the trivial punctures, is a Gassner fixed point, or $\rho$ has exactly $n-3$ trivial peripheral monodromies and, after filling the punctures, is an irreducible pair of pants representation.
\end{itemize}
\end{thm}
The case $n=3$ is immediate since $\PMod(S_{0,3})$ is trivial. We start with the case $S=S_{0,4}$, then proceed to the case of spheres with more punctures. 

\subsection{Four-punctured spheres}
Recall from Definition~\ref{def:pseudo-reflections} that a matrix $C\in \mathrm{GL}_3\C$ is a \emph{pseudo-reflection} if $C - 1$ has rank one. The next result classifies irreducible representations $\pi_1 S_{0,4}\to \mathrm{PU}(2,1)$ with peripheral monodromies being pseudo-reflections (without fixed point assumption). Abdulrahim observed in~\cite[discussion after Theorem~2]{abdulrahim} that Formanek's result~\cite[Theorem~5 (b)]{formanek} can be interpreted as a classification of free group representations $F_r\to \mathrm{GL}_r\C$ that map all $r$ generators to pseudo-reflections. The next result shows that, when $r=3$ and with the extra hypothesis that the product of the images of the generators is also a pseudo-reflection, then the representation must be Gassner.

\begin{prop}\label{prop:irreducible-fixed-pts-4-punctures-spheres}
    Let $\rho\colon \pi_1 S_{0,4}\to \mathrm{PU}(2,1)$ be an irreducible representation whose peripheral monodromies are pseudo-reflections. Then, $\rho$ is a Gassner fixed point (Definition~\ref{def:gassner_fixed_point}).
\end{prop}
\begin{proof}
  Let $C_1,C_2,C_3,C_4\in \U(2,1)$ be lifts of the peripheral monodromies of $\rho$ that are pseudo-reflections. In particular, there exists $\mu\in \C^*$ such that 
\[C_1C_2C_3C_4 = \mu 1.\]  
    We will think of the quadruple $(C_1,C_2,C_3,C_4)$ as a representation of the rank-4 free group $F_4$ into $\U(2,1)$.
    Let us first observe that we can consider representations $F_4\to \U(2,1)$ up to conjugation in $\mathrm{GL}_3 \C$ instead of $\U(2,1)$.
    \begin{claim}\label{claim:conjugate-in-U(2,1)-iff-conjugate-in-GL_3C}
    Two irreducible representations $\rho_1, \rho_2\colon F_4\to \U(2,1)$ are conjugate in $\U(2,1)$ if and only if they are conjugate in $\mathrm{GL}_3\C$. (The conclusion remains true if $F_4$ is replaced by any finitely generated group.)
    \end{claim}
    \begin{proofclaim}
    Suppose that $\rho_1 = P \rho_2 P^{-1}$, with $P\in \mathrm{GL}_3\C$, and suppose after rescaling if necessary that $\det P = 1$.
    The representation $\rho_1$ preserves the hermitian form $J$, as well as $K = (P^{-1})^\ast J P$, since $\rho_2$ preserves $J$.
    For any $A$ in the image of $\rho_1$, we have $A^* J A = J$ thus $ J^{-1}{\left (A^*\right )}^{-1} =  AJ^{-1}$.
    We also have $A^* K A = K$, thus $\left (A^*\right )^{-1} K = KA$.
    Therefore
    \[A\left (J^{-1}K\right ) = J^{-1} \left ( A^*\right )^{-1} K = \left (J^{-1} K\right ) A.\]
    Schur's lemma gives that $J^{-1} K = \lambda I$, for some $\lambda\in \C^*$, that is $P^* J P = \lambda J$. Since $P^* J P$ and $J$ are hermitian, $\lambda\in \R$.
    But $|\det P|^2 \det J = \lambda^3 \det J$, and $\lambda = 1$. Therefore $P\in \U(2,1)$.
    \end{proofclaim}
For each $i\in\{1,2,3,4\}$, the matrix $C_i$ is a pseudo-reflection. Therefore, there exists  $(A_i, B_i)\in M_{3\times 1}(\C)\times M_{1\times 3}(\C)$ such that $C_i = 1 - A_iB_i$.
Consider the matrices $M = (m_{ij})_{1\leqslant i,j\leqslant 4}$ where $m_{ij} = B_iA_j$, and $M'\in M_{3\times 3}(\C)$ obtained from $M$ by removing its last row and last column. Observe that the choices of $A_i$ and $B_i$ are not unique: one can replace $(A_i, B_i)$ with $(a_iA_i, a_i^{-1}B_i)$ for $a_1, a_2,a_3, a_4\in \C^*$. It changes $M$ by conjugating it with the diagonal matrix $\mathrm{diag}(a_1, a_2, a_3, a_4)$.
\begin{claim}\label{claim:M_invertible}
    The vectors $A_1, A_2, A_3$, as well as $B_1, B_2, B_3$, are linearly independent. Equivalently, the matrix $M'$ is invertible.
\end{claim}
\begin{proofclaim}
    The span $W$ of $A_1, A_2, A_3$ is preserved by $C_1, C_2$ and $C_3$. 
    Since $C_1C_2C_3C_4 = \mu 1$, it is also preserved by $C_4$ and $\dim W = 3$.
    If the $B_i$ are not independent, there exists $v\neq 0$ in $\cap_{i\leqslant 3} \ker B_i$ (where here $B_i$ is seen as a linear form). Then for every $i\in\{1,2,3\}$, we have $C_i v = \left (1 - A_iB_i
    \right )v = v$. Since $C_1C_2C_3C_4 = \mu 1$, $C_4$ must also preserve the line spanned by $v$, contradicting irreducibility.
    The last statement follows from:
    \[M' = 
    \begin{pmatrix}
        B_1\\
        B_2\\
        B_3
    \end{pmatrix}\begin{pmatrix}
        A_1 & A_2 & A_3
    \end{pmatrix}.\qedhere\]
\end{proofclaim}

Observe that conjugating $C_i$ by $P$ replaces $B_i$ with $B_iP$ for all $i$:
\[P^{-1} C_i P = P^{-1}\left (1 - A_i B_i\right )P = 1 - \left (P^{-1}A_i\right ) (B_i P).\]
We can thus assume by Claim~\ref{claim:M_invertible} that $\left (B_1, B_2, B_3\right )$ is the canonical basis of $\mathbb C^3$ (this statement was also proved by Formanek~\cite[Theorem~5 (c)]{formanek}).
Therefore, setting $\lambda_i = 1 - m_{ii}$, we get 
\[
    C_1 = \begin{pmatrix}
    \lambda_1 & 0 & 0\\
    -m_{21} & 1 & 0\\
    -m_{31} & 0 & 1
\end{pmatrix}, \quad C_2 = \begin{pmatrix}
    1 & -m_{12} & 0\\
    0 & \lambda_2 & 0\\
    0 & -m_{32} & 1
\end{pmatrix}, \quad
C_3 = \begin{pmatrix}
    1 & 0 & -m_{13}\\
    0 & 1 & -m_{23}\\
    0 & 0 & \lambda_3
\end{pmatrix}.
\]
We compute:
\[
C_1C_2C_3 = \begin{pmatrix}
    \lambda_1 & -\lambda_1m_{12} & \lambda_1(m_{12}m_{23}-m_{13}) \\
    -m_{21} & \lambda_2 + m_{12}m_{21} & -\lambda_2m_{23}-m_{12}m_{21}m_{23}+m_{21}m_{13} \\
    -m_{31} & -m_{32}+m_{12}m_{31} & \lambda_3 + m_{13}m_{31} + m_{23}m_{32} - m_{12}m_{31}m_{23}
\end{pmatrix}
\]
The matrix $C_4 = \mu(C_1C_2C_3)^{-1}$ is a complex reflection, so $C_4 - 1$ has rank one. Hence, the matrix $Q$ defined by $Q = C_1C_2C_3 - \mu 1$ has rank one.
We may simplify $Q$ by taking linear combinations of its columns. Indeed we have $Q' = QR$, where
\[
Q' = \begin{pmatrix}
    \lambda_1-\mu & -m_{12}\mu & -m_{13}\mu\\
    -m_{21} & \lambda_2-\mu & -m_{23}\mu\\
    -m_{31} & -m_{32} & \lambda_3-\mu
\end{pmatrix}, \quad
    R=\begin{pmatrix}
    1 & m_{12} & m_{13}\\
    0 & 1 & m_{23}\\
    0 & 0 & 1
\end{pmatrix}.
\]
Since $Q'$ has rank one, all its $2\times 2$ minors vanish.
It follows that

\begin{equation}
\begin{cases}
    \mu m_{12} m_{21} &= (\lambda_1 - \mu)(\lambda_2 - \mu),\\
    \mu m_{12}m_{23} &= -m_{13}(\lambda_2 - \mu),\\
    \mu m_{32} m_{23} &= (\lambda_2 - \mu)(\lambda_3 - \mu),\\
    \mu m_{31} m_{13} &= (\lambda_1 - \mu)(\lambda_3 - \mu).
\end{cases}
\label{eq:minors}
\end{equation}

\begin{claim}\label{claim:lamda_neq_mu}
    For every $1\leqslant i\leqslant 4$, we have $\lambda_i \neq \mu$.
\end{claim}
\begin{proofclaim}
    First, note that the matrix $C_4$ admits $\lambda_4$ as an eigenvalue. Indeed, with $v = A_4$, we have 
    \[C_4 v = (1 - A_4B_4)v = v - vm_{44} = \lambda_4v.\]
    Now, if $\lambda_4 = \mu$, then $C_1C_2C_3$ admits $1$ as eigenvalue because $C_4 = \mu(C_1C_2C_3)^{-1}$.
    However, 
    \[\left (C_1C_2C_3 - 1\right )R = -M'.\]
    This is a contradiction, since by Claim~\ref{claim:M_invertible}, $\det M'\neq 0$. We can cyclically permute the $C_i$ in $C_1C_2C_3C_4 = \mu 1$ to show that, similarly, $\lambda_i\neq \mu$ for every $1\leqslant i\leqslant 4$.
\end{proofclaim}

From \cref{eq:minors} and Claim~\ref{claim:lamda_neq_mu}, we get that $m_{ij}\neq 0$ for $1\leqslant i\neq j\leqslant 3$.
We can thus assume, up to conjugating $M'$ by a diagonal matrix, that 
\begin{align*}
    m_{12} &= \lambda_1^{-1}(\mu - \lambda_1)\\
    m_{23} &= \lambda_2^{-1}(\mu - \lambda_2)
\end{align*}
\cref{eq:minors} then gives
\[
\begin{cases}
    m_{21} = \lambda_1(1 - \lambda_2\mu^{-1}),\\
    m_{13} = \lambda_2^{-1}\mu(\lambda_1^{-1}\mu - 1),\\
    m_{32} = \lambda_2(1 - \lambda_3 \mu^{-1}),\\
    m_{31} = \lambda_1\lambda_2\mu^{-2}(\mu - \lambda_3).
\end{cases}
\]

Set $z_5 = \mu^{-1}$ and $z_i= \lambda_i^{-1}\mu$ for $1\leqslant i\leqslant 3$. We also set $z_4 = \left (z_1 z_2 z_3 z_5\right )^{-1}$.
We thus get
\begin{align*}
C_1= \begin{pmatrix}
    z_2z_3z_4 & 0 & 0\\
    z_3z_4(1-z_2) & 1 & 0\\
    z_4(1-z_3) & 0 & 1\\
\end{pmatrix}, \quad&
C_2= \begin{pmatrix}
    1 & 1-z_1 & 0\\
    0 & z_1z_3z_4 & 0\\
    0 & z_3z_4(1 - z_3) & 1\\
\end{pmatrix}, \quad & C_3= \begin{pmatrix}
    1 & 0 & z_2(1-z_1)\\
    0 & 1 & 1-z_2\\
    0 & 0 & z_1z_2z_4\\
\end{pmatrix},
\end{align*}
Since $C_1C_2C_3C_4 = \mu 1$, and $C_i\in U(2,1)$ for all $i$, we must have $|\mu| = 1$. For every $1\leqslant i\leqslant 3$, we have $C_i\in U(2,1)$ thus $|z_iz_5| = 1$ and $|z_i| = 1$, which also implies $|z_4|=1$. By picking $\alpha = (\alpha_i)_{1\leqslant i\leqslant 5}\in \Omega$ such that $z_i = \exp(2i\pi \alpha_i)$, then $\rho$ is conjugated to $\rho_\alpha$.
\end{proof}

\begin{cor}\label{cor:irreducible_S04}
     Let $\rho\colon \pi_1S_{0,4}\to \mathrm{PU}(2,1)$ be a corked irreducible fixed point. 
     Then, $\rho$ is a Gassner fixed point.
\end{cor}
\begin{proof}
    Since $\rho$ is irreducible, it is regular and by Theorem~\ref{thm:excluding-regular-elliptic-regular-parabolic-screw-parabolic-regular-hyperbolic} all its peripheral monodromies are pseudo-reflections. The conclusion follows from Proposition~\ref{prop:irreducible-fixed-pts-4-punctures-spheres}.
\end{proof}

\subsection{Spheres with more punctures}
We now turn to the classification of irreducible fixed point $\pi_1 S_{0,n}\to \mathrm{PU}(2,1)$ for $n\geq 5$.
One way of constructing such fixed points is to start with a fixed point $\pi_1 S_{0,4}\to \mathrm{PU(2,1)}$ and pre-compose it with the moprhism $\pi_1 S_{0,n}\to \pi_1 S_{0,4}$ induced by filling $n-4$ punctures.
We show that every irreducible fixed point $\pi_1 S_{0,n} \to \mathrm{PU}(2,1)$ is obtained by this construction.

\begin{prop}\label{prop:irreducible_many_points}
    Let $\rho\colon \pi_1 S_{0,n} \to \mathrm{PU}(2,1)$ be an irreducible fixed point. Then $\rho$ maps the peripheral curves around at least $n-4$ punctures to $1$.
\end{prop}

We first show that if one can find two peripheral loops whose image by $\rho$ acts irreducibly on $\mathbb C^3$, then one can in fact find $n-3$ peripheral loops that are killed by $\rho$. It will be convenient to fix a presentation $\pi_1S_{0,n}=\langle c_1,\ldots,c_n : \prod_i c_i =1 \rangle$ and denote the peripheral monodromies of $\rho$ by $C_i=\rho(c_i)$.

\begin{lem}\label{lem:two_irreducible}
    Suppose that $\rho\colon \pi_1 S_{0,n}\to \PU(2,1)$ is a fixed point and that there exist $1\leqslant i < j\leqslant n$ such that $ \langle C_i,C_j\rangle$ is an irreducible subgroup of $\PU(2,1)$. Then, there exists $J\subset \{1, \ldots, n\}$ of cardinality $n-3$ such that $C_k = 1$ for $k\in J$.
\end{lem}
\begin{proof}
    We proceed by induction on $n\geq 4$. When $n=4$, by Corollary~\ref{cor:irreducible_S04}, either a peripheral monodromy of $\rho$ is trivial, or $\rho$ is conjugate to a Gassner fixed point. However, any two peripheral monodromies of a Gassner fixed point generate a reducible subgroup of $\PU(2,1)$ by Lemma~\ref{lem:gassner-fixed-points-reducible-on-all-pairs-of-pants}, so the latter case cannot happen.

    Let now $n\geq 5$ and $\rho$ be as in the statement. If $C_k=1$ for some $k$, then the induction hypothesis applies and we are done. So, we may assume by contradiction that $C_k\neq 1$ for all $k=1,\ldots,n$. Let $C_i$ and $C_j$ be such that $\langle C_i,C_j\rangle$ is an irreducible subgroup of $\PU(2,1)$. Up to reordering the punctures, we may assume that $i=1$ and $j=2$. For every $k\geq 3$, consider the $4$-punctured sub-sphere $S^{(k)}$ of $S_{0,n}$ whose peripheral curves are represented by $c_1$, $c_2$, $c_k$, and $(c_1c_2c_k)^{-1}$. The restriction $\rho^{(k)}$ of $\rho$ to $S^{(k)}$ is still a fixed point and $\langle C_1,C_2\rangle$ is still irreducible. So, applying the induction hypothesis to $\rho^{(k)}$, we conclude that $C_1C_2C_k=1$ (since we are assuming $C_k\neq 1$ for every $k$). Therefore, $C_k=(C_1C_2)^{-1}$ for every $k\geq 3$. Now, consider $k_2>k_1\geq 3$, and let $S$ be the 4-punctured sub-sphere of $S_{0,n}$ whose peripheral curves are represented by $c_1$, $c_{k_1}$, $c_{k_2}$, and $(c_1c_{k_1}c_{k_2})^{-1}$. Since $\langle C_1,C_{k_1}\rangle = \langle C_1, (C_1C_2)^{-1}\rangle = \langle C_1,C_2\rangle$, the induction hypothesis applies as before and gives $C_1C_{k_1}C_{k_2}=1$, which implies $C_1=C_2^{-2}$, contradicting the irreducibility of $\langle C_1,C_2\rangle$.
\end{proof}

In order to prove Proposition~\ref{prop:irreducible_many_points}, we will need to find subsurfaces $S$ in $S_{0,n}$, that are $4$-punctured spheres, such that the restriction of $\rho$ to $\pi_1 S$ is still irreducible.
The following lemma allows us to do so.
\begin{lem}\label{lem:irreducible_four_sphere}
    Given $n\geq 3$, suppose that $A_1, \ldots, A_n\in \U(2,1)$ do not preserve any non-zero proper subspace $V\subset\mathbb C^{2,1}$. Then, there exists $1\leqslant i< j< k\leqslant n$ such that $A_{i}, A_{j}, A_{k}$ do not either.
\end{lem}
\begin{proof}
    Let $I\subset \{1, \ldots, n\}$ be of the smallest cardinality such that $\langle A_i\rangle_{i\in I}$ is irreducible. This means that, for every $i\in I$, the subgroup $\langle A_j\rangle_{j\in I\setminus \{i\}}$ is reducible, i.e.~there exists a subspace $L_i\subset \mathbb C^{2,1}$ such that $\dim L_i\in \{1,2\}$ and $A_j L_i = L_i$ for $j\in I\setminus \{i\}$. Replacing $L_i$ by its orthogonal if necessary, we can assume that $\dim L_i = 1$ and write $L_i = \mathbb Cp_i$ for each $i\in I$. For the sake of contradiction, assume that $|I|\geq 4$.
    \begin{claim}
        If $J\subset I$ and $|J| = 3$, then $(p_j)_{j\in J}$ is a basis of $\mathbb C^{2,1}$.
    \end{claim}
    \begin{proofclaim}
        If $p_{j_0} = {\scriptstyle\sum_{j\neq j_0}} \lambda_j p_j$, then $W = \mathrm{span}\{p_j\mid j\neq j_{0}\}$ would be invariant by $\langle A_i\rangle_{i\in I}$ and by irreducibility $W = \mathbb C^{2,1}$, which is impossible.
    \end{proofclaim}
    Relabelling if necessary, assume that $\{1,2,3,4\}\subset I$. 
    Since $A_1 L_i = L_i$ for $i=2,3,4$, the matrix $A_1$ is diagonal in the basis $(p_2, p_3, p_4)$.
    \begin{claim}
        If a basis $(e_1, e_2, e_3)$ of $\C^{2,1}$ consists of eigenvectors of a non-scalar matrix $A\in \U(2,1)$, then one of the $e_i$ is orthogonal to the other two.
    \end{claim}
    \begin{proofclaim}
        This follows, for example, from the classification of conjugacy classes in $\PU(2,1)$ (see Table~\ref{tab:conjugacy-classes}).
    \end{proofclaim}
    In our case, the minimality of $I$ implies that no $A_i$ is scalar. So, let us assume that $p_2$ is orthogonal to $p_3$ and $p_4$, permuting the indices if necessary. The matrix $A_2$ preserves $L_3, L_4$ and thus also $L_2 = \mathrm{span}\{p_3, p_4\}^\perp$. Hence, $L_2$ is preserved by all the matrices $A_i$ for $i\in I$, contradiction.
    \end{proof}

We are now ready to prove Proposition~\ref{prop:irreducible_many_points}.

\begin{proof}[Proof of Proposition~\ref{prop:irreducible_many_points}.]
    If there exist $1\leqslant i < j\leqslant n$ such that $\langle C_i, C_j\rangle$ is irreducible, then the result follows from Lemma~\ref{lem:two_irreducible}. We now assume that it is not the case. Proceeding by induction, with the case $n=4$ being immediate, it suffices to show that when $n\geq 5$, then $C_i = 1$ for some $i$. So, assume for sake of contradiction that $C_i\neq 1$ for all $i=1,\ldots,n$.
   
    Since $\rho$ is assumed to be irreducible, there exist $1\leq i < j < k \leq n$ such that $\langle C_i,C_j,C_k\rangle$ is irreducible by Lemma~\ref{lem:irreducible_four_sphere}. In other words, for any choice of lifts of $C_i,C_j,C_k$ to $\U(2,1)$, the induced subgroup of $\U(2,1)$ does not preserve any proper non-zero subspace $V\subset \mathbb C^3$. Let $S$ denote the $4$-punctured sub-sphere of $S_{0,n}$ whose peripheral curves are represented by $c_i$, $c_j$, $c_k$, and $(c_ic_jc_k)^{-1}$. The restriction of $\rho$ to $S$ is still an irreducible fixed point, and so, by Corollary~\ref{cor:irreducible_S04}, either some peripheral monodromy of $\rho|_{\pi_1 S}$ is trivial or $\rho|_{\pi_1 S}$ is conjugate to a Gassner fixed point $\rho_\alpha$. In the former case, since we are assuming that no peripheral monodromy is trivial, we would have $C_iC_jC_k=1$. However, that would imply $\langle C_i,C_j\rangle =\langle C_i,C_j,C_k\rangle$ is irreducible, contradicting our assumption. So, $\rho|_{\pi_1 S}$ is conjugate to a Gassner fixed point $\rho_\alpha$ and $C_iC_jC_k\neq 1$ is a pseudo-reflection.

    Consider first the case $n=5$, and assume that, up to changing the generators by a mapping class, that $(i,j,k)=(1,2,3)$ and $C_1C_2C_3$ is a pseudo-reflection. We consider the following 4-punctured sub-spheres of $S_{0,5}$:
    \begin{itemize}
        \item $S^{(1)}$ with peripheral curves $c_1$, $c_2c_3$, $c_4$, and $c_5$;
        \item $S^{(2)}$ with peripheral curves $c_1c_2c_3c_2^{-1}$, $c_2$, $c_4$, and $c_5$;
        \item $S^{(3)}$ with peripheral curves $c_1c_2$, $c_3$, $c_4$, and $c_5$.
    \end{itemize}
    The restrictions $\rho^{(\ell)}=\rho|_{\pi_1S^{(\ell)}}$ are fixed points for $\ell=1,2,3$. Since $c_2c_3$, $c_1c_2c_3c_2^{-1}$, and $c_1c_2$ are interior simple closed curves in $S$ and $\rho|_{\pi_1S}$ is a Gassner fixed point, their images by $\rho$ are non-trivial and not pseudo-reflections by Lemma~\ref{lem:gassner-fixed-points-interior-scc-are-RE-or-SP}. Therefore, $\rho^{(\ell)}$ cannot be a Gassner fixed point, and its kernel does not contain any peripheral curve. So, by Corollary~\ref{cor:irreducible_S04}, all the restrictions $\rho^{(\ell)}$ are reducible, hence pants-like by Theorem~\ref{thm:reducible-fixed-points-are-alomst-abelian}. However, using that $c_2c_3$, $c_1c_2c_3c_2^{-1}$, and $c_1c_2$ are interior simple closed curves in $S$ again and $\rho|_{\pi_1S}$ is a Gassner fixed point, we have, by Lemma~\ref{lem:gassner-fixed-points-interior-scc-are-RE-or-SP},
    \[
    [C_1,C_2C_3]\neq 1,\quad [C_1C_2C_3C_2^{-1},C_2]\neq 1,\quad [C_1C_2,C_3]\neq 1.
    \]
    This means that, for every $\ell=1,2,3$, we have either $C_4\in Z_{\rho^{(\ell)}}$ or $C_5\in Z_{\rho^{(\ell)}}$. By the pigeonhole principle, there exist $\ell_1\neq \ell_2$ and $m\in\{4,5\}$ such that $C_m\in Z_{\rho^{(\ell_1)}}\cap Z_{\rho^{(\ell_2)}}$, which implies that $C_m\in Z_\rho$, contradicting the irreducibility of $\rho$. 

    Now, let us consider the case $n\geq 6$. As before, we will assume that $(i,j,k)=(1,2,3)$. For every $\ell\geq 4$, consider the $5$-punctured sub-sphere $S^{(\ell)}$ whose peripheral curves are represented by $c_1$, $c_2$, $c_3$, $c_\ell$, and $(c_1c_2c_3c_\ell)^{-1}$. The restriction $\rho^{(\ell)}$ of $\rho$ to $\pi_1S^{(\ell)}$ is still a fixed point and $\rho^{(\ell)}$ is irreducible because its image contains $\langle C_1,C_2,C_3\rangle$. So, by our treatment of the case $n=5$ applied to $\rho^{(\ell)}$, we conclude that $C_1C_2C_3C_\ell =1$, or equivalently $C_\ell = (C_1C_2C_3)^{-1}$ for every $\ell\geq 4$. Now, for $\ell_2>\ell_1\geq 4$, consider the $5$-punctured sub-sphere $S$ of $S_{0,n}$ whose peripheral curves are $c_2$, $c_3$, $c_{\ell_1}$, $c_{\ell_2}$, and $(c_2c_3c_{\ell_1}c_{\ell_2})^{-1}$. The image of the restriction of $\rho$ to $\pi_1S$ is still irreducible because it contains $\langle C_2,C_3,C_{\ell_1}\rangle = \langle C_2,C_3,(C_1C_2C_3)^{-1}\rangle = \langle C_1,C_2,C_3\rangle$. So, the case $n=5$ implies $C_2C_3C_{\ell_1}C_{\ell_2}=1$, or equivalently $C_2=(C_3C_1^2)^{-1}$. This further implies that $\langle C_1,C_3\rangle = \langle C_1,C_2,C_3\rangle$ is irreducible, contradicting our assumption.
\end{proof}

\subsection{Proof of Theorem~\ref{thm:irreducible-fixed-points-genus-zero}}
\begin{proof}
    The statement follows directly from Corollary~\ref{cor:irreducible_S04} and Proposition~\ref{prop:irreducible_many_points}.
\end{proof}


\section{Mapping class group fixed points in genus at least one}\label{sec:fixed-points-genus-at-lesast-one}
\subsection{Overview}
The goal of this section is to classify all fixed points $\pi_1 S_{g,n}\to\PU(2,1)$ where $g\geq 1$ and $n\geq 0$. We will prove the following.
\begin{thm}\label{thm:fixed-points-genus-at-least-1}
    Let $\rho\colon\pi_1S_{g,n}\to\PU(2,1)$ be a fixed point with $g\geq 1$. Then one of the following holds:
    \begin{itemize}
        \item $g=1$, $n\geq 1$, and $\rho$ is irreducible, maps the peripheral loops around exactly $n-1$ peripheral curves to the identity, and factorizes through a Burau fixed point (Defintion~\ref{def:Burau-fixed-points-torus});
        \item $g=1$, $n\geq 1$, and $\rho$ is an extended quaternionic fixed point (Definition~\ref{def:extended-quaternionic-fixed-point});
        \item $\rho$ is the trivial representation.
    \end{itemize}
\end{thm}
\subsection{Once-punctured torus}\label{sec:once-punctured-torus}
We start by classifying all fixed points $\rho\colon\pi_1S_{1,1}\to\PU(2,1)$. As a consequence (Corollary~\ref{cor:fixed-points-closed-torus}), we obtain a classification of fixed points for the closed torus.
\subsubsection{The Burau fixed points}\label{sec:Burau-fixed-points-torus}
Because of the isomorphism $\PMod(S_{1,2})\cong B_4/Z(B_4)$ stated in Equation~\eqref{eq:PMod_1,2=B_4/Z(B_4)}, the reduced Burau representations $\beta'_{4,t}\colon B_4/Z(B_4)\to \PU(2,1)$ lead to mapping class group fixed points $\rho_t\colon\pi_1 S_{1,1}\to\PU(2,1)$ by restriction from $\PMod(S_{1,2})$ to $\pi_1S_{1,1}$ (by the same reasoning as in Lemma~\ref{lem:restriction-of-PMod-representations-are-fixed-points}). 
\begin{defn}\label{def:Burau-fixed-points-torus}
The representations $\rho_t\colon\pi_1 S_{1,1}\to\PU(2,1)$ obtained this way are called \emph{Burau fixed points}.
\end{defn}
Recall from Section~\ref{sec:Burau-fixed-points-sphere} that in order for the image of $\rho_t$ to belong to $\PU(2,1)$, we must have $t=e^{2i\theta}$ with $\theta\in (\pi/4,3\pi/4)$.

In order to work out the explicit form of $\rho_t$, we fix a presentation
\[
\pi_1S_{1,1} = \langle a,b,c : [a,b]=c\rangle.
\]
We will use the embedding $\PMod(S_{1,2})\subset \Aut(\pi_1S_{1,1})$ (whose explicit form is provided by Equation~\eqref{eq:embedding-PMod_1,2-in-Aut(pi_1S_1,1)}) to think of mapping classes in $\PMod(S_{1,2})$ as automorphism of $\pi_1 S_{1,1}$. The inclusion $\pi_1S_{1,1}\subset \PMod(S_{1,2})$ given by point-pushing maps $a$ and $b$ to the inner automorphisms of $\pi_1S_{1,1}$ given by conjugation by $a$, respectively $b$. In terms of the Dehn twist generators of $\PMod(S_{1,2})$ and using~\eqref{eq:embedding-PMod_1,2-in-Aut(pi_1S_1,1)}, we verify that
\[
a\mapsto \tau_3\tau_1^{-1},\quad b\mapsto (\tau_1\tau_2)(\tau_1\tau_3^{-1})(\tau_1\tau_2)^{-1}.
\]
Using the reduced Burau representation in the form~\eqref{eq:Burau-reduced-n=4-sage} with $t=s^2$, we obtain
\begin{align}\label{eq:burau-fixed-point-torus}
    \rho_t(a)&=\begin{pmatrix}
        -s^{-2} & s^{-1} & 0\\ 0 & 1 & 0\\ 0 & s & -s^2
    \end{pmatrix},\quad
    \rho_t(b)=\begin{pmatrix}
        1-s^{-2} & 0 & -s^2\\
        s-s^{-3} & -s^2 & 0\\
        -s^{-4} & 0 & 0
    \end{pmatrix},\\
    \rho_t(c) &= \begin{pmatrix}
        s^{-4} & 0 & s^6-s^4+1-s^{-2}\\
        0 & s^{-4} & s^7-s^3+s-s^{-3}\\
        0 & 0 & s^8
    \end{pmatrix}.\nonumber
\end{align}
We now recall some classical properties of Burau fixed points (see also~\cite[Section~3]{FGG26}).
\begin{lem}\label{lem:properties-Burau-fixed-points-torus}
    The Burau fixed points have the following properties. Recall that $t=e^{2i\theta}$ for $\theta\in (\pi/4,3\pi/4)\setminus\{\pi/2\}$.
    \begin{enumerate}
        \item For all $t$, $\rho_t(a)$ and $\rho_t(b)$ are regular elliptic.
        \item If $\theta\in \{\pi/3, 2\pi/3\}$, then $\rho_t(c)$ is a vertical translation.
        \item If $\theta\in (\pi/4,\pi/3)\cup (2\pi/3,3\pi/4)$, then $\rho_t(c)$ is a scalar rotation.
        \item If $\theta\in (\pi/3, 2\pi/3)\setminus\{\pi/2\}$, then $\rho_t(c)$ is a complex reflection.
    \end{enumerate}
\end{lem}
\begin{proof}
    Note that $\rho_t(a)$, $\rho_t(b)$, and $\rho_t(c)$ have determinant $1$ in the form~\eqref{eq:burau-fixed-point-torus}. We observe that the traces of $\rho_t(a)$ and $\rho_t(b)$ are equal to $1-s^2-s^{-2}=1-2\cos(2\theta)$ and therefore lies in the interior of Goldman's deltoid as long as $\theta\neq k\pi/2$. This shows that $\rho_t(a)$ and $\rho_t(b)$ and regular elliptic. The nature of $\rho_t(c)$ can be determined by analogous computations as in Lemma~\ref{lem:properties-Burau-fixed-points}.
\end{proof}

\begin{lem}\label{lem:Burau-useful}
    Let $\rho\colon\pi_1S_{1,1}\to\PU(2,1)$ be a Burau fixed point and let $a,b\in\pi_1S_{1,1}$ representing simple closed curves intersecting once. Then, $\rho(a)$ does not commute with $\rho(ba^{-1}b^{-1})$.
\end{lem}
\begin{proof}
    By Lemma~\ref{lem:properties-Burau-fixed-points-torus}, $\rho(a)$
    is regular elliptic. So, if $\rho(a)$ commutes with $\rho(ba^{-1}b^{-1})$, then $\rho(b)$ must fix the unique fixed point $X$ of $\rho(a)$ (see Table~\ref{tab:conjugacy-classes}). However, this would imply $\Im(\rho)\subset \stab(X)$, contradicting the irreducibility of $\rho$.
\end{proof}
\subsubsection{Irreducible fixed points}
Irreducible fixed points $\pi_1S_{1,1}\to\PU(2,1)$ always lift to homomorphisms $\PMod(S_{1,2})\to\PU(2,1)$ for the same reason as in Lemma~\ref{lem:regular-and-abelian-fixed-points-are-restrictions}. The following result is thus a special case of Formanek's classification~\cite{formanek}.

\begin{lem}\label{lem:irreducible-fixed-point-once-punctured-torus}
    Let $S=S_{1,1}$ be a once-punctured torus and $\rho\colon \pi_1 S\to \PU(2,1)$ be an irreducible fixed point. Then, $\rho$ is a Burau fixed point.
\end{lem}
\begin{proof}
    Let $\iota$ be the hyperelliptic involution of $S$, i.e.~the only non-trivial central element in $\Mod_{1,1}$. It induces an injection $f\colon \pi_1 S\to \pi_1^\mathrm{orb}(S / \iota)$, where $S / \iota$ is a sphere with one puncture and three orbifold singularities of order $2$. 
    We fix presentations
    \begin{align*}
        \pi_1 S &= \langle a, b, c : [a,b] = c\rangle\\
        \pi_1^\mathrm{orb}(S / \iota) &= \langle q_1, q_2, q_3, q_4 : q_2^2 = q_3^2 = q_4^2 = q_1q_2q_3q_4 = 1\rangle
    \end{align*}
    such that $f$ is defined by
\begin{align*}
    a &\mapsto q_2q_3\\
    b &\mapsto q_4 q_3\\
    c&\mapsto q_1^{-2}.
\end{align*}
    The map $\iota_*\colon \pi_1 S\to \pi_1 S$ is given by $a\mapsto a^{-1}$ and $b\mapsto b^{-1}$. It corresponds to the conjugation by $q_3$.
    Since $\rho$ is invariant and irreducible, it extends uniquely to a representation $\widetilde \rho\colon \pi_1^\mathrm{orb}(S / \iota)\to \mathrm{PU}(2,1)$ under which $q_3$ is mapped to the unique $g\in \mathrm{PU}(2,1)$ such that $g\left (\rho\circ \iota_*\right ) g^{-1} = \rho$.
    Let $p\colon \pi_1 S_{0,4} = \langle q_1, q_2, q_3, q_4 : q_1q_2q_3 q_4 = 1\rangle \to \pi_1^\mathrm{orb}(S / \iota)$ be the obvious map and $\bar \rho = \widetilde \rho\circ p$.
    The representation $\bar \rho$ is fixed by the subgroup of $\mathrm{Mod}_{0,4}$ stabilising the puncture corresponding to $q_1$ and isomorphic to $\Mod_{1,1}/\langle \iota\rangle$, hence is also stabilized by $\PMod_{0,4}$.
    Since it is irreducible (the image of $\tilde\rho$ contains the image of $\rho$), by Corollary~\ref{cor:irreducible_S04}, there exists $1\leq i \leq 4$ such that $\bar \rho(q_i) = 1$ or $\bar\rho$ is conjugate to a Gassner fixed point $\rho_\alpha$.
    If there exists $2\leqslant i\leqslant 4$, such that $\bar \rho(q_i) = 1$, then $\bar \rho(q_2) = \bar \rho(q_3) = \bar\rho(q_4) = 1$, and $\bar \rho$ is trivial. 
    Similarly, if $\bar \rho(q_1) = 1$, then $\rho$ is abelian since $\bar \rho(q_1)^2 = \rho(c)^{-1}$.
    Therefore, there exists $\alpha = \left ( z_1, \ldots, z_5 \right )\in \Omega$ such that $\bar \rho$ is conjugate to $\rho_\alpha$.
    Since the $\bar \rho(c_i)$ are involutions for $2\leqslant i\leqslant 4$, we have $z_iz_5 = -1$ and $z_2 = z_3 = z_4$.
    In particular, $\alpha$ is of the form $\alpha = (-t^2, t^{-1}, t^{-1}, t^{-1}, -t)$.
    It follows that 
    \[
    \rho(a) = \begin{pmatrix}
1 & 1+t^2 & 0\\
0 & -1 & 0\\
0 & -t+1 & 1 \end{pmatrix} \begin{pmatrix}
1 & 0 & t + t^{-1}\\
0 & 1 & 1-t^{-1}\\
0 & 0 & -1
\end{pmatrix} = \begin{pmatrix}
1 & 1+t^2 & 1+t^2\\
0 & -1 & -1+t^{-1}\\
0 & 1-t & -t+1-t^{-1}
\end{pmatrix}.
    \]
Similarly, we compute
    \[\rho(b) = \begin{pmatrix}
-t^{-2} & -1-t^{-2} & 0\\
-t^{-1}+t^{-2} & 1-t^{-1}+t^{-2} & 0\\
-1+t^{-1} & -1+t^{-1} & -t
\end{pmatrix}.\]
Now take $s$ such that $t = s^{2}$ and consider \[P = \begin{pmatrix}
    0 & 0 & -1 - s^4\\
    -1 & 0 & 1\\
    1 & -s & s^2
\end{pmatrix}.\]
The matrix $P$ is invertible if $t\neq \pm i$. But we cannot have $t = \pm i$, since $\rho$ is then reducible.
We check that $\rho(a) = P \rho_t(a) P^{-1}$ and $\rho(b) = P\rho_t(b)P^{-1}$, where $\rho_t$ is given by~\eqref{eq:burau-fixed-point-torus}. This means that $\rho$ is a Burau fixed point (see Section~\ref{sec:Burau-is-Gassner}).
\end{proof}

\subsubsection{Reducible fixed points}
We will show that there are two kinds of reducible fixed points $\pi_1S_{1,1}\to\PU(2,1)$: the trivial representation and a representation valued in the quaternion group $H_8\subset \mathrm{P}(\U(2)\times \U(1))\subset\PU(2,1)$:
\begin{equation}\label{eq:quaternion-group}
H_8=\{1,\sigma, I,J,K,\sigma I, \sigma J,\sigma K\},
\end{equation}
where
\[
\sigma=\begin{pmatrix}
    -1 & 0 & 0\\
    0 & -1 & 0\\
    0 & 0 & 1
\end{pmatrix}, \quad I = \begin{pmatrix}
    i & 0 & 0\\
    0 & -i & 0\\
    0 & 0 & 1
\end{pmatrix}, \quad J= \begin{pmatrix}
    0 & i & 0\\
    i & 0 & 0\\
    0 & 0 & 1
\end{pmatrix}, \quad K = \begin{pmatrix}
    0 & - 1 & 0\\
    1 & 0 & 0\\
    0 & 0 & 1
\end{pmatrix}.
\]
\begin{defn}\label{def:quaternionic}
    A \emph{quaternionic fixed point} is any representation conjugate to $\rho\colon\pi_1S_{1,1}\to\PU(2,1)$ defined by
    \begin{equation}\label{eq:quaternionic-representation}
        \rho(a)=I,\quad
        \rho(b)=J.
    \end{equation}
\end{defn}

\begin{lem}\label{lem:properties-quaternionic-fixed-point}
    The conjugacy class of a quaternionic fixed point is indeed fixed by $\PMod(S_{1,1})$. Furthermore, their peripheral monodromies are scalar rotations of angle $\pi$.
\end{lem}
\begin{proof}
    If $\rho$ is quaternionic, then, using~\eqref{eq:quaternionic-representation}, we compute
    \[
    \rho([a,b])=\begin{pmatrix}
        -1 & 0 & 0\\
        0 & -1 & 0\\
        0 & 0 & 1
    \end{pmatrix},
    \]
    which is a scalar rotation of angle $\pi$. Now, according to~\cite[Figure~4.10]{farb-margalit}, $\PMod(S_{1,1})$ is generated by two Dehn twists represented by the automorphisms of $\pi_1S_{1,1}$ given by
    \begin{equation}\label{eq:twists-PMod_1,1}
    T_a \colon \begin{cases}
    a\mapsto a\\
    b\mapsto ba^{-1},
    \end{cases} \quad 
    T_b \colon \begin{cases}
    a\mapsto ab\\
    b\mapsto b.
    \end{cases}
    \end{equation}
    A simple computation gives $(\rho\circ T_a) (b) =(\rho\circ T_b)(a) = K$, as well as $X(\rho\circ T_a)X^{-1}=\rho$ and $Y(\rho\circ T_b)Y^{-1}=\rho$, where
    \[
        X=\begin{pmatrix}
            e^\frac{3i\pi}{4} & 0 & 0\\
            0 & e^\frac{5i\pi}{4} & 0\\
            0 & 0 & 1
        \end{pmatrix}, \quad
        Y=\begin{pmatrix}
            \frac{-1}{\sqrt{2}}&\frac{-i}{\sqrt{2}}&0\\
            \frac{-i}{\sqrt{2}}&\frac{-1}{\sqrt{2}}&0\\
            0&0&1
        \end{pmatrix}.\qedhere
    \]
\end{proof}
\begin{rem}
    Every quaternionic fixed point $\rho\colon\pi_1S_{1,1}\to H_8$ extends to $\widehat{\rho}\colon B_4/Z(B_4)\rightarrow \PU(2,1)$ by
    \[
    \sigma_1\mapsto X,\quad \sigma_2\mapsto Y,\quad \sigma_3\mapsto IX,
    \]
    where $\sigma_1,\sigma_2,\sigma_3$ are Artin's generators of $B_4$. Note that the image of $\widehat{\rho}$ is isomorphic to the binary octahedral group $\widetilde{S}_4$ of cardinal $48$. The representation $\widehat{\rho}$ is accounted for in Formanek's classification in~\cite[Theorem~12~(a)]{formanek} (with $y=e^{5i\pi/4}$ and $z=i$).    
\end{rem}

\begin{prop}\label{prop:reducible-fixed-point-once-punctured-torus}
    Let $\rho\colon\pi_1S_{1,1}\to\PU(2,1)$ be a reducible fixed point. Then, either $\rho$ is the trivial representation, or $\rho$ is a quaternionic fixed point.
\end{prop}
\begin{proof}
    Since $\rho$ is a fixed point and $\PMod(S_{1,1})$ acts transitively on oriented non-peripheral simple closed curves, they are mapped to the same conjugacy class by $\rho$. We will show that this conjugacy class $\mathcal{C}$ must be trivial, unless $\rho$ is quaternionic.
    \begin{claim}\label{claim:no-SP-no-VT}
        The class $\mathcal{C}$ does not contain screw-parabolic elements or vertical translations.
    \end{claim}
    \begin{proofclaim}
        If $\gamma\in\pi_1S_{1,1}$ represents a non-peripheral simple
        closed curve, so does $\gamma^{-1}$. However, screw-parabolic elements and vertical translations are not conjugate to their inverse.
    \end{proofclaim}
    \begin{claim}\label{claim:no-PR}
        The class $\mathcal{C}$ does not contain pseudo-reflections.
    \end{claim}
    \begin{proofclaim}
    As $a$, $b$, $b^{-1}a^{-1}$ all represent non-peripheral simple closed curves and multiply to $1$, they can only be mapped to pseudo-reflections if their images commute by Lemma~\ref{lem:product-pseudo-reflections}.
    Moreover, their images must be conjugate to their inverses, so they are either all complex reflections of angle $\pi$
    or all scalar rotations of angle $\pi$.
    However, the product of two commuting scalar rotations of angle $\pi$ is trivial, while the product of two commuting complex reflections of angle $\pi$ is either trivial or a scalar rotation of angle $\pi$. This leads to a contradiction in both cases.
    \end{proofclaim}
    
    It remains to study when non-peripheral simple closed curves are all mapped to hyperbolic, regular parabolic or regular elliptic elements. Since $\rho$ is reducible, the image of $\rho$ is contained in $\stab(X)$ for $X\in\CH^2$, $\stab(Y)$ for $Y\in\partial\CH^2$, or in $\stab(L)$ for a $\C$-line $L$ by Lemma~\ref{lem:maximal-reducible-subgroups-of-PU(2,1)}.
    \begin{claim}\label{claim:stab(Y)-implies-trivial}
        If $\Im(\rho)\subset\stab(Y)$, then $\rho$ is trivial.
    \end{claim}
    \begin{proofclaim}
        We may assume that $\rho$ is upper triangular. Since the images under $\rho$ of $a$, $b$, $ba^n$ (that all represent non-peripheral simple closed curves) are all conjugate, their eigenvalues must be of modulus $1$. Since they cannot be screw-parabolic, nor vertical translations, nor complex reflections by Claims~\ref{claim:no-SP-no-VT} and~\ref{claim:no-PR}, they must be regular parabolic or trivial. Assuming $\rho$ is not trivial, we write $\rho(\gamma)=A(\zeta_\gamma, v_\gamma)$ using the notation~\eqref{eq:Heisenberg-translations} for $\gamma\in\{a,b,ba^{-1}\}$. Recall that $Z(A(\zeta,v))$ consists of all $A(\lambda\zeta,t)$ for $\lambda\in\R$. Since $\rho$ is conjugate to $\rho\circ T_a$, where $T_a$ is the automorphism from~\eqref{eq:twists-PMod_1,1}, $A(\zeta_b,v_b)$ is conjugate to $A(\zeta_{ba^{-1}},v_{ba^{-1}})$ by an element of the form $A(\lambda\zeta_a,v_a)$, forcing $\zeta_a=0$ because $\zeta_{ba^{-1}}=\zeta_b-\zeta_a$. So $\rho(a)$ is a vertical translation, contradiction.
    \end{proofclaim}
    \begin{claim}
        If $\Im(\rho)\subset\stab(X)$, then $\rho$ is trivial or a quaternionic fixed point.
    \end{claim}
    \begin{proofclaim}
        Since Claim~\ref{claim:no-PR} excludes pseudo-reflection, assuming $\rho$ is not trivial, we may assume that $\mathcal{C}$ consists of regular elliptic elements. Since $\rho(a)$ is conjugate to its inverse and not a pseudo-reflection, it is of the form $E(\theta,2\pi-\theta)$ with $\theta\in(0,\pi)$, using the notation~\eqref{eq:elliptic-elements}. As in Claim~\ref{claim:stab(Y)-implies-trivial}, $\rho(b)$ and $(\rho\circ T_a)(b)=\rho(ba^{-1})$ must be conjugate by an element commuting with $\rho(a)$, hence of the form $E(\alpha,\beta)$. This forces $\rho(b)$ to be of the form
        \[
        \rho(b)=\begin{pmatrix}
            0 & x & 0\\
            y & 0 & 0\\
            0 & 0 & 1
        \end{pmatrix}.
        \]
        Since $\rho(b)^2$ is conjugate to $\rho(a)^2$, we must have $e^{2i\theta}=e^{-2i\theta}$, forcing $\theta=\pi/2$ and $\rho(a)=I$. This also implies $xy=-1$. Conjugating $\rho$ by
        \[
        \begin{pmatrix}
            (-ix)^{-1/2} & 0 & 0\\
            0 & (-ix)^{1/2} & 0\\
            0 & 0 & 1
        \end{pmatrix},
        \]
        it becomes of the form~\eqref{eq:quaternionic-representation}, proving that $\rho$ is a quaternionic fixed point.
    \end{proofclaim}
    \begin{claim}
        If $\Im(\rho)\subset\stab(L)$, then $\rho$ is trivial.
    \end{claim}
    \begin{proofclaim}
        Claims~\ref{claim:no-SP-no-VT} and~\ref{claim:no-PR} exclude pseudo-reflections and screw-parabolic, assuming $\rho$ is not trivial, we may assume that $\mathcal{C}$ consists of regular elliptic elements or of hyperbolic elements.
        
        If $\rho(a)$, $\rho(b^{\pm 1})$, and $\rho(ab^{\pm 1})$ are all regular elliptic, they all act on $L$ as rotations of angle either $\theta$ or $2\pi-\theta$ for some $\theta\in(0,\pi)$ only depending on $\mathcal{C}$. Let $\varepsilon\in\{\pm 1\}$ such that the rotation angle of $\rho(a)$ on $L$ is opposite to that of $\rho(b^\varepsilon)$. Note that $\rho(a)$ and $\rho(b^\varepsilon)$ cannot fix the same point, for otherwise $\rho(ab^\varepsilon)$ is trivial or a complex reflection. Therefore, the restriction of $\rho(ab^\varepsilon)$ to $L$ is the product of two rotations around different points with opposite angles. This implies that $\rho(ab^\varepsilon)$ is not elliptic, contradiction.

        If $\rho(a)$ and $\rho(b)$ are hyperbolic, we may assume that they fix different points, for otherwise Claim~\ref{claim:stab(Y)-implies-trivial} applies. In that case, $\rho(ba^n)$ is hyperbolic for all $n\geq 1$ and its translation length cannot be the same for all $n\geq 1$, contradiction.
        \end{proofclaim}
        This finishes the proof.
\end{proof}

\begin{cor}\label{cor:fixed-points-closed-torus}
    If $\rho\colon\pi_1S_{1,0}\to\PU(2,1)$ is a fixed point, then $\rho$ is trivial.
\end{cor}
\begin{proof}
    Note that $\rho$ is necessarily reducible. Since $\PMod(S_{1,0})=\PMod(S_{1,1})$, $\rho$ must be one of the reducible fixed points of the once-punctured torus with trivial peripheral monodromy. By Proposition~\ref{prop:reducible-fixed-point-once-punctured-torus} and Lemma~\ref{lem:properties-quaternionic-fixed-point}, $\rho$ may only be the trivial representation.
\end{proof}

\subsection{Tori with at least two punctures}\label{sec:twice-punctured-torus}
We now move on to the classification of all fixed points $\rho\colon\pi_1S_{1,n}\to\PU(2,1)$ when $n\geq 2$. We fix a presentation
\[
\pi_1 S_{1,n}=\langle a,b,c_1,\ldots,c_n : [a,b]=\prod c_i\rangle.
\]
We will use the following terminology.
\begin{defn}\label{def:extended-quaternionic-fixed-point}
    An \emph{extended quaternionic fixed point} is any representation conjugate to a representation $\rho\colon\pi_1S_{1,n}\to\PU(2,1)$ of the form
    \begin{equation}\label{eq:extended-quaternionic}
    \rho(a)= I,\quad \rho(b)= J, \quad \rho(c_i)\in \{1,\sigma\},
    \end{equation}
    where $I,J,\sigma$ are the elements of the quaternion group $H_8$ from~\eqref{eq:quaternion-group}. (Note that $\rho$ must map an odd number of $c_i$ to $\sigma$ in order to be a morphism.)
\end{defn}

\begin{rem}
    A representation $\rho\colon\pi_1S_{1,n}\to\PU(2,1)$
    is an extended quaternionic fixed point if and only if, up to conjugation, $\Im(\rho)=H_8$ and $\rho$ maps all peripheral curves to $Z(H_8)=\{1,\sigma\}$.
\end{rem}

\begin{lem}\label{lem:extended-quaternionic-are-fixed}
    The conjugacy class of an extended quaternionic fixed point is indeed fixed by $\PMod(S_{1,n})$.
\end{lem}
\begin{proof}
    A family of Dehn twists generators of $\PMod(S_{1,n})$ is given by the Dehn twists along the simple closed curves represented by $a,b,ac_1^{-1},\ldots,ac_n^{-1}$ (see~\cite[Figure~4.10]{farb-margalit}). These Dehn twists are represented by the following automorphisms of $\pi_1S_{1,n}$.
    \begin{equation}\label{eq:twists-PMod_1,n}
    T_a \colon \begin{cases}
    a\mapsto a\\
    b\mapsto ba^{-1}\\
    c_i\mapsto c_i,
    \end{cases} \quad 
    T_b \colon \begin{cases}
    a\mapsto ab\\
    b\mapsto b\\
    c_i\mapsto c_i,
    \end{cases} \quad
    T_{ac_i^{-1}} \colon \begin{cases}
    a\mapsto a\\
    b\mapsto a^{-1}c_ib\\
    c_j\mapsto c_ia^{-1}c_jac_i^{-1} & \text{if }j<i\\
    c_i\mapsto c_i\\
    c_j \mapsto a^{-1}c_ic_jc_i^{-1}a & \text{if }j>i.
    \end{cases}
    \end{equation}
    Given the extended quaternionic fixed point $\rho\colon\pi_1S_{1,n}\to\PU(2,1)$ from~\eqref{eq:extended-quaternionic}, since $1$ and $\sigma$ are central in $H_8$, we have $\rho(c_i)\in Z_\rho$ for every $i$. In particular, $\rho\circ T_a$ and $\rho\circ T_b$ are both conjugate to $\rho$ by the same element as in the proof of Lemma~\ref{lem:properties-quaternionic-fixed-point}. Furthermore, note that 
    \[
    \rho\circ T_{ac_i^{-1}}\colon \begin{cases}
    a\mapsto \rho(a)=I\\
    b\mapsto \rho(a^{-1}bc_i)\in \{K,\sigma K\}\\
    c_j\mapsto \rho(c_j) \in \{1,\sigma\}.
    \end{cases}
    \]
    Since $I(\sigma K)I^{-1}=K$, we see that $\rho\circ T_{ac_i^{-1}}$ is always conjugate to $\rho$ (either by $X$ or by $XI$, in the notation of Lemma~\ref{lem:properties-quaternionic-fixed-point}).
\end{proof}

\begin{prop}\label{prop:fixed-points-torus}
    Let $\rho\colon\pi_1S_{1,n}\to\PU(2,1)$ be a fixed point with $n\geq 2$. Then either:
    \begin{itemize}
        \item $\rho$ is the trivial representation;
        \item $\rho$ is irreducible, maps the peripheral curves around exactly $n-1$ punctures to the identity, and factorizes through a Burau fixed point;
        \item $\rho$ is an extended quaternionic fixed point.
    \end{itemize}
\end{prop}
We divide the proof of Proposition~\ref{prop:fixed-points-torus} into two parts, studying the nature of the restriction to a once-punctured sub-torus. Note that $\PMod(S_{1,n})$ acts transitively on separating simple closed curves in $S_{1,n}$ that separates it into a once-punctured sub-torus and an $(n+1)$-punctured sphere. In particular, the restrictions of a fixed point $\rho\colon S_{1,n}\to\PU(2,1)$ to any once-punctured sub-torus are all conjugate.

\begin{lem}\label{lem:irred-multipunctured-tori}
    If the restriction of $\rho$ to one (equivalently any) once-punctured sub-torus is irreducible, then $\rho$ maps the peripheral curves around exactly $n-1$ punctures to the identity and restricts to a Burau fixed point on each once-punctured sub-torus.
\end{lem}
\begin{proof}
    We will make use of the simple closed  curves represented by $a,b,ab,[a,b]$.
    Cutting $S_{1,n}$ along $[a,b]$ separates it into a once-punctured sub-torus $S^{(1)}=S_{1,1}$ and an $(n+1)$-punctured sphere $S^{(2)}=S_{0,n+1}$. Since the restriction of $\rho$ to $\pi_1 S^{(1)}$ is an irreducible fixed point by assumption, it is a Burau fixed point by~\ref{lem:irreducible-fixed-point-once-punctured-torus}. Therefore, by Lemma~\ref{lem:properties-Burau-fixed-points-torus}, $\rho([a,b])$ is a pseudo-reflection and $\rho(a)$, $\rho(b)$, and $\rho(ab)$ are regular elliptic. 
    
    Now, we cut $S_{1,n}$ along the simple closed curve represented by $\gamma\in\{a,b,ab\}$ to obtain an $(n+2)$-punctured sub-sphere $S_\gamma$ with peripheral curves $\gamma,\gamma',c_1,\ldots,c_n$ (where $\gamma'$ is conjugate to $\gamma^{-1}$). Since $\rho(\gamma)$ is regular elliptic, the restriction of $\rho$ to $\pi_1S_\gamma$ has to be pants-like by Theorems~\ref{thm:reducible-fixed-points-are-alomst-abelian} and~\ref{thm:irreducible-fixed-points-genus-zero} (because Gassner fixed points cannot have regular elliptic peripheral monodromies by Lemma~\ref{lem:Gassner-peripheral-monodromies-are-pseudo-refelections}). So, $n-1$ among the peripheral curves $\gamma,\gamma',c_1,\ldots,c_n$ are mapped into $Z(\Im(\rho|_{\pi_1 S_\gamma}))$.
    By Lemma~\ref{lem:Burau-useful}, $\rho(\gamma)$ and $\rho(\gamma')$ do not commute, so $\gamma$ and $\gamma'$ are not among these $n-1$ curves. This means that the images of $n-1$ peripheral curves among $c_1,\ldots,c_n$ commute with all $\rho(c_i)$ and $\rho(\gamma)$.

    We apply the previous observation $\gamma=a$ and $\gamma=b$, which implies that the images of at least $n-2$ peripheral curves among $c_1,\ldots,c_n$ commute with all $\rho(c_i)$ and with $\rho(a)$ and $\rho(b)$. The images of those curves therefore lie in $Z_\rho=1$.
    Up to permuting the punctures, assume that $\rho(c_3)=\cdots=\rho(c_n)=1$. Furthermore, we also know that among $\rho(c_1)$ and $\rho(c_2)$, at least one of those commutes with two among $\rho(a),\rho(b),\rho(ab)$, forcing it to be trivial.
\end{proof}
We now consider the fixed points that restrict to a quaternionic fixed point.
\begin{lem}\label{lem:quaternionic-extended}
    If the restriction of $\rho$ to one (equivalently any) once-punctured sub-torus is a quaternionic fixed point, then $\rho$ is an extended quaternionic fixed point.
\end{lem}
\begin{proof}
    We consider again the simple closed curves $a,b,ab,[a,b]$ and the sub-surfaces $S^{(1)}=S_{1,1}$ and $S^{(2)}=S_{0,n+1}$ obtained by cutting $S_{1,n}$ along $[a,b]$. Since the restriction of $\rho$ to $\pi_1S^{(1)}$ is a quaternionic fixed point, $\rho([a,b])$ is a scalar rotation of angle $\pi$ by Lemma~\ref{lem:properties-quaternionic-fixed-point}. Up to conjugating $\rho$, we may assume that $\rho(a)=I$, $\rho(b)=J$, and $\rho([a,b])=\sigma$.
    
    Cutting $S_{1,n}$ along the simple closed curve represented by $a$, we obtain an $(n+2)$-punctured sub-sphere $S_a=S_{0,n+2}$ with peripheral curves $bab^{-1},a^{-1},c_1,\ldots,c_n$. Since $\rho(a)=I$ is regular elliptic and the restriction of $\rho$ to $\pi_1S_a$ is a fixed point, it must be pants-like (for the same reason as in Lemma~\ref{lem:irred-multipunctured-tori}). So, the images of $n-3$ peripheral curves among $bab^{-1},a^{-1},c_1,\ldots,c_n$ belong to $Z(\Im(\rho|_{\pi_1S_a}))$.

    Assume first that $\rho(a)=I$ commutes with all $\rho(c_i)$. Since $I$ is regular elliptic, $\rho(c_i)$ are of form $E(\alpha_i,\beta_i)$ (using the notation~\eqref{eq:elliptic-elements}). Since $ac_i^{-1}$ and $a$ both represent non-separating simple closed curves on $S_{1,n}$ and since $\PMod(S_{1,n})$ acts transitively on non-separating simple closed curves, $\rho(ac_i^{-1})=I\rho(c_i^{-1})$ and $\rho(a)=I$ are conjugate. A simple computation shows that this is only possible if $(\alpha_i,\beta_i)\in\{(0,0),(\pi,\pi)\}$, or equivalently $\rho(c_i)\in\{1,\sigma\}$. In order for $\rho$ to be a morphism, an odd number of $\rho(c_i)$ are equal to $\sigma$, showing that $\rho$ is an extended quaternionic fixed point.
    
    If $\rho(a)=I$ does not commute with all $\rho(c_i)$,
    then neither does $\rho(bab^{-1})=\sigma I$. So, $n-1$ among $\rho(c_1),\ldots,\rho(c_n)$ must commute with $\rho(a)$ and all $\rho(c_i)$. Up to permuting the punctures, assume that they are $\rho(c_2),\ldots,\rho(c_n)$. The same reasoning as before shows that  $\rho(c_2),\ldots,\rho(c_n)\in\{1,\sigma\}$. Since $\rho([a,b])=\sigma$, it follows that all $\rho(c_1)\in\{1,\sigma\}$,
    which contradicts the fact that $\rho(a)$ does not commute with all $\rho(c_i)$.
\end{proof}

\begin{proof}[Proof of Proposition~\ref{prop:fixed-points-torus}]
    Let $\rho\colon\pi_1 S_{1,n}\to\PU(2,1)$ be a fixed point. Let $S_{1,1}\subset S_{1,n}$ be a once-punctured sub-torus. The restriction of $\rho$ to $\pi_1 S_{1,1}$ is either trivial, a Burau fixed point, or a quaternionic fixed point.
    
    If the restriction is trivial, then $\rho$ maps a non-separating simple closed curve in $S_{1,n}$ to the identity. Since $\PMod(S_{1,n})$ acts transitively on non-separating simple closed curves and $\rho$ is a fixed point, they all have conjugate images, hence trivial. In particular, all of $a,b,ac_1^{-1},\ldots,ac_n^{-1}$ have trivial images and $\rho$ is trivial.
    
    If the restriction is a Burau fixed point, then all but one peripheral monodromy of $\rho$ are trivial by Lemma~\ref{lem:irred-multipunctured-tori}.

    If the restriction is a quaternionic fixed point, $\rho$ is an extended quaternionic fixed point by Lemma~\ref{lem:quaternionic-extended}. 
\end{proof}

\subsection{Surfaces of genus at least two}\label{sec:higher-genus}
As we are about to see, if $g\geq 2$ then all fixed points are trivial.
\begin{prop}\label{prop:fixed-points-genus-at-least-two}
    If $\rho\colon\pi_1 S_{g,n}\to\PU(2,1)$ is a fixed point with $g\geq 2$, then $\rho$ is the trivial representation.
\end{prop}
\begin{proof}
We fix a presentation
\[
\pi_1S_{g,n}=\langle a_1,b_1,\ldots,a_g,b_g,c_1,\ldots,c_n : \prod [a_i,b_i]=\prod c_j\rangle.
\]
Cutting $S_{g,n}$ along the simple closed curves represented by $a_1,\ldots,a_{g-1}$, we obtain a sub-surface $S'=S_{1,n+2(g-1)}$. Note that since $g\geq 2$, $S'$ is a torus with at least two punctures.
Since the restriction of $\rho$ to $\pi_1S'$ is still a fixed point, it is either trivial, factorizes through a Burau fixed point, or is an extended quaternionic fixed point by Proposition~\ref{prop:fixed-points-torus}.

If the restriction is trivial, then $\rho(a_1)=1$ and since $\PMod(S_{g,n})$ acts transitively on non-separating simple closed curves, all of $a_1,\ldots,a_g,b_1,\ldots,b_g,ac_1^{-1},\ldots,ac_n^{-1}$ have trivial image, proving that $\rho$ is trivial.

If the restriction factorizes through a Burau representation, then all but one of its peripheral monodromies are trivial. In particular $\rho(a_1)=1$, and we conclude that $\rho$ is trivial as before.

If the restriction is an extended quaternionic fixed point, then $\rho(a_1)=1$ or $\rho(a_1)=\sigma$ (up to conjugating $\rho$). Furthermore, $\rho(a_g)$ is regular elliptic by Lemma~\ref{lem:properties-quaternionic-fixed-point}, so $\rho(a_1)$ is not conjugate to $\rho(a_g)$, contradiction.
\end{proof}

\subsection{Proof of Theorem~\ref{thm:fixed-points-genus-at-least-1}}
\begin{proof}
    The statement follows directly from Proposition~\ref{prop:reducible-fixed-point-once-punctured-torus}, Corollary~\ref{cor:fixed-points-closed-torus}, Proposition~\ref{prop:fixed-points-torus}, and Proposition~\ref{prop:fixed-points-genus-at-least-two}.
\end{proof}

\appendix
\setcounter{section}{0}
\renewcommand\sectionname{Appendix}
\addcontentsline{toc}{section}{Appendix}
\renewcommand{\thesection}{\Alph{section}}

\section{Pseudo-reflections and rigidity}\label{apx:horn-pseudo-reflections}
A consequence of our classification is that all the irreducible fixed points $\pi_1S_{0,4}\to\PU(2,1)$ are rigid (Theorem~\ref{thm:intro-fixed-point-iff-rigid-spheres}). The purpose of this appendix is to give geometric intuition and a sketch of proof that irreducible representations of $4$-punctured spheres with pseudo-reflections as peripheral monodromies are automatically rigid, hence fixed points.
\subsection{Interior simple closed curves}
\begin{lem}\label{lem:interior-scc-are-regular-elliptic-or-screw-parabolic}
    Let $\rho\colon\pi_1S_{0,4}\to\PU(2,1)$ be a regular representation (i.e.~$Z_\rho=\{1\}$) that maps all peripheral curves to pseudo-reflections.
    Then $\rho$ maps every non-peripheral simple closed curve
    to regular elliptic or screw parabolic elements.
\end{lem}
\begin{proof}
    Let $\gamma$ be a non-peripheral simple closed curve in $S_{0,4}$. There is a presentation of $\pi_1S_{0,4}$ as
    \[
    \pi_1 S_{0,4}= \langle c_1,c_2,c_3,c_4\colon \prod c_i=1\rangle
    \]
    such that $\gamma$ is represented by the fundamental group element $c_1c_2$. We will proceed by contradiction and assume that $\rho(c_1c_2)$ is neither regular elliptic nor screw parabolic.

    First, assume that $\rho(c_1c_2)=1$. In that case, the image of $\rho$ is the subgroup of $\PU(2,1)$ generated by $\rho(c_1)$ and $\rho(c_3)$. Both being pseudo-reflections, the intersection of their centralizers is non-trivial by Corollary~\ref{cor:intersection-centralizers-pseudo-reflections-non-trivial}, which contradicts the regularity of $\rho$.

    Next, assume that $\rho(c_1c_2)\neq 1$ is a pseudo-reflection. In that case, Lemma~\ref{lem:product-pseudo-reflections} would imply that $\rho(c_1)$ and $\rho(c_2)$ commute. However, $\rho(c_1c_2)$, $\rho(c_3)$
    and $\rho(c_4)^{-1}$ are pseudo-reflections, and so $\rho(c_3)$ and $\rho(c_1c_2)$ would also commute. This would imply $\rho(c_1c_2)\in Z_\rho$, contradicting the regularity of $\rho$.

    Now, assume that $\rho(c_1c_2)$ is regular parabolic. In particular, $\rho(c_1c_2)$ does not stabilize any $\C$-line. This means that $\rho(c_1)$ and $\rho(c_2)$ cannot be vertical translations or scalar rotations, as their product with a pseudo-reflection always stabilizes a $\C$-line. If $\rho(c_1)$ and $\rho(c_2)$ are complex reflections with respective mirrors $M_1$ and $M_2$, then $M_1$ and $M_2$ must intersect at a unique point $Y\in\partial\CH^2$ in order for $\rho(c_1c_2)$ to be parabolic. The same conclusion holds for $\rho(c_3)$ and $\rho(c_4)$ with the same boundary point $Y$, so that the image of $\rho$ stabilizes the horosphere foliation $H$ at $Y$. Hence, $\mathcal{VT}_Y\subset Z_\rho$ (see Lemma~\ref{lem:maximal-subgroups-of-PU(2,1)-with-non-trivial-centralizer}), contradicting the regularity of $\rho$.
    
    Finally, assume that $\rho(c_1c_2)$ is hyperbolic. In particular, $\rho(c_1c_2)$ stabilizes a unique $\C$-line $L$. If one of $\rho(c_1)$ or $\rho(c_2)$ is a vertical translation or a scalar rotation, then in order for $\rho(c_1c_2)$ to be hyperbolic, $\rho(c_1)$ and $\rho(c_2)$ must both stabilize $L$. This can be checked case by case using Table~\ref{tab:conjugacy-classes}. If $\rho(c_1)$ and $\rho(c_2)$ are complex reflections with respective mirrors $M_1$ and $M_2$, then $M_1$ and $M_2$ must be ultraparallel because $\rho(c_1c_2)$ does not fix any point in $\CH^2$ and does not stabilize any horosphere foliation. Furthermore, $L$ is the common perpendicular to $M_1$ and $M_2$, hence again stabilized by both $\rho(c_1)$ and $\rho(c_2)$. We can draw the same conclusion for $\rho(c_3)$ and $\rho(c_4)$, so that $\Im(\rho)\subset \stab(L)$, contradicting the regularity of $\rho$.   
\end{proof}

Recall from Lemma~\ref{lem:properties-Burau-fixed-points} that Burau fixed points have pseudo-reflection peripheral monodromies and map every non-peripheral simple closed curves to regular elliptic elements. For Gassner representations, interior simple closed curves can be mapped to either screw-parabolic or regular elliptic elements by Lemma~\ref{lem:gassner-fixed-points-interior-scc-are-RE-or-SP}.

\subsection{Rotation numbers and product of pseudo-reflections}
Given two conjugacy classes of pseudo-reflections,
we want to describe all possible conjugacy classes that can be obtained as the product of two elements in these conjugacy classes. To this effect, we introduce rotation numbers as handy morphisms from reducible subgroups of $\PU(2,1)$ to $\U(1)$.
\begin{defn}
    Let $X\in\CH^2$ and fix an identification $\stab(X)=\mathrm{P}(\U(2)\times \U(1))$. The \emph{$X$-rotation number} is
    the morphism $\eta_X\colon\stab(X)\to \U(1)$
    defined by $\eta_X([A,\lambda])=\lambda^{-2}\det A$.
\end{defn}
\begin{defn}
    Let $L\subset\CH^2$ be a $\C$-line and fix an identification $\stab(L)=\mathrm{P}(\U(1)\times \U(1,1))$.
    The \emph{$L$-rotation number} is the morphism
    $\eta_L\colon\stab(L)\rightarrow \U(1)$
    defined by $\eta_L([\lambda,A])=\lambda^{-2}\det(A)$.
\end{defn}
A direct computation gives the following.
\begin{fact}
        The $X$-rotation number of $A\in\stab(X)$ conjugate to $E(\alpha,\beta)$ (in the notation of~\eqref{eq:elliptic-elements}) is $\eta_X(A)=e^{i(\alpha+\beta)}$.
\end{fact}
\begin{fact}
    The $L$-rotation number of $A\in\stab(L)$ can be computed as follows.
    \begin{enumerate}
        \item If $A$ is conjugate to $E(\alpha,\beta)$ and restricts to a rotation of angle $\alpha$ on $L$, then $\eta_L(A)=e^{i(\alpha-2\beta)}$.
        \item If $A$ is a vertical translation, then $\eta_L(A)=1$.
        \item If $A$ is screw parabolic or hyperbolic, and restricts to a rotation of angle $\theta$ around $L$, then $\eta_L(A)=e^{-2i\theta}$.
    \end{enumerate}
\end{fact}
\begin{prop}\label{prop:apx-irreducible-and-pseudo-reflections-imply-rigid}
    Irreducible representations $\rho\colon\pi_1S_{0,4}\to\PU(2,1)$ whose peripheral monodromies are pseudo-reflections must be rigid.
\end{prop}
\begin{proof}[Sketch of proof]
Let $\rho$ be such a representation.
First, since a subgroup of $\PU(2,1)$ generated by two pseudo-reflections always has non-trivial centralizer by Corollary~\ref{cor:intersection-centralizers-pseudo-reflections-non-trivial}, the restriction of $\rho$ to every pair of pants in $S_{0,4}$ is reducible. Furthermore, since $\rho$ is irreducible, $\rho(c_1c_2)$ is either regular elliptic or screw-parabolic by Lemma~\ref{lem:interior-scc-are-regular-elliptic-or-screw-parabolic}. In both cases, by irreducibility of $\rho$, the group $\langle \rho(c_1),\rho(c_2)\rangle$ is a subgroup of unique group of the form $\stab(X)$, $\stab(L)$, or $\stab(H)$ by Lemma~\ref{lem:maximal-subgroups-of-PU(2,1)-with-non-trivial-centralizer}.

\begin{claim}
    If $\varphi$ is a small perturbation of $\rho$ with the same peripheral monodromies, then $\varphi(c_1c_2)$ is conjugate to $\rho(c_1c_2)$.
\end{claim}
\begin{proofclaim}[Sketch of proof]
    Assume first that $\rho(c_1c_2)$ is regular elliptic. Then, when $\varphi$ is close enough to $\rho$, $\varphi(c_1c_2)$ is also regular elliptic. Furthermore, if $\varphi$ is close enough to $\rho$, then $\langle \rho(c_1),\rho(c_2)\rangle$ and $\langle\varphi(c_1),\varphi(c_2)\rangle$ either both stabilize unique points $X_\rho$, $X_\varphi$, or both stabilize unique $\C$-lines $L_\rho$, $L_\varphi$. Denoting the stabilized object by $W\in \{X,L\}$, we obtain a first equation
    \[
    \eta_{W_\varphi}(\varphi(c_1c_2))=\eta_{W_\varphi}(\varphi(c_1))\eta_{W_\varphi}(\varphi(c_2))=\eta_{W_\rho}(\rho(c_1))\eta_{W_\rho}(\rho(c_2)).
    \]
    The same argument with $c_3,c_4$ instead of $c_1,c_2$, and $V$ as stabilized quantity instead of $W$, gives the second equation
    \[
    \eta_{V_\varphi}(\varphi(c_1c_2))=\eta_{V_\varphi}(\varphi(c_3c_4))^{-1}=\eta_{V_\varphi}(\varphi(c_3))^{-1}\eta_{V_\varphi}(\varphi(c_4))^{-1}=\eta_{V_\rho}(\rho(c_3))^{-1}\eta_{V_\rho}(\rho(c_4))^{-1}.
    \]
    Note that $W\neq V$ since $\rho$ is irreducible. The system of equations 
    \[
    \begin{cases}
        \eta_{W_\varphi}(\varphi(c_1c_2)) = \eta_{W_\rho}(\rho(c_1))\eta_{W_\rho}(\rho(c_2))\\
        \eta_{V_\varphi}(\varphi(c_1c_2)) = \eta_{V_\rho}(\rho(c_3))^{-1}\eta_{V_\rho}(\rho(c_4))^{-1}
    \end{cases}
    \]
    determines the conjugacy class of $\varphi(c_1c_2)$ up to a finite set of solutions because $V\neq W$ (one of them being the conjugacy class of $\rho(c_1c_2)$), which thus forces $\varphi(c_1c_2)$ to be conjugate to $\rho(c_1c_2)$.

    The case where $\rho(c_1c_2)$ is screw-parabolic is dealt similarly. Note that we have no analogous rotation number defined on $\stab(H)$, but since screw-parabolic elements form a one-parameter family of conjugacy classes, the $L$-rotation number equation on the complementary pair of pants is sufficient to get the desired conclusion.
\end{proofclaim}
We now know that $\rho(c_1c_2)$ lies in a conjugacy class that is locally independent of $\rho$ as long as peripheral monodromies are preserved. We will consider the case where it is regular elliptic (the screw-parabolic case can be dealt with similarly) and we normalize so that $\rho(c_1c_2)=\varphi(c_1c_2)$. This means that $\langle \rho(c_1), \rho(c_2)\rangle$ and $\langle \rho(c_3), \rho(c_4)\rangle$ are contained in $\stab(X)$, $\stab(L_1)$, or $\stab(L_2)$, where $X$, $L_1$, and $L_2$ are the invariant objects of $\rho(c_1c_2)$. Since $\rho$ is irreducible, the two groups are not contained in the same stabilizer and so there are six possible cases. Assume that $\langle \rho(c_1),\rho(c_2)\rangle\subset \stab(W)$ and $\langle \rho(c_3),\rho(c_4)\rangle\subset \stab(V)$ with $W\neq V\in \{X,L_1,L_2\}$. Note that $\rho(c_1c_2)=\rho(c_3c_4)^{-1}\in\stab(V)\cap\stab(W)$, and that $\langle\varphi(c_1),\varphi(c_2)\rangle\subset W$ and $\langle \varphi(c_3),\varphi(c_4)\rangle \subset V$. Moreover, 
\[
Z(\rho(c_1c_2))=\mathcal{CR}_{L_1}\mathcal{CR}_{L_2}=\mathcal{SR}_X\mathcal{CR}_{L_2}=\mathcal{CR}_{L_1}\mathcal{SR}_X,
\]
and since $V\neq W$, we conclude that $Z(\rho(c_1c_2))=Z(\stab(W))Z(\stab(V))$.

Since $W,V$ are isomorphic to either $\mathrm{P}(\U(2)\times\U(1))$ or $\mathrm{P}(\U(1)\times\U(1,1))$, and pair of pants representations are rigid in $\U(2)$ and $\U(1,1)$ character varieties, the restrictions of $\rho$ and $\varphi$ to $\langle c_1,c_2,(c_1c_2)^{-1}\rangle$ and $\langle c_3,c_4,(c_3c_4)^{-1}\rangle$ are conjugate. The conjugation on both pairs of pants is given by an element that lies in the centralizer of the image of the other pair of pants. This means that $\rho$ and $\varphi$ are globally conjugate.
\end{proof}

\bibliographystyle{amsalpha}
\bibliography{references.bib}
\end{document}